\documentclass{report}

\usepackage{amsfonts,epsfig,amsmath,eufrak, amscd}

\usepackage{xy}
\xyoption{all}
\CompileMatrices

\newtheorem{claim}{Claim}[chapter]
\newtheorem{theorem}[claim]{Theorem}
\newtheorem{definition}[claim]{Definition}
\newtheorem{corollary}[claim]{Corollary}

\newtheorem{construction}{Construction}[chapter]

\newtheorem{proposition}[claim]{Proposition}
\newtheorem{lemma}[claim]{Lemma}
\newtheorem{subclaim}{Subclaim}[claim]

\begin{document}
\pagenumbering{roman}
\begin{center}
{\Large \bf Minimal Triangulations of 3-manifolds}
\vskip 0.2in
\centerline{By}
\vskip 0.1in
\centerline{\large  Alexandre Barchechat}
\centerline{\large B.A. (UC Davis) 1996}
\centerline{\large M.A. (UC Davis) 1998}
\vskip 0.3in
\centerline{\large  DISSERTATION}
\vskip 0.1in
\centerline{\large Submitted in partial satisfaction of the requirements for the
degree of}
\vskip 0.1in
\centerline{\large  DOCTOR OF PHILOSOPHY}
\vskip 0.1in
\centerline{\large in}
\vskip 0.1in
\centerline{\large MATHEMATICS}
\vskip 0.2in
\centerline{\large in the}
\vskip 0.1in
\centerline{\large  OFFICE OF GRADUATE STUDIES}
\vskip 0.1in
\centerline{\large of the}
\vskip 0.1in
\centerline{\large  UNIVERSITY OF CALIFORNIA}
\vskip 0.1in
\centerline{\large  DAVIS}
\end{center}
\vskip 0.1in
{\large Approved:}
\begin{center}
\centerline{\underbar{\hskip 2.5in}}
\vskip 0.15in
\centerline{\underbar{\hskip 2.5in}}
\vskip 0.15in
\centerline{\underbar{\hskip 2.5in}}
\vskip 0.2in
\centerline{\large Committee in Charge}
\vskip 0.2in
\centerline{\large 2003}
\end{center}
\newpage
\large
\tableofcontents


\newpage
\large
{\Large \bf ACKNOWLEDGMENTS} \\

\vspace{2cm}

\hspace{1cm} I would like to thank, first of all, my thesis advisor Joel Hass. He has had faith in me since the beginning of my graduate years, and I could never thank him enough for that. His patience for my ``quick'' questions has been extremely valuable to me. I would also like to thank Bill Thurston for my numerous and enlightening meetings with him. His global view of mathematics has helped me to broaden mine in many different ways. I thank Abby Thompson for introducing me to algebraic topology and low-dimensional topology. She was (and still is!) an excellent teacher. 

\hspace{1cm} I would like to thank the staff on the fifth floor of Kerr Hall and especially Celia Davis for all the help she has given me. I would also like to thank all the graduate students and the professors in the mathematics department for making such a great environment for research. In particular, I thank Sunny Fawcett for helping me with my figures in my dissertation. She has taught me how to insert latex commands in Xfig, and without her, my final draft would not look nearly the way it does now. I would like to thank my parents, my brother, my sister in law and my parents in law for supporting me through these hard years of graduate school. 

\hspace{1cm} Finally, and most importantly, I thank my life partner and my dearest love Nicole Hoover. She is the one that made it all possible and worth it. It is through her that I found my inspiration. I also thank her for putting up with my moods through the toughest times.


\newpage 
\begin{center}
\underline{\bf \Large Abstract} 
\vskip 0.1in
\centerline{Minimal Triangulations}
\vskip 0.1in
\centerline{by}
\vskip 0.1in
\centerline{Alexandre Charles Barchechat}
\vskip 0.1in
\centerline{Doctor of Philosophy in Mathematics}
\vskip 0.1in
\centerline{University of California at Davis}
\vskip 0.3in
\end{center}
\hspace{1cm} The study of three-dimensional manifolds, through normal surface theory, started over seventy years ago by Hellmuth Kneser ~\cite{Kn:gnus}. It was then further developed in the sixties by Wolfgang Haken ~\cite{Hak:gnus}. In this thesis, we use normal surface theory to understand certain properties of minimal triangulations of compact orientable 3-manifolds. We start with a brief review of this well known theory and of the theory of I-bundles over compact surfaces. We then show, in Chapters 3, how to collapse non-trivial normal 2-spheres in a triangulated orientable closed 3-manifold and how to obtain an induced triangulation on the resulting summands. This collapsing process has already been described by William Jaco and Hyam Rubinstein but only for normal 2-spheres with special properties, whereas our process can be applied to an arbitrary normal 2-sphere.

\hspace{1cm} In Chapter 4, we describe our main result: any closed orientable reducible 3-manifold equipped with a minimal triangulation with more than 4 tetrahedra contains a non-trivial normal 2-sphere which intersects at most 2 distinct tetrahedra in quadrilaterals. In Chapter 5, we describe the collapsing process of a non-trivial normal disk in a compact orientable irreducible triangulated 3-manifold with nonempty boundary. This process has been described by W. Jaco and H. Rubinstein for $\partial$-irreducible 3-manifolds only. We also generalize the main result of Chapter 4 to 3-manifolds with nonempty boundary. Chapter 6 is devoted to explaining, in full detail, Casson's Algorithm and its generalization to compact orientable 3-manifolds with non-empty boundary. 

\hspace{1cm} Finally, we use Casson's Algorithm and our main result to show that it takes polynomial time, with respect to the number of tetrahedra of the triangulation, to check if a closed orientable 3-manifold equipped with a minimal triangulation is reducible or not.


\newpage
\pagestyle{myheadings} 
\pagenumbering{arabic}
\markright{  \rm \normalsize CHAPTER 1. \hspace{0.5cm}
 Minimal Triangulations }
\large 
\chapter{Introduction}
\thispagestyle{myheadings}

\section{Background on Normal Surface Theory}
\hspace{1cm} All the 3-manifolds considered in this thesis are in the piecewise linear category, i.e. every 3-manifold will be associated with a triangulation. A \textbf{\textit{triangulation}} of a compact orientable 3-manifold $M$ is a set $\Delta$ of pairwise disjoint tetrahedra, together with a family of homeomorphisms, $\Phi$, where the domain and image of each homeomorphism consist of faces of tetrahedra. The identification space, $\Delta$/$\Phi$, is homeomorphic to $M$. If $\Delta$ consist of n tetrahedra, we write $\mathbf {|M|} = n$. Let $\phi$ be a homeomorphism from a face $F_1$ to a faces $F_2$. By abuse of language, we will describe the identification space $\Delta$/$\phi$ as the  \textbf{\textit{gluing}} of $F_1$ with $F_2$. Here is an example of a 2-tetrahedra triangulation of the lens space L(3,1). 

\bigskip
\begin{figure}[h]
\hspace{2in} \psfig{figure=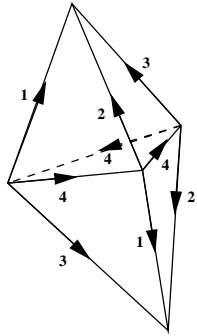, width= 1in} \caption{A triangulation of L(3,1).}
\end{figure}
\bigskip

\hspace{1cm} We are interested in surfaces in $M$ which can be described using the triangulation of $M$. The largest known class of such surfaces is call the class of normal surfaces. Normal surfaces were introduced by Kneser (~\cite {Kn:gnus}) in the late 20's, but they were truly developed only in the early 60's by Haken (see ~\cite{Hak:gnus}). Kneser proved that if $M$ contains an essential 2-sphere, i.e. a sphere which does not cut off a ball, then $M$ contains a normal sphere which is essential. Haken generalized this result by proving that if $M$ contains a surface $F$, then $M$ contains a normal surface $G$ equivalent to $F$. We will define this equivalence later, but the point here is that, up to this equivalence, any embedded surface can be made normal.

\hspace{1cm} We call the  \textbf{\textit{2-skeleton}} ($T^{(2)}$), the  \textbf{\textit{1-skeleton}} ($T^{(1)}$), and the  \textbf{\textit{0-skeleton}} ($T^{(0)}$) of $M$ the respective identification spaces (faces of $\Delta$)/$\Phi$, (edges of $\Delta$)/$\Phi$, and (vertices of $\Delta$)/$\Phi$. Let $R: M \times I \rightarrow M$ be an isotopy of $M$. $R$ is called a  \textbf{\textit{normal isotopy}} if it is invariant in each tetrahedron $\Delta_i$, i.e. $R$ ($\Delta_i$, t)= $\Delta_i$ for all $t$ $\in I$.  \textbf{\textit{Normal surfaces}} are embedded surfaces in $M$ which intersect each tetrahedron in planes in general position with respect to $T^{(2)}$. It is not hard to see that there are 7 different isotopy classes of planes for each tetrahedron: 4 triangles types and 3 quadrilaterals types. We call these triangles and quadrilaterals  \textbf{\textit{elementary disks}}. A normal surface can thus be described as an ordered set of elementary disks in each tetrahedron. Using this definition, we see that a normal surface truly describes a whole class of embedded surfaces which are all normal isotopic. Hence, by a normal surface, we mean its normal isotopy class.

\bigskip
\begin{figure}[h]
\hspace{1in} \psfig{figure=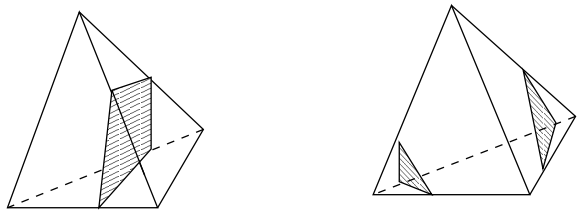, width= 3in} \caption{Some elementary disks}
\end{figure}
\bigskip

\hspace{1cm} For future reference, we call a normal 2-sphere \textbf{\textit{trivial}} if it intersects the tetrahedra in triangles only, and we call it non-trivial otherwise.

\hspace{1cm} The elementary disks have to satisfy some properties. First of all, in order for a normal closed surface to be embedded, it has to satisfy the  \textbf{\textit{quadrilateral property}}: if there is a quadrilateral type in a tetrahedron, then no other types of quadrilateral can exist in that same tetrahedron. Secondly, if $F$ is a face in the triangulation of $M$, then the elementary disks belonging to the two (or the single) tetrahedra, having $F$ in common, must match up. Each matching can be seen as an equation. Indeed, the number of quadrilaterals plus the number of triangles, in a tetrahedron, intersecting a given face in parallel arcs, must equal the number of quadrilaterals plus the number of triangles in the adjacent tetrahedron sharing this face, and whose intersections with the face are parallel arcs of the same type. The matchings are called the \textbf{\textit{matching equations}}. 

\bigskip
\begin{figure}[h]
\hspace{1in} \psfig{figure=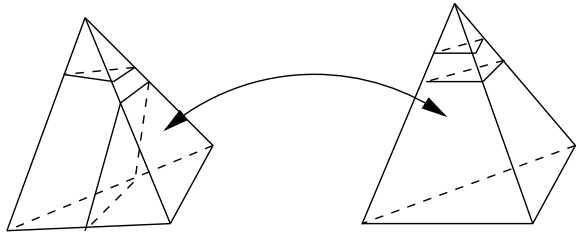, width= 3in} \caption{If two faces are identified, the edges from the elementary disks must match.}
\end{figure}
\bigskip

\hspace{1cm} If $M$ has $t$ tetrahedra, there are exactly $6t$ matchings to be satisfied, 3 for each face of the triangulation. We have noticed above that a normal surface is described by a set of $7t$ elementary disks, 7 for each tetrahedron. So one way to think of a normal surface is to look at it as a non-negative integer valued vector with $7t$ entries satisfying the quadrilateral property and the matching equations. Given a normal surface $F$, we may refer to it as $x_F$, where $x_F = (x_1, x_2,$... $, x_{7t})$.

\vspace{.2in}

\begin{center}
Matching equations\\
$x_{i} + x_{j} = x_{k} + x_{l}$\\
$x_{i} \geq 0$, where $1 \leq i \leq 7t$
\end{center}

\vspace{.2in}

\hspace{1cm} Let $F$ and $G$ be two normal surfaces which intersect each other. Let $x_{F}$ and $x_{G}$ be their corresponding vectors. A natural question would be to ask what kind of surface is represented by the sum of the two vectors. Let us describe first a cut-and-paste operation called \textbf{\textit{regular exchange}}. Suppose that $F$ and $G$ intersect in some tetrahedron with the only requirement that their quadrilateral types are the same. We can perform an operation, described in one of many ways in the figure below, such that the resulting pieces are disjoint elementary disks.

\bigskip
\begin{figure}[h]
\hspace{1in} \psfig{figure=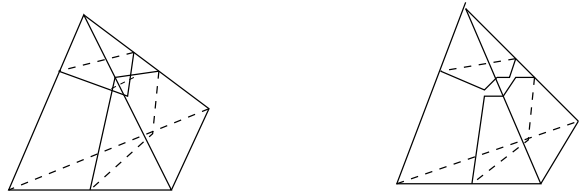, width= 3in} \caption{Before and after a regular exchange.}
\end{figure}
\bigskip

\hspace{1cm} It is not hard to see that the surface obtained in this manner is in the same isotopy class as a normal surface and is represented by the sum of the vectors $x_{F}$ and $x_{G}$. This operation is called the \textbf{\textit{surface addition}}, \textbf{\textit{normal sum}}, or \textbf{\textit{Haken sum}} of $F$ and $G$, and it is well-defined if and only if $F$ and $G$ have the same quadrilateral types in each tetrahedra. Moreover, for each pair of elementary disks, there is only one way of performing a regular exchange. This shows that for each pair of normal surfaces, there is at most one possible Haken sum. See  ~\cite {JT:gnus},~\cite{ Ha:gnus}, ~\cite{JR:gnus}, ~\cite{JS:gnus}, and ~\cite{Mat:gnus} for more details. A useful complexity that is used in normal surface theory is the \textbf{\textit{weight}} of a normal surface $F$, denoted \textbf{\textit{wt(F)}}. It is defined to be the number of intersection points of $F$ with the one-skeleton of the triangulation. One can easily check that the Haken sum preserves the weight and the Euler characteristic, i.e. if $A$ and $B$ are normal surfaces and $A + B = C$, then $wt(A) + wt(B) = wt(C)$ and $\chi(A) + \chi(B) = \chi(C)$. 

\hspace{1cm} Given a normal surface, one may ask if it can be written as the normal sum of two normal surfaces. A normal surface is called \textbf{\textit{fundamental}} if it cannot be written as the sum of two non isotopically parallel surfaces. Obviously, a fundamental surface has to be connected. Haken proved that the set of fundamental surfaces is finite and that they can be found algorithmically. The proof is relatively easy, but the result is strong. It says that any normal surface can be constructed by surface addition from this finite set of normal surfaces. From this, he proved several other results such as: if $M$ is the complement of a non-trivial knot, then it contains an embedded essential disk which is fundamental.

\hspace{1cm} Consider the set of non-negative real solutions to the matching equations and satisfying the quadrilateral property. It is known, through linear programming, that this set forms a cone in $\mathbf{R}^{7t}$. We intersect this cone with the set of solutions to the equation: $\sum_{i=1}^{7t}x_i = 1$. We obtain a convex polyhedron called the \textbf{\textit{projective solution space}} of $M$ with respect to its triangulation $T$. It is denoted \textit{P(M, T)}. For each normal surface $S$, there corresponds a rational vector \=S in \textit{P(M, T)} called the projective class of $S$. Conversely, any rational vector in \textit{P(M, T)} can be multiplied by a rational number to obtain a vector representing a normal surface. 

\hspace{1cm} It can be shown that the vertices of \textit{P(M, T)} have rational entries. Let $v$ be a vertex of \textit{P(M, T)} and let $k$ be the smallest integer such that $k \cdot v$ is an integral solution. We call $k \cdot v$ a \textbf{\textit{vertex solution}}. In particular, an integral solution $F$ is a vertex solution if and only if the integral solutions, $X$ and $Y$, to the equation $n \cdot F = X + Y$ are multiples of $F$. We call $F$ a \textbf{\textit{vertex surface}} if it is connected, 2-sided, and if its representative on \textit{P(M, T)} is a vertex. Note, if $F$ is a vertex surface, then either $F$ is also a vertex solution or it is the double of a vertex solution. Let $F$ be a normal surface. Consider the ray emanating from the origin and passing through $x_{F}$. We intersect this ray with $\sum_{i=1}^{7t} x_i = 1$ to obtain a point in $\mathbf{R}^{7t}$. If $F$ is a vertex surface, then $x_F$ represents a vertex of the convex polyhedron described above. Several results, related to vertex surfaces, have been found (see ~\cite{JO:gnus} and ~\cite{JT:gnus}). For example, if $M$ contains an incompressible surface or an essential sphere or disk, then it must contain one which is a vertex surface. Algorithmically, vertex surfaces are more interesting since there are fewer of them than fundamental surfaces and also they are easier to find.

\hspace{1cm} Let $F$ be an embedded closed surface in $M$. A \textbf{\textit{compression disk}} for $F$ is a disk $D$, embedded in $M$, such that $D \cap F = \partial D$. Doing a compression on $F$ along $D$ means cutting $F$ along $\partial D$ and attaching two parallel copies of $D$ along the two boundary components of the resulting surface.

\bigskip
\begin{figure}[h]
\hspace{1in} \psfig{figure=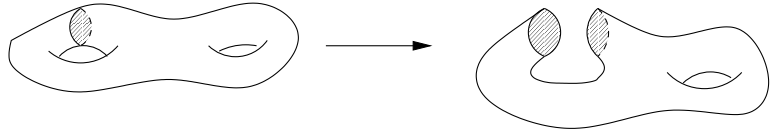, width= 4in} \caption{A compression along an embedded disk.}
\end{figure}
\bigskip

\hspace{1cm} We mention, at the beginning of the introduction that any embedded surface can be represented as a normal surface, up to some equivalence. We now define this equivalence. We say that two normal surfaces are equivalent if one can be obtained from the other one by a series of compressions and normal isotopies. It is not hard to prove (see ~\cite{Ca:gnus}) that any embedded surface $S$, which is in general position with the triangulation, can be made normal up to this equivalence. In particular, if $S$ is incompressible in an irreducible 3-manifold $M$, then it is isotopic to a normal surface. Indeed, if $S$ intersects a face of the 2-skeleton of the triangulation of $M$ in a simple closed curve $\alpha$, then $\alpha$ bounds a disk on $S$. Let's call $D$ the innermost disk on a face of the 2-skeleton whose boundary is $\alpha$. We do a surgery along $D$, and we obtain a surface $S'$. Because $S$ is incompressible, it is homeomorphic to $S'$. Because $M$ is irreducible, $S'$ and $S$ are isotopic.

\vspace{.8in}

\section{Prime Decomposition}

\hspace{1cm} Throughout this thesis, we will make use of the fact that, in a closed orientable 3-manifold, there is a finite collection of independent pairwise disjoint 2-spheres. By a collection of independent 2-spheres, we mean a collection of disjoint embedded 2-spheres such that the closure of its complement does not contain a connected component homeomorphic to a punctured 3-sphere. In 1929, H. Kneser~\cite{Kn:gnus} proved the existence of such a collection of spheres. The proof rely on the existence of triangulations of 3-manifolds: 

{\it \hspace{.2in} \underline{Theorem }~\cite{Kn:gnus}: Let $M$ be a compact 3-manifold. There is an integer $c(M)$ such that if $k>c(M)$ and {\bf S}= $\{ S_1,$ ...,$ S_k \}$ is a set of pairwise disjoint 2-spheres, then a connected component of the complement of {\bf S} is homeomorphic to a punctured 3-sphere.}

\hspace{1cm} W. Haken proved a similar theorem for a collection of pairwise disjoint embedded surfaces. Sometimes, the two theorems together are called the Kneser-Haken finiteness theorem. About 40 years ago, it was already a breakthrough to know that, in a closed 3-manifold, there exist only a finite number of non-parallel disjoint embedded surfaces. Since the development of computers, programs to classify 3-manifolds have led mathematicians to try to find better bounds on the number $k(M)$. For instance, H. Kneser has shown that $k(M)$ does not have to be bigger than $6t + 2dim H_2 (M; \mathbf{Z}_2)$, and Allen Hatcher gives $8t + dimH_1 (M; \mathbf{Z}_2)$ as an upper bound for $k(M)$.

\hspace{1cm} In the early sixties, John Milnor~\cite{Mi:gnus} first proved that this collection of 2-spheres is, in some sense, unique. It is unique in the sense that, if $M \cong M_{1} \#$... $\# M_{k}$ and $M \cong N_{1} \#$... $\# N_{l}$ where the $M_{i}$'s and $N_{j}$'s are prime, then $k=l$ and $M_{i} \cong N_{j}$ after reordering the $M_{i}$'s. Today, several books and notes on low dimensional topology include the proof of the existence and uniqueness of prime decomposition of a closed orientable 3-manifold (e.g. ~\cite{Hat 1:gnus}, ~\cite{Ca:gnus}, ~\cite{Mi:gnus}). This major result is extremely useful in normal surface theory. It says that this unique collection of 2-spheres can be made into a collection of normal 2-spheres. In particular, W. Jaco and J. Tollefson proved that there exists such a collection of normal 2-spheres where all the 2-spheres are vertex surfaces (~\cite{JT:gnus}).

\newpage
\pagestyle{myheadings}
\markright{  \rm \normalsize CHAPTER 2. \hspace{0.5cm}
  Minimal Triangulations }
\chapter{Classification of Homeomorphism Types of the Total Space of I-Bundles over Compact Surfaces}
\thispagestyle{myheadings}

\section{Definition}
  
\hspace{1cm} First of all, let us mention that, for our purpose, this classification will be made up to homeomorphism of $I$-bundles as {\it 3-manifolds}.

\begin{definition}

Let $S$ be a compact surface. We will first assume that $S$ is closed. Throughout the chapter, $I$ is identified with the closed interval $\lbrack -1, 1 \rbrack$. An $I$-bundle $E$ over a surface S consists of a map $p: E \rightarrow S$ called the $projection$ map with the following properties: 

\hspace{1cm} 1) Each point $x$ in $S$ has an open neighborhood $U$ such that $p^{-1}(U)$ is homeomorphic to $ U \times I$. In other words, $E$ has the property that locally, it looks like the product of an open set of $S$ with $I$.  (This is actually the definition of a fibre bundle).

\hspace{1cm} 2) Let $\cup_{i} U_{i}$ be any open cover of $S$ (since $S$ is compact, we can assume that this cover is finite) such that $\phi_{i} : p^{-1}(U_{i}) \rightarrow U_{i} \times I$ is a homeomorphism for each $i$. For $E$ to be an $I$-bundle over $S$, we require the $\phi_{i}$'s to have the following property: when $U_{i}$ and $U_{j}$ have nonempty intersection, the map $\psi_{ij} = \phi_{i} \circ \phi_{j}^{-1} : (U_{i} \cap U_{j}) \times I \rightarrow (U_{i} \cap U_{j}) \times I$ has the form $\psi_{ij}(u, x) = (u, \gamma_{ij}(u)x)$, where $\gamma_{ij}: U_{i} \cap U_{j} \rightarrow \mathbf{ {\it Z_{2}}}$ is a continuous map and $\gamma_{ij}(u)x = x$ or $-x$. 
\end{definition}
\medskip

\hspace{1cm} Intuitively, the maps $\gamma_{ij}$ describe how the pieces $U_{i} \times I$ are glued together (by preserving the fiber structure) in the space $E$. The reason $Im (\gamma_{ij}) \subseteq \mathbf{Z_{2}}$ is because there are only two homeomorphisms (up to isotopy) of the unit interval.

\hspace{1cm} We would like to look at $I$-bundles from a geometric point of view. Hence, let $T$ be a triangulation of $S$ such that for every triangle $t_{i}$, the map $\phi_{i} : p^{-1}(t_{i}) \rightarrow t_{i} \times I$ is a homeomorphism. We can find such a triangulation by triangulating each $U_{i}$ of the open cover with the above properties (here we assume that each $U_{i}$ is closed, but this is only to simplify the notation). Denote the triangles by $U'_{1}$, ..., $U'_{n}$. Now it is easy to see that for any $i, j$, $U'_{i} \cap U'_{j}$ either consists of a single point, an edge, or is empty. Suppose that for each intersecting $U'_{i}$ and $U'_{j}$, the map $\psi_{ij}$ is the identity map. Then $E$ is homeomorphic to the trivial product $S \times I$. 

\hspace{1cm} On the other hand, suppose that $U'_{i} \cap U'_{j} \neq \emptyset$ for some $i, j$ and that $\gamma_{ij} (u, x) = (u, -x)$ for some $u \in (U_{i} \cap U_{j})$ and $x \in I$. Because $\gamma_{ij}$ is continuous and  $G$ is a finite group, we must have $\gamma_{ij} (u, x) = (u, -x)$ for all $u \in (U_{i} \cap U_{j})$ and $x \in I$. Hence, if $U_{i}$ and $U_{j}$ have an edge in common, we must have $\gamma_{ij} (u, x) = (u, -x)$ for all $x$ on this edge, $e$. We will say that $E$ has an \textbf{\textit{inversion}} along $e$ in the triangle $U_{i}$. 

\section{$I$-Bundles over Closed Orientable Surfaces}

 \begin{claim}
The set of edges of inversions of an $I$-bundle over a closed surface $S$ is homeomorphic to a \textbf{\textit{closed}} one-dimensional complex. 
\end{claim}

\hspace{1cm} Here, a closed one-dimensional complex $C$ is a compact subset of $S$ such that the boundary of a regular neighborhood of any point in it is homeomorphic to the finite disjoint union of at least 2 points.

\hspace{1cm} \underline {{\bf Proof :}}
Let $T$ be a triangulation of $S$ with the properties 1 and 2 above. Suppose now that there is a vertex $v$ in $T$ such that only one edge $e$ emanating from $v$ is an inversion. We want to show that there does not exist any neighborhood of $v$ having property 1). Indeed, let $U(v)$ be a regular neighborhood of $v$, homeomorphic to a disk. For convenience, we can think of $U(v)$ as the union of triangles in the second barycentric subdivision of $T$, having $v$ as one of their vertices.

\hspace{1cm} If $U(v)$ had property 1, then $p^{-1}(U(v))$ would be homeomorphic to $U(v) \times I$. In particular, we would have  $p^{-1}(\partial U(v)) \cong (\partial U(v)) \times I \cong S^{1} \times I$ which is an annulus. But, by assumption, $\partial U(v)$ contains exactly one point from e. Hence, $p^{-1}(\partial U(v))$ is homeomorphic to a Mobius band which contradicts property 1.

\bigskip
\begin{figure}[h]
\hspace{1.8in} \psfig{figure=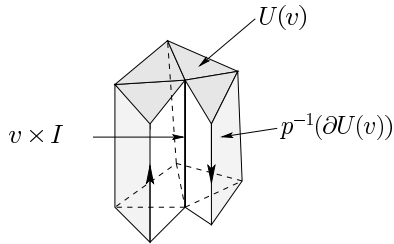, width= 2in} \caption{The boundary of a regular neighborhood of $v \times I$ which is {\it not} homeomorphic to $S^1 \times I$.}
\end{figure}
\vspace{.2in}

\begin{claim} ~\label{CC2}
An $I$-bundle over a closed surface $S$ is determined by a closed \textbf{\textit{embedded}} curve on $S$.
\end{claim}

\hspace{1cm} \underline {{\bf Proof :}} Let $T$ be a triangulation of $S$ with properties 1 and 2. Suppose that the set of inversions of $E$ is not embedded. We want to show that there is another $I$-bundle $E'$, homeomorphic to $E$, through a fiber preserving homeomorphism, with an embedded set of inversions. Suppose that there is a vertex $u$ on $S$ such that there are more than one inversion edges emanating from it. Following the proof of Claim 1, we see that there need to be an even number of those edges. Without loss of generality, suppose that there are four of them: $e_{1}, e_{2}, e_{3}$, and $e_{4}$. One way to visualize the $I$-bundle topologically is as follows: cut $S$ along the inversions, take the trivial product of each resulting piece $V$ with the unit interval $I$, and identify each piece $(\partial V) \times I$ via the map $\phi (v, x) = (v, -x)$. 

\hspace{1cm} Let's take a closer look at what happen locally to this cut-and-paste along the four edges. After cutting along the $e_{i}$'s, suppose we obtain 4 pieces $A, B, C, D$ (this assumption is made only to simplify the notation). As in figure below, we see that $(e_{1} \times I) \cup (e_{2} \times I) \subset (B \times I)$, $(e_{1} \times I) \cup (e_{4} \times I) \subset (A \times I)$, $(e_{3} \times I) \cup (e_{2} \times I) \subset (C \times I)$, and $(e_{4} \times I) \cup (e_{3} \times I) \subset (D \times I)$. Let us call $\phi_{i}$ the inversion map $\phi : e_{i} \times I \rightarrow e_{i} \times I$.

\bigskip
\begin{figure}[h]
\hspace{1.7in} \psfig{figure=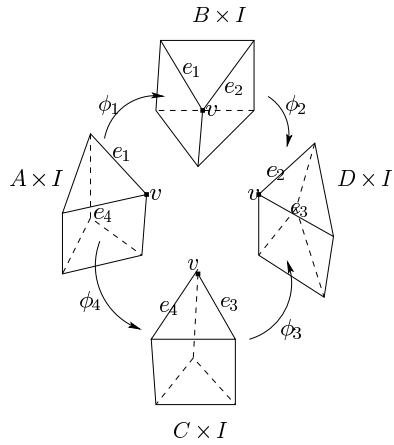, width= 2in} \caption{The inversion maps $\phi_i$ which determine the $I$-bundle.}
\end{figure}
\bigskip

\hspace{1cm} Consider the composition map $\phi_{2} \circ \phi_{1}$ restricted on  $v \times I$:  $\phi_{2} \circ \phi_{1}: v \times I \rightarrow v \times I$. It is the identity map and we call it $\phi_{0}$. Because $\phi_{0}$ is the identity map, it is not part of an inversion and hence, we can glue back $A \times I$ with $D \times I$ along $v \times I$ via $\phi_0$. We now define some new inversions: $\theta_{1} : (e_{1} \times I) \cup (e_{2} \times I) \rightarrow  (e_{1} \times I) \cup (e_{2} \times I) $ by $\theta_{1}(e_{1} \times I) = \phi_{1}(e_{1} \times I)$ and $\theta_{1}(e_{2} \times I) = \phi_{1}(e_{2} \times I)$. Note that $\phi_{1}(v \times I) = \phi_{2}(v \times I)$ so $\theta_{1}$ is well-defined at $u \times I$.

\bigskip
\begin{figure}[h]
\hspace{2.05in} \psfig{figure=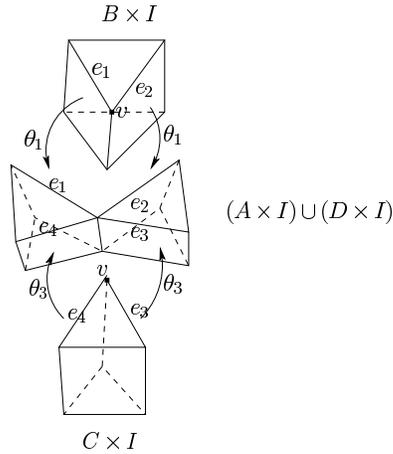, width= 2in} \caption{Reducing the number of points of self-intersection of the curves of inversion.}
\end{figure}
\bigskip

\hspace{1cm} In some ways, we can think of this construction as some kind of regular exchange at the point of intersection. Contrarily to a regular exchange with normal surfaces, there are always two ways to perform this exchange. Indeed, we could have also defined a map $\theta_{2} = \phi_{3} \circ \phi_{2}$. 

\hspace{1cm} After taking the union of $B \times I$, $C \times I$, and $\partial (Nbhd(A \cup D)) \times I$, and redefining the maps $\theta_1$ and $\theta_3$, it is not hard to see that this construction reduces by one the number of intersection points of the curves of inversion. Since these inversions are polygonal curves lying on the one skeleton of some triangulation of $S$, the number of points of self-intersection is finite, and using induction, we complete the proof of the claim.

\bigskip

\begin{claim}
An $I$-bundle over a closed surface $S$ is determined by a closed \textbf{\textit{connected}} embedded curve.
\end{claim}

\hspace{1cm} \underline {{\bf Proof :}} Let $T$ be a triangulation of $S$ with properties 1 and 2. Suppose that the set of inversions consists of two embedded closed curves $\alpha$ and $\beta$. Following the idea of Claim~\ref{CC2}, we can do the reverse construction to obtain one non-embedded closed inversion. We can then do the following exchange:

\bigskip
\begin{figure}[h]
\hspace{1.5in} \psfig{figure=fig2.4.eps, width= 2.7in} \caption{An exchange on the curves of inversion.}
\end{figure}
\vspace{.2in}

\hspace{1cm} It is clear that the number of connected components of inversion is strictly reduced by one. Applying induction on the number of intersection points gives us the desired result.

\bigskip

\begin{claim} ~\label {CC4}
Suppose that the curve of inversion is separating. Then the $I$-bundle is homeomorphic to the trivial product $S \times I$.
\end{claim}

\bigskip

\hspace{1cm} \underline {{\bf Proof :}} Let $T$ be a triangulation of $S$ with properties 1 and 2. Suppose that the inversion $\alpha$ is a separating curve on $S$. As before, cut $S$ along $\alpha$ and take the trivial product of the two resulting pieces ($A$ and $B$) with the interval $I$. 

\hspace{1cm} Instead of applying the inversion on the annuli, we first apply the following fiber-preserving homeomorphism: $\phi : A \times I \rightarrow A \times I$ with $\phi (u, x) = (u, -x)$. Looking at the new pieces, $B \times I$ and $\phi(A \times I)$, our new inversion becomes $\alpha$' = $\alpha \circ \phi$ which is the identity map. Hence, there are no inversions and we have a trivial $I$-bundle.

\bigskip

\hspace{1cm} As we have said earlier, our classification of $I$-bundles will be up to homeomorphism of the total space and not up to $I$-bundle isomorphism. Let us clarify this distinction. Let $E$ be an $I$-bundle with compact base surface $S$. $E$ can be foliated by compact leaves of codimension one. One leaf is homeomorphic to $S$, and all the other leaves are homeomorphic to double covers of $S$. These double covers are completely determined by the curve of inversion on $S$. Conversely, given a double cover of $S$, the curve of inversion is then completely determined by the cover map. Hence, we see that there are at least two ways to describe an $I$-bundle: either by giving a curve of inversion on a surface $S$, or by giving a double cover of that surface. 

\hspace{1cm} Let $E$ and $E'$ be two $I$-bundles over $S$ which are homeomorphic as 3-manifolds. Then we must have $\partial E \cong \partial E'$. As we have just explained, $E$ and $E'$ are completely determined by double covers over $S$, and so we have the commutative diagram:

 $$\xymatrix{
{\partial E} \ar[d]_{i} \ar[r]^{f'} & {\partial E'} \ar[d]_{i'} \\
{E}\ \  \ar[r]^>>>>>>{f}& {E'}.}
$$
\vspace{.2in}


\hspace{1cm} About 80 years ago, it was proved by Nielsen that closed compact surfaces are completely determined by their fundamental groups. Moreover, the fundamental group of an I-bundle is isomorphic to the fundamental group of its base. We conclude that the above diagram is equivalent to the following two:

$$\xymatrix{
{\pi_{1}(\partial E)} \ar[d]_{i_{*}} \ar[r]^{f'_{*}} & {\pi_{1}(\partial E')} \ar[d]_{i'_{*}} \\
{\pi_{1}(E)}\ \  \ar[r]^>>>>>>{f_{*}}& {\pi_{1}(E')}.} $$

\vspace{.2in}

$$\xymatrix{
{\pi_{1}(\text{double cover of $S$})} \ar[d]_{p_{*}} \ar[r]^{f'_{*}} & {\pi_{1}(\text{double cover of $S$})} \ar[d]_{p'_{*}} \\
{\pi_{1}(S)}\ \  \ar[r]^>>>>>>{f_{*}}& {\pi_{1}(S)}.}
$$
\vspace{.2in}

\hspace{1cm} On the other hand, let $E$ and $E'$ be two isomorphic $I$-bundles. Because $E$ and $E'$ are determined by some double cover of $S$, this means that they are isomorphic as covering spaces. In particular, we have the commutative diagram:

$$\xymatrix{
{\pi_{1}(\partial E)} \ar[d]_{p_{*}} \ar[r]^{f'_{*}} & {\pi_{1}(\partial E')} \ar[d]_{p'_{*}} \\
{\pi_{1}(E)}\ \  \ar[r]^>>>>>>{id_{*}}& {\pi_{1}(E')}.} $$

\vspace{.2in}

Let us rewrite this diagram in the following way:

$$\xymatrix{
{\pi_{1}(\partial E)} \ar[d]_{p_{*}} \ar[r]^{f'_{*}} & {\pi_{1}(\partial E')} \ar[d]_{p'_{*}} \\
{\pi_{1}(S)}\ \  \ar[r]^>>>>>>{id_{*}}& {\pi_{1}(S')}.} $$

\vspace{.2in}

\hspace{1cm} We can now clearly see the difference between $I$-bundle isomorphisms and $I$-bundle homeomorphisms of the total spaces. In the former case the homeomorphism of $\pi_{1}(S)$ is the identity, and in the latter case the homeomorphism is arbitrary. Clearly, an $I$-bundle isomorphism is a homeomorphism, and so there are more classes of $I$-bundle isomorphisms than there are classes of $I$-bundle homeomorphisms. 
 Note, since the fundamental group of an $I$-bundle is isomorphic to the fundamental group of the base space (see ~\cite{Sp:gnus}), two $I$-bundles over non-homeomorphic surfaces are non-homeomorphic. On the other hand, if the base surfaces are homeomorphic, two non-isomorphic double cover of $S$ could , potentially, give rise to homeomorphic manifolds. 

\bigskip

\hspace{1cm} We have proved that an $I$-bundle over a closed surface $S$ is determined by a unique (up to isotopy) closed non-separating connected embedded curve on a surface. In other words, a double cover over $S$ can be represented by one closed non-separating connected embedded curve on a surface. We may ask the following question: Are there classes of curves which give rise to homeomorphic $I$-bundles? The answer is yes. 

Let's first suppose that $S$ is orientable. By Claim~\ref{CC4}, we know that there is only the trivial $I$-bundle over $S^{2}$. So let us suppose $S$ has genus $g$ with $g \geq 1$. 

\bigskip

\begin{claim}
There are only \textbf{\textit{two}} non-homeomorphic $I$-bundles over a closed orientable surface.
\end{claim}

\hspace{1cm} \underline {{\bf Proof :}} We already know the existence of the trivial bundle. We claim that the other $I$-bundle is one with an inversion along any connected embedded closed non-separating curve. 

\medskip

\begin{subclaim}
Let $\alpha$ and $\beta$ be two connected embedded closed non-separating curves on $S$. Then there is a homeomorphism of $S$ sending $\alpha$ to $\beta$.
\end{subclaim}

\hspace{1cm} \underline {{\bf Proof of subclaim:}} Let $\chi (S) = -2n$ for some $n$. Cut along $\alpha$. We obtain a twice punctured surface $S'$ with $\chi (S') = -2n$ and boundary components $S'_{1}$ and $S'_{2}$. Similarly we cut along $\beta$ and obtain $S"$ with $\chi (S") = -2n$ and boundary components $S"_{1}$ and $S"_{2}$. By the classification of compact surfaces, $S'$ and $S"$ are homeomorphic via an orientation preserving homeomorphism, say  $h'$. After some isotopies, we can assume that $h'$ leaves the boundary components fixed (He$\lbrack 7 \rbrack$). We now identify $S'_{1}$ and $S'_{2}$ together (resp. $S"_{1}$ and $S"_{2}$), through a orientation reversing homeomorphism f, in order to make the following diagram commute:

$$\xymatrix{
{S'_{1}} \ar[d]_{f} \ar[r]^{h} & {S"_{1}} \ar[d]_{f} \\
{S'_{2}}\ \  \ar[r]^>>>>>>{h}& {S"_{2}}.} $$

\vspace{.2in}

\hspace{1cm} Note, it is possible to make the diagram commute since $h'$ is the identity on the boundary components. This proves the subclaim.

\medskip

\hspace{1cm} Let $E$ and $E'$ be two $I$-bundles over $S$ with curve of inversions $\alpha$ and $\beta$ respectively. Consider the homeomorphism $h$ of $S$ from the subclaim, sending $\alpha$ to $\beta$. We lift $h$ to a homeomorphism $h"$ from $\partial E$ to $\partial E'$ in a canonical way:

\vspace{.2in}

Let $\alpha_{1}$ (resp. $\beta_{1}$) and $\alpha_{2}$ (resp. $\beta_{2}$) be the preimages, under the covering projection, of $\alpha$ (resp. $\beta$). We cut $\partial E$ (resp. $\partial E'$) along these curves. Because $\alpha$ (and $\beta$) is the only curve of inversion, we obtain two disjoint copies of $S'$ (and $S"$). We now use $h'$ to define $h"$. 

\hspace{1cm} Let us summarize what we have shown: we cut one copy of $S$ along $\alpha$ and one copy along $\beta$. We took the trivial product $S' \times I$ and $S" \times I$. We then used $h$ to define a homeomorphism $h"$ from $S' \times$ to $S" \times I$ such that $h"$ is the identity on $(\partial S') \times I$. This gave us a homeomorphism from $E$ to $E'$. Note, if $S$ is closed orientable surface of genus $g$, then $\partial E$ is either two copies of $S$ or a closed orientable surface of genus $2g-1$.

\bigskip

\bigskip

\indent

\section{$I$-Bundles over Closed Non-orientable Surfaces}

\hspace{1cm} Suppose now that $S$ is non-orientable with genus $g$. Here, the genus represents the number of $\mathbf{RP}^{2}$ summands. It seems natural to think that there are two kinds of non-homeomorphic non-trivial $I$-bundles over non-orientable surfaces. One with a two-sided non-separating curve of inversion and another one with a one-sided non-separating curve of inversion (by one-sided we mean a curve $\alpha$ which has a regular neighborhood $N$ such that $N-\alpha$ is connected). We are going to see that this is not the case. 

\medskip

\hspace{1cm} Let $S$ be non-orientable of odd genus ( $g \geq 3$), and $\chi (S) = 1 - 2n$ for some $n$. Let $\alpha$ be a one-sided curve on $S$. If we cut along $\alpha$, the resulting surface $S'$ has one boundary component and the same Euler characteristic. Hence, $S'$ may either be orientable with genus $n$ or non-orientable with genus $2n$. Similarly, if $S$ has even genus, then cutting along a two-sided curve may result in either an orientable or non-orientable surface. We will see that we get non-homeomorphic $I$-bundles depending on the orientability of the resulting surface after cutting along the curve of inversion. 

\bigskip

 First, let us treat two special cases separately:

\hspace{1cm} \underline {{\bf Case 1:}} $S$ is homeomorphic to $\mathbf{RP}^{2}$. In this case, $\pi_{1}( S ) = \mathbf{Z_{2}}$. We conclude that, up to homotopy, only one class of non-separating curves exists, and it is one-sided. Hence, the only nontrivial $I$-bundle over $S$ has boundary a 2-sphere and is homeomorphic to a punctured $\mathbf{RP}^{3}$.

\medskip

\hspace{1cm} \underline {{\bf Case 2:}} $S$ is homeomorphic to a Klein bottle. Here we have $\chi(S) = 0$. Suppose that our $I$-bundle is non-trivial. Let us cut along any non-separating two-sided curve of inversion. We obtain a twice punctured surface $S'$ with $\chi(S') = 0$. The only possibility is a twice punctured 2-sphere, i.e. an annulus. Hence, there can only be one $I$-bundle $E$ , up to homeomorphism, with a two-sided curve of inversion, and $\partial E$ is homeomorphic to a torus. Let us now cut along a one-sided curve. We obtain a once punctured surface $S"$ with $\chi(S") = 0$. The only possibility is a punctured projective plane. So there can only be one $I$-bundle $E'$ with a one-sided curve of inversion. Here we see that  $\partial E'$ is homeomorphic to a Klein bottle. The three $I$-bundles described above all have non-homeomorphic boundaries and hence, they are all non-homeomorphic.

\bigskip

\hspace{1cm} Let $S$ be a closed non-orientable surface with genus greater than 2, i.e. $S= \#_{n} \mathbf{RP}^{2}$ where $n \geq 3$. We will show that there are \textbf{\textit{four}} non-homeomorphic $I$-bundles. Without loss of generality, we will assume that $n$ is odd. Let $\alpha_{1}$ be a two-sided curve, $\alpha_{2}$ one-sided such that S/$\alpha_{2}$ is orientable with genus $(n-1)/2$, and $\alpha_{3}$ one-sided such that S/$\alpha_{3}$ is non-orientable with genus $(n-1)$. Contrarily to case 2, we cannot distinguish these $I$-bundles just by looking at their boundaries because the boundary of $E$ with inversion along $\alpha_{1}$ is homeomorphic to the boundary of $E$ with inversion along $\alpha_{3}$ (we will denote these spaces by ($E$; $\alpha_{1}$) and ($E$; $\alpha_{3}$)).

\bigskip
\begin{figure}[h]
\hspace{1in} \psfig{figure=fig2.11.eps, width= 4in} \caption{Two homeomorphic surfaces whose double covers are non-homeomorphic as covering spaces but homeomorphic as surfaces.}
\end{figure}
\bigskip

\hspace{1cm} Looking at the first homology does not help since the homology of an $I$-bundle is isomorphic to the homology of the base surface (see ~\cite{Sp:gnus}). First of all ($E$; $\alpha_{1}$) and ($E$; $\alpha_{2}$) are non-homeomorphic since they have non-homeomorphic connected boundaries. We want to show that there does not exist homeomorphisms h and $h'$ of $S$ and $\partial E$ making the diagram below commute. 

\hspace{1cm} To differentiate ($E$; $\alpha_{1}$) from ($E$; $\alpha_{3}$), we need to introduce another topological invariant. We have explained earlier that an $I$-bundle over $S$ is determined by a double covering over $S$. Consider the projection map $p$ from the double cover, \~S, onto the base $S$. This map induces a map $p_{*}$ on the fundamental groups, and it is a monomorphism (see Ma$\lbrack 17 \rbrack$). Hence, $\pi_{1}($\~S$)$ can be viewed as a subgroup of $\pi_{1}(S)$. Because $p$ is a double cover, $\pi_{1}($\~S$)$ has index two, and hence it is a normal subgroup. This make the cover $p$ a regular cover. We thus have $\pi_{1}(S)$/$\pi_{1}($\~S$) \cong \mathbf{Z}_{2}$. By the first isomorphism theorem, we conclude that $p$ is determined by a homomorphism $\phi :\pi_{1}(S) \rightarrow \mathbf{Z}_{2}$ with $Ker(\phi) \cong \pi_{1}($\~S$)$. In other words, the double cover is completely determined by a short exact sequence $0 \rightarrow \pi_{1}($\~S$) \rightarrow \pi_{1}(S) \rightarrow \mathbf{Z}_{2} \rightarrow 0$. Because $\mathbf{Z}_{2}$ is an abelian group, we abelianize $\pi_{1}($\~S$)$ and $\pi_{1}(S)$ to obtain a new short exact sequence on the first homology: $0 \rightarrow H_{1}($\~S$) \rightarrow H_{1}(S) \stackrel{\phi}{\rightarrow} \mathbf{Z}_{2} \rightarrow 0$, where $\phi$ is thought of an element of the first cohomology group of $S$. To summarize everything, an $I$-bundle is determined by an element of $H^{1}(S; \mathbf{Z}_{2})$ or by a double covering.  The question now is the following: given two elements of $H^{1}$, corresponding to  ($E$; $\alpha_{1}$) and ($E$; $\alpha_{3}$), how do we check if they give rise to non-homeomorphic $I$-bundles?

\medskip

\hspace{1cm} Suppose ($E$; $\alpha_{1}$) and ($E$; $\alpha_{3}$) are homeomorphic. We have the following commutative diagrams.

$$\xymatrix{
{H_{1}(\partial (E; \alpha_3))} \ar[d]_{p_{*}} \ar[r]^{f'_{*}} & {H_{1}(\partial (E; \alpha_1))} \ar[d]_{p'_{*}} \\
{H_{1}(E; \alpha_3)}\ \  \ar[r]^>>>>>>{f_{*}}& {H_{1}(E; \alpha_1)}.} $$

\vspace{.3in}

\hspace{1cm} Because $f$ is a homeomorphism, it induces isomorphisms, $f_{*}$ and $f^{*}$, on the homology and cohomology groups. We want to reach a contradiction by showing that $f^{*}$ is not an isomorphism on the second cohomology. First, we construct two homomorphisms, $\phi$ and $\phi'$, such that the diagram below is commutative:

$$\xymatrix{
{H_{1}(E; \alpha_3)} \ar[d]_{\phi} \ar[r]^{f_{*}} & {H_{1}(E; \alpha_1)} \ar[d]_{\phi'} \\
{\mathbf{Z}_2}\ \  \ar[r]^>>>>>>>>{f"_{*}}& {\mathbf{Z}_2}.} $$

\vspace{.3in}

\hspace{1cm} From this, we will conclude that $f^{*}(\phi) = \phi'$. We will then show that $\phi \cup \phi = 0$ and $\phi' \cup \phi' \neq 0$, where $\cup$ represents the cup product. Because $f^{*}$ preserves the cup product, we will have $f^{*}(\phi' \cup \phi') = 0$ which contradicts $f^{*}$ being an isomorphism.

\vspace{.1in}

\hspace{1cm} To define $\phi$ and $\phi'$, we only need to define them on basis elements of $H_{1}(S; \mathbf{Z}_{2})$. The idea is to choose the basis elements judiciously. Let $\beta_{1}$ be a curve on $S$, intersecting $\alpha_{1}$ once. We extend $\beta_{1}$ to a basis, $\{ \beta_{1}, \gamma_{1}$, ..., $\gamma_{n} \}$, such that none of the $\gamma_{i}$' intersect $\alpha_{1}$. We explain why such a basis can be found: we cut $S$ along $\alpha_{1}$ and obtain a twice punctured non-orientable surface $S'$. Consider a basis $B$ for the first homology of $S'$. It is clear that the elements of the basis can be chosen disjoint from the boundary components. To obtain $S$ from $S'$, we glue the two boundary components together, which creates a non-trivial element, say $\beta_{1}$, in $H_{1}(S)$. So $B$ together with $\beta_{1}$ constitute a basis for $H_{1}(S)$. Note that $p^{-1}(\gamma_{i})$ represents a closed curve in \~S for all $i$, whereas $p^{-1}(\beta_{1})$ represents a path. Because the diagram in figure 1 is commutative, $f(\gamma_{i})$ must represent, for all $i$, a curve intersecting $\alpha_{3}$ an even number of times. Indeed, if $f(\gamma_{i})$ intersected $\alpha_{3}$ an odd number of times, for some $i$, then $p^{'-1}(f_{*}(\gamma_{i}))$ would represent a path in \~S which would contradict the continuity of $\partial f$. For the same reason, $f(\beta_{1})$ must intersect $\alpha_{3}$ an odd number of times. Consider the two bases of $H_{1}(S; \mathbf{Z}_{2})$: $ B_{1} = \{ \beta_{1}, \gamma_{1}$, ..., $\gamma_{n} \}$ and $ B_{3} = \{ f_{*}(\beta_{1}), f_{*}(\gamma_{1})$, ..., $f_{*}(\gamma_{n}) \}$ (the second set must represent a basis since $f_{*}$ is an isomorphism). 

\hspace{1cm} We now define the maps $\phi$ and $\phi'$ on $B_{1}$ and $B_{3}$ respectively. We define $\phi(\gamma)$ (respectively $\phi'(\gamma)$) to be 1  for all closed curves $\gamma$ that are lifted to a path, and 0 otherwise. Another equivalent way to define the maps would be: $\phi(\gamma) = 1$ for all closed curves $\gamma$ intersecting $\alpha_{1}$ an odd number of times. The fact that $\phi$ and $\phi'$ make the diagram commute comes directly from our choice of basis. The fact that they are homomorphism is also straightforward. We have explained above that $\phi$ and $\phi'$ represent element of the first cohomology of S with $\mathbf{Z}_{2}$ coefficient. Now, because the diagram in figure 2 is commutative, we conclude that $f^{*}(\phi) = \phi'$. First of all, using the universal coefficient theorem and the fact that $H_{0}(S)$ has no torsion, we see that $H_{1}$ is dual to $H^{1}$. From this, we see that $\phi$ (resp. $\phi'$) is dual to $\beta_{1}$ (resp. $f_{*}(\beta_{1})$). We can think of a representative of the cohomology class of $\phi$ as a curve intersecting $\beta_{1}$ an odd number of times and $\gamma_{i}$ an even number of times, and a representative of $\phi'$ as a curve intersecting $f_{*}(\beta_{1})$ an odd number of times and $f_{*}(\gamma_{i})$ an even number of times. 

\hspace{1cm} We now would like to show that $\phi \cup \phi = 0$ and $\phi' \cup \phi' \neq 0$. We give a geometrical interpretation of $\phi$ and $\phi'$ based on Hatcher's notes (see ~\cite{Hat 2:gnus}). Let $S$ be triangulated in the following way:

\bigskip
\begin{figure}[h]
\hspace{1.8in} \psfig{figure=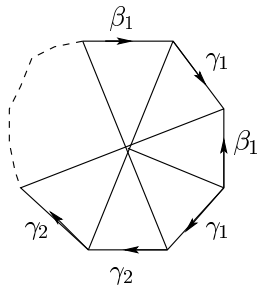, width= 1.3in} \caption{Triangulation of $S$.}
\end{figure}
\bigskip

\hspace{1cm} A representative of $\phi$ (which we will keep calling $\phi$) is a function defined on $H^{1}(S)$, i.e. it assigns to the edges of the triangulation the values 0 or 1. For $\phi$ to be a cocycle, we require that the sum of the values of the edges surrounding each triangle is 0 mod 2. This requirement makes it possible to visualize $\phi$ as a curve which intersects once each triangle with value 1. We have explained above that $\phi$ can be described as a curve intersecting $\beta_{1}$ an odd number of times and $\gamma_{i}$ an even number of times. A natural candidate for $\phi$ is $\alpha_{1}$. Similarly, a natural candidate for $\phi'$ is $\alpha_3$. We show below a choice of representative of $\phi$ and $\phi'$.

\bigskip
\begin{figure}[h]
\hspace{1.8in} \psfig{figure=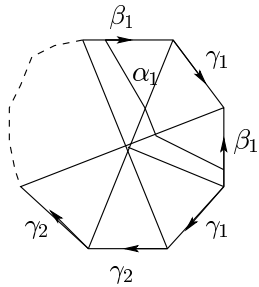, width= 1.3in} \caption{A representative of $\phi$.}
\end{figure}
\bigskip

\hspace{1cm} Now, to compute $\phi \cup \phi$ we first need to orient the edges in the above way. This gives us a consistent way to calculate the cup product in each triangle. Consider two disjoint parallel copies of $\phi$. The cup product $\phi \cup \phi$ is a function which assigns the values 0 or 1 to each triangle. This value is obtained by looking at the two edges, $e_{1}$ and $e_{2}$, having compatible directions, and by multiplying the number of intersection points of the first copy of $\phi$ with $e_{1}$ by the number of intersection points of the second copy with $e_{2}$. The diagrams below show that $\phi \cup \phi$ is the zero function and that $\phi' \cup \phi'$ is non-zero. This proves that $f^{*}$ is not an isomorphism and hence that ($E$; $\alpha_{1}$) is \textbf{\textit{not}} homeomorphic to ($E$; $\alpha_{3}$). 

\bigskip
\begin{figure}[h]
\hspace{1in} \psfig{figure=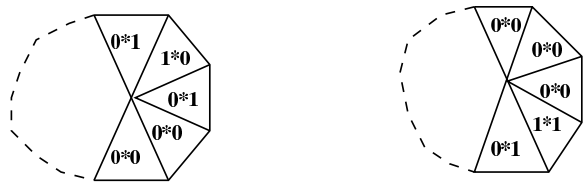, width= 3in} \caption{The cup product of $\phi$ with itself is zero, and the cup product of $\phi'$ with itself is non-zero.}
\end{figure}
\bigskip

\section{$I$-bundles over Compact Surfaces with Nonempty Boundary}

\hspace{1cm} We would like to generalize these results to compact surfaces with nonempty boundary. The proofs being similar to the ones for closed surfaces, we will not give them rigorously.

\hspace{1cm} Let $S$ be an orientable compact surface with $n$ boundary components. Contrary to the case of closed surfaces, an inversion does not have to be connected nor closed. Indeed, there is another kind of inversion: an arc with endpoints in the boundary components of $S$. Following the proofs of Claim 3 and 4, one can show that an $I$-bundle is determined by the disjoint union of one closed connected embedded curve $\alpha$ and two-sided non-separating arcs, $a_{1}$, ..., $a_{k}$ (zero relative homology), such that each connected boundary component of $S$ contains at most one end point from some $a_{i}$. Let $S$ have $n$ connected boundary components. 

\hspace{1cm} \underline{{\bf Case 1:}} The inversion is represented by a single closed connected embedded curve. This is clearly the same as with closed surfaces. If $S$ is orientable, then there are two non-homeomorphic $I$-bundles. If $S$ is non-orientable, there are four. 

\hspace{1cm} \underline{{\bf Case 2:}} The inversion is represented only by two-sided non-separating arcs, $a_{1}$, ..., $a_{k}$. There is exactly one $I$-bundle for each fixed $k$. Moreover, if two inversions are represented by a different number of arcs as above, then they give rise to non-homeomorphic $I$-bundles. This can be seen by counting the number of connected boundary components of the double cover of $S$.

\hspace{1cm} Combining case 1 and case 2, we see that, for an orientable surface $S$, there are $2 \cdot \lfloor n/2 \rfloor$ non-homeomorphic $I$-bundles over $S$. Suppose $S$ has genus $g \geq 1$, $n$ boundary components, and one closed curve and $k$ arcs of inversion. Then the double cover of $S$, \~S, is either disconnected (i.e. two disjoint copies of S) or it has genus $2g-1$ and has $n-k$ boundary components. If $S$ is non-orientable with genus $2g+1$, there are $4 \cdot \lfloor n/2 \rfloor$ non-homeomorphic $I$-bundles over $S$. In this case, \~S is either disconnected, or it has genus $4g$ and $2n-k$ boundary components and is non-orientable, or it has genus $2g$ and $2n-k$ boundary components and is orientable.

\newpage
\pagestyle{myheadings} 
\markright{  \rm \normalsize CHAPTER 3. \hspace{0.5cm} 
  Minimal Triangulations}
\chapter{Normal 2-Spheres in Compact Orientable 3-Manifolds}
\thispagestyle{myheadings}
\section{Statement of the Main Theorem}

\hspace{1cm}  This work is directly inspired by a theorem, stated first by Jaco and Rubinstein, which appeared in ~\cite{JR:gnus}. The theorem is the following: 

\hspace{1cm} (Jaco-Rubinstein) Let $M$ be a closed orientable triangulated 3-manifold. If $M$ contains a non-trivial 2-sphere, then either $M \cong M_{1} \# M_{2}$, with $|M_{1}| + |M_{2}| < |M|$, or $M \cong M_{1} \# r_{1}(S^{1} \times S^{2}) \# r_{2}\mathbf{RP}^{3} \# r_{3}L(3, 1)$, with $|M_{1}| < |M|$.

\hspace{1cm} The proof of this theorem involves a construction called ``crushing'' a 2-sphere, which consists of collapsing two parallel copies of a normal 2-sphere to points. The construction Jaco and Rubinstein describe (see ~\cite{JR:gnus}) is done with normal 2-spheres having certain properties. In fact, essential surgery surfaces, which we will define later on, seem to be an obstruction for this ``crushing''. What we show in this chapter is that the ``crushing'' can be done for any normal 2-sphere, i.e. there are no obstructions to this ``crushing'' process. From an algorithmic point of view, this makes a big difference if the 3-manifold is reducible. Indeed, given a normal 2-sphere, we show that the ``crushing'' process takes time polynomial in the number of tetrahedra. In ~\cite{JR:gnus}, it seems that, the time needed to find a normal 2-sphere with certain properties (no obstructions) may be polynomial if the manifold is irreducible. On the other hand, given a normal 2-sphere in a reducible 3-manifold, it is not clear how long it takes to construct a normal 2-sphere with no obstruction to the crushing process.

\hspace{1cm} Finally, we would like to clarify the fact that Proposition~\ref{Pro1}, Proposition~\ref{Pro2}, and a similar version of Theorem~\ref{MT} have already been proven in ~\cite{JR:gnus}. We give here proofs, obtained independently, which rely on the fact that there are no obstructions to ``crushing" a normal 2-sphere.

\begin{theorem} \label{MT} Let $M$ be a closed orientable triangulated 3-manifold with $t$ tetrahedra ($|M| = t $). Let S be a non-trivial normal 2-sphere. Then M  is homeomorphic to $ M_{1} \# M_{2}$...$\# M_{k}$ $\# r_{1}(S^{1}\times S^{2})$ $\# r_{2}\mathbf {RP}^{3}$ $\#r_{3}L(3,1)$, where $r_{1}$, $r_{2}$, $r_{3}$, $k$ $\geq 0$, $|M_{1}| +...+ |M_{k}| < |M|$ and the $M_{i}$'s are closed orientable triangulated 3-manifolds. 
\end{theorem}

\medskip
\hspace{1cm} If we cut $M$ along $S$, we obtain a 3-manifold $M \backslash S$ with two 2-spheres as boundary (this has been proved in Chapter 2). What we would like to do is to collapse each of the two 2-spheres to a point, and obtain a {\it well-defined} triangulation for the resulting 3-manifold. Note that this is topologically equivalent to gluing a 3-ball on each 2-sphere. A well-defined triangulation of a closed 3-manifold is a union of tetrahedra such that each face of each tetrahedron is identified to a unique face of another tetrahedron (possibly the same tetrahedron) and the boundary of the neighborhood of each vertex is homeomorphic to a 2-sphere. See the introduction for a more rigorous definition of a triangulation.

\vspace{.2in}

\smallskip

\hspace{1cm} After cutting $M$ along $S$, we obtain 7 different types of polyhedra in the cell decomposition of $M \backslash S$: tetrahedra, truncated tetrahedra (which can have 0, 1, 2, 3 or 4 truncations), prisms, truncated prisms, tips, $I \times$ quadrilateral, and $I \times$ triangles. The last two polyhedra will be called \textbf{\textit{I-bundles}}. Note that some types of polyhedra may be combinatorially the same (e.g. a tetrahedron and a tip), but for topological reasons we will consider them different.

\bigskip
\begin{figure}[h]
\hspace{1in} \psfig{figure=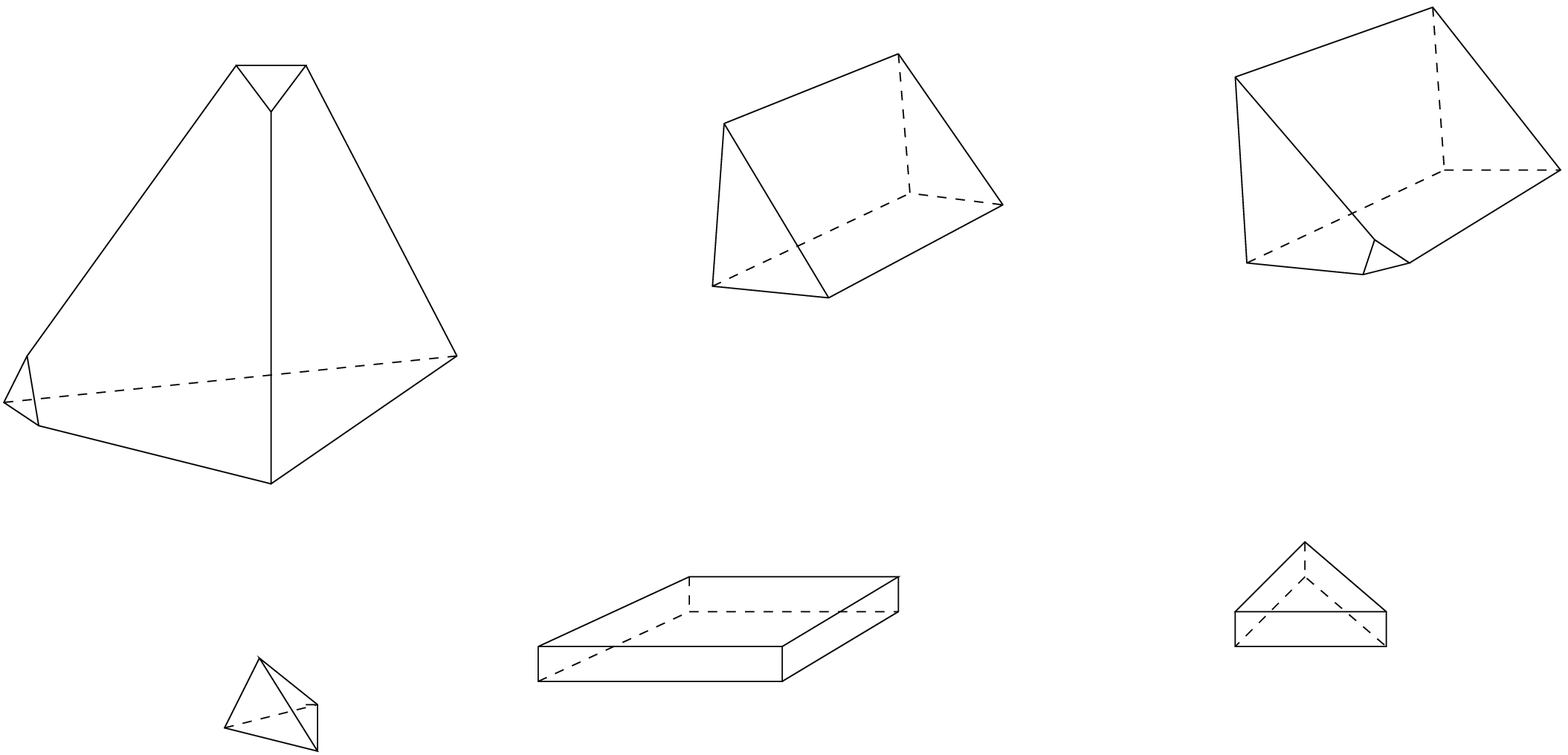, width= 3in} \caption{The seven types of polyhedra in $M \backslash S$}
\end{figure}
\bigskip

\hspace{1cm} Our goal is to cut $M$ along $S$ and change the cell decomposition of $M \backslash S$ to obtain a triangulation for some summands of $M$. The main difficulty in this problem is to change the cell decomposition of $M \backslash S$ by getting rid of all the polyhedra except the truncated tetrahedra. 

\hspace{1cm} \underline {Summary of the algorithm to cut $M$ along a 2-sphere and to retriangulate the}

\underline{ resulting manifold:}

\hspace{1cm} Step 1: We collapse the prisms and truncated prisms one at a time. 
        This collapsing may result in a number of connected summands homeomorphic to $L(3,1)$ or $S^{3}$. As we will see later, the collapsing of these polyhedra may also result in the collapsing of other embedded 2-spheres.

\hspace{1cm} Step 2: We then collapse the tips. This may result in $S^{3}$ summands only.

\hspace{1cm} Step 3: We finally collapse the $I$-bundles. This may result in $S^{3}$ and $\mathbf{RP}^{3}$ summands.

\hspace{1cm} Step 4: We are left with tetrahedra and truncated tetrahedra. We collapse each truncated tetrahedron to a tetrahedron by collapsing each triangle of the original 2-sphere to a point.

\hspace{1cm} Step 5: We count the number of $S^{1} \times S^{2}$ summands.

\medskip

Here are some definitions which will be needed in the proof of the theorem.

\section{Definitions}

\begin{definition}
$M \backslash S$ will denote the resulting manifold after cutting $M$ along $S$, i.e. $M \backslash S= \overline{M-Nbhd(S)}$ where $Nbhd(S)$ denotes a regular neighborhood of $S$.

\hspace{1cm} To describe the process of collapsing the prisms, we foliate them by intervals. To do so, we define a prism $P$ as a quotient space: $(I \times J) \times K / (a, 1, c_{1}) \sim (a, 1, c_{2})$, where $I, J,$ and $K$ are unit intervals. Leaves correspond to $(a, b, K) / (a, 1, c_{1}) \sim (a, 1, c_{2})$.

\bigskip
\begin{figure}[h]
\hspace{1in} \psfig{figure=fig2.eps, width= 3in} \caption{Induced foliation of a prism}
\end{figure}
\bigskip

\hspace{1cm} A similar foliation can be defined for truncated prisms: we define the truncated prism $P$ as a quotient space $(H \times K) / (a, c_{1}) \sim (a, c_{2})$ where $H$ is a hexagon, $K$ is a unit interval, and $a$ is any point of a given edge of $H$. 

\bigskip
\begin{figure}[h]
\hspace{1in} \psfig{figure=fig4.eps, width= 3in} \caption{Collapsing of a truncated prism}
\end{figure}
\bigskip

\textbf{\textit{Collapsing}} a (truncated) prism means taking the quotient space of a prism by identifying each leaf to a point.

\bigskip
\begin{figure}[h]
\hspace{1in} \psfig{figure=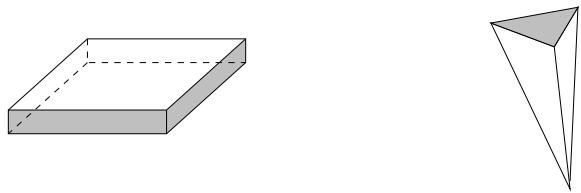, width= 3in} \caption{Sides of an $I$-bundle and the face of a tip}
\end{figure}
\bigskip

\hspace{1cm} We will call the \textbf{\textit{sides}} of an $I$-bundle or a tip, the faces which were originally subsets of the 2-skeleton of $M$.
We will call the \textbf{\textit{face(s)}} of an $I$-bundle or a tip, the face(s) which where originally embedded subsets of $S$.

\hspace{1cm} Similarly, we define the \textbf{\textit{sides}}, the \textbf{\textit{top}} face, and the \textbf{\textit{bottom}} face of a (truncated) prism $P$. Our point here is to give a name for the two hexagonal faces of a truncated prism.

\bigskip
\begin{figure}[h]
\hspace{1in} \psfig{figure=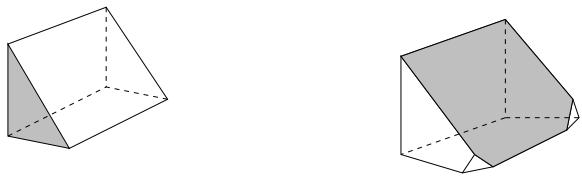, width= 3in} \caption{A side and the top of a (truncated) prism}
\end{figure}
\bigskip

\hspace{1cm} A polyhedra $P_{1}$ is \textbf{\textit{adjacent}} to a (truncated) prism $P_{2}$ if a side of $P_{2}$ is a subset of $P_{1} \cap P_{2}$. There are at most 2 adjacent polyhedra to each prism.

\bigskip
\begin{figure}[h]
\hspace{1.7in} \psfig{figure=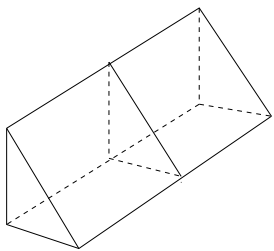, width= 1.4in} \caption{Two adjacent prisms}
\end{figure}
\bigskip

\hspace{1cm} A \textbf{\textit{surgery annulus}} is an annulus $A$ embedded in the 2-skeleton of $M$ with the following properties:

- $\partial A$ = $\beta_{1}$ $\cup$ $\beta_{2}$ , where $\beta_{1}$ and $\beta_{2}$ are composed of normal arcs of the same type.

- $A \cap S = \partial A$.
\smallskip

\hspace{1cm} $A$ \textbf{\textit{surgery Mobius band}} is a Mobius band $B$ embedded in the 2-skeleton of $M$ with the following properties:

- $\partial B$ = $\beta$, where $\beta$ is composed of pairs of normal arcs.

- $B\cap S = \partial B$.

\bigskip
\begin{figure}[h]
\hspace{1in} \psfig{figure=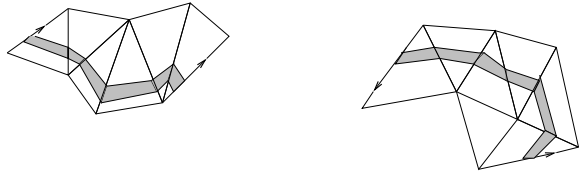, width= 3in} \caption{A surgery annulus and a surgery Mobius band in $T^{(2)}$}
\end{figure}
\bigskip

\hspace{1cm} By \textbf{\textit{surgery disk}} we mean a disk $D$ embedded in the 2-skeleton of $M$ with the following properties:

- $\partial D$ is composed of normal arcs and $D \cap M$ is composed of triangles only. 

- $D \cap S = \partial D$. 
\end{definition}

\bigskip
\begin{figure}[h]
\hspace{2.2in} \psfig{figure=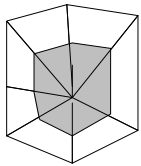, width= .7in} \caption{A surgery disk in $T^{(2)}$}
\end{figure}
\bigskip

\hspace{1cm} Note, a surgery annulus, Mobius band, or disk is always defined with respect to an embedded normal surface. These surgery surfaces will be a major ingredient in the proof of Theorem~\ref{MT}. First, we need the following results. They will be used in the proof of the theorem to describe what happens when surgery surfaces are collapsed.

\begin{lemma} \label{L1}
 Let $M$ be a 3-manifold with boundary a 2-sphere $S$. Collapsing $S$ to a point is topologically equivalent to collapsing it to any simply connected graph. 
\end{lemma}

\hspace{1cm} \underline{{\bf Proof :}} What we mean by two collapsings being topologically equivalent is that the two collapsings induce homeomorphic 3-manifolds. Let $M$ be as above. Consider $N$, a regular neighborhood of a point. $N$ is homeomorphic to a ball. Similarly, a neighborhood $N'$ of a simply connected graph is homeomorphic to a ball. Hence, gluing a ball on $\partial(M \ N)$ is topologically equivalent to gluing a ball on $\partial (M \ N)$.

\begin{lemma} \label{L2}
Let M be a closed orientable 3-manifold. Let S be an embedded 2-sphere and A an embedded annulus such that $A \cap S = \partial A$. Collapsing S to a point and A to an edge is topologically equivalent to cutting M along two disjoint 2-spheres and collapsing each boundary component to a point.
\end{lemma}

\hspace{1cm} \underline{{\bf Proof:}} $A$ has an $I$-bundle structure $S^{1} \times I$, so when we collapse $A$ to an edge we think of the retraction $\phi : S^{1} \times$ I $\rightarrow $ \{pt\} $\times I$. Let $\beta$ and $\beta'$ be the two boundaries of $A$, and let $D$ and $D'$ be the two disjoint disks on $S$ bounded by $\beta$ and $\beta'$, respectively. Let $A'= S-(D \cup D')$. Let $D_{1}$ (resp. $D_{2}$) be a parallel copy of $D$ (resp. $D'$) with boundary $\beta_{1}$ (resp. $\beta_{2}$) parallel to $\beta$ (resp. $\beta'$). Let $N_{1} = Nbhd(\beta) \cap A$ and $N_{2} = Nbd(\beta') \cap A$. Let $A = N_{1} \cup N \cup N_{2}$ as in the figure below. Without loss of generality, we choose $N_1$ such that $\partial N_1 = \beta \cup \beta_1$ and $\partial N_2 = \beta' \cup \beta_2$. 

\bigskip

There are 4 possible different  embeddings of $A$ in $M$.

\bigskip

\hspace{1cm} {\bf Case 1:} $N_{1}$ and $N_{2}$ are on the same side of $S$, i.e. $N_{1}$ and $N_{2}$ intersect the same connected component of $Nbhd(S) - S$, and the embedding of $A$ preserves orientation. We cut along the two 2-spheres $S$ and $D_{1} \cup N \cup D_{2}$. To simplify the notation, we assume that the two 2-spheres are separating. After cutting along $S$ and  $D_{1} \cup N \cup D_{2}$, we obtain three summands $M_1$, $M_2$, and $M_3$ having boundary components $S$,  $D_{1} \cup N \cup D_{2}$, and  $S$ and $D_{1} \cup N \cup D_{2}$, respectively. Our concern is with $M_3$. Note that $M_3$ contains the two annuli $N_1$ and $N_2$. We notice that, collapsing $A$ to an edge is topologically equivalent to collapsing $N_1$, $N_2$, and $N$ to edges. Moreover, because $D_1$ is parallel to $D$, and $D_2$ is parallel to $D'$, we can also cut along these disks and collapse them to points without changing the topology of the manifold. 

\bigskip
\begin{figure}[h]
\hspace{1.2in} \psfig{figure=fig3.4.eps, width= 2.9in} \caption{$M_3$ after collapsing $D_1 \cup N \cup D_2$ and $S$.}
\end{figure}
\vspace{.2in}

\hspace{1cm} Let us recapitulate what we do. We have the 3-manifold $M_3$ with 2 connected boundary components $S$ and  $D_{1} \cup N \cup D_{2}$. We collapse $S$ to a point and the other 2-sphere to an edge. Indeed, $D_1$ and $D_2$ are collapsed to points and $N$ to an edge. When collapsing the two 2-spheres, we need to keep track of the two annuli $N_1$ and $N_2$. The reason we keep track of these annuli is because collapsing the branched surface $S \cup A$ is topologically equivalent to collapsing $S$,  $D_{1} \cup N \cup D_{2}$, $N_1$, and $N_2$. After collapsing $S$ and  $D_{1} \cup N \cup D_{2}$, we see that $N_1$ and $N_2$ are homeomorphic to 2-spheres bounding balls and touching each other in two distinct points. We collapse the two 2-spheres to edges according to the collapsing of $N_1$ and $N_2$. Because they bound balls, the topology of the manifold is unchanged and we proved the lemma. Note, we cannot collapse these two 2-spheres $N_1$ and $N_2$ to points for this would lead to a 3-dimensional object which is not a manifold.

\bigskip

\hspace{1cm} {\bf Case 2:} $N_{1}$ and $N_{2}$ are on the same side of $S$ and the embedding of $A$ reverses orientation. $M \backslash S$ has two boundary components: $D_{1} \cup D'_{1} \cup A'_{2}$ and $D_{2} \cup D'_{2} \cup A_{2} \cup A_{1} \cup A'_{1}$ which are both homeomorphic to 2-spheres (an Euler characteristic  argument can be used to show this). Consider the torus $T = A \cup D_{1} \cup D_{2}$ embedded in $M$. Consider a regular neighborhood $N$ of $T$. Because $M$ is orientable, $N$ must be homeomorphic to a trivial $I$-bundle over $T$. Hence, $\partial N$ consists of two parallel copies of $T$. On the other hand, because $A$ reverses orientation the neighborhood of $T$ is connected which contradicts the triviality of the $I$-bundle. Therefore, this case cannot occur in an orientable 3-manifold. 

\bigskip
\begin{figure}[h]
\hspace{1.5in} \psfig{figure=fig3.5.eps, width= 2in} \caption{A torus such that the boundary of its regular neighborhood is connected.}
\end{figure}
\bigskip
\vspace{.2in}

\hspace{1cm} {\bf Case 3:} $N_{1}$ and $N_{2}$ are on the opposite sides of $S$ and the embedding of $A$ preserves orientation. In this case, one of the summands of the resulting manifold is homeomorphic to $S^1 \times S^2$. Indeed, a regular neighborhood of $S \cup A$ is homeomorphic to a twice punctured $S^1 \times S^2$. We collapse both 2-spheres to edges and, by Lemma~\ref{L1}, we obtain two summands and a third summand homeomorphic to $S^1 \times S^2$. See Claim ~\ref{C1} for more explanation.

\bigskip
\begin{figure}[h]
\hspace{1.5in} \psfig{figure=fig3.3.eps, width= 2in} \caption{A properly embedded annulus $A$ in $M \backslash S$ such that the two components of $\partial A$ lie on different boundary components of $M \backslash S$.}
\end{figure}
\bigskip
\vspace{.2in}

\hspace{1cm} {\bf Case 4:} $N_{1}$ and $N_{2}$ are on the opposite sides of $S$ and the embedding of $A$ reverses orientation. $M \backslash S$ has two boundary components: $D_{2} \cup A_{2} \cup A'_{2} \cup D_{1} $ and $D'_{2} \cup A_{1} \cup A'_{1} \cup D'_{1}$. This case is similar to case 2. Consider the Klein bottle $K = A \bigcup A'$. Since $M$ is orientable, a regular neighborhood of $K$ must be non-trivial (in particular, it must have an orientable boundary). Because $A$ is orientation reversing, the boundary of a regular neighborhood of $K$ is disconnected and hence, it is homeomorphic to two copies of a Klein bottle. This contradicts the orientability of $M$. 

\bigskip
\begin{figure}[h]
\hspace{1.5in} \psfig{figure=fig3.6.eps, width= 2in} \caption{A Klein bottle such that the boundary of its regular neighborhood is connected.}
\end{figure}
\vspace{.2in}

\hspace{1cm} {\bf Case 5:} We include this case for future reference. Suppose now that there are two disjoint 2-spheres $S_{1}$ and $S_{2}$ such that $\beta_{1} \subset S_{1}$ and $\beta_{2} \subset S_{2}$. We cut along the branched surface $S_{1} \cup S_{2} \cup A$, collapse each copy of $S_1$ and $S_2$ to points and $A$ to an edge. We obtain a 3-manifold homeomorphic to one obtained by cutting $M$ along three 2-spheres, and collapsing each of the 2-spheres o points. Indeed, a regular neighborhood of $S_1 \cup S_2 \cup A$ is homeomorphic to a four-times punctured 3-sphere. We cut along each of these 2-spheres which are composed of parallel copies of $D$, $D'$, $D_1$, $D'_1$, and $A$, and we collapse each of them to edges or points depending on if a parallel copy of $A$ is a subset of them.

\bigskip
\begin{figure}[h]
\hspace{1.7in} \psfig{figure=fig3.2.eps, width= 2in} \caption{An embedded annulus whose boundary components belong to disjoint 2-spheres.}
\end{figure}
\vspace{.2in}
Q.E.D.

\begin{lemma} \label{L3}
Let $M$ be a closed orientable 3-manifold. Let $S$ be an embedded 2-sphere and $D$ an embedded disk such that $D \cap S = \partial D$. Cutting along $S \cup D$, collapsing each connected component of $S$ to a point, and each copy of $D$ to an edge is topologically equivalent to cutting along three 2-spheres and collapsing each of them to points.
\end{lemma}

\hspace{1cm} \underline{{\bf Proof:}} First, we foliate $D$ as in the figure below. 

\bigskip
\begin{figure}[h]
\hspace{1in} \psfig{figure=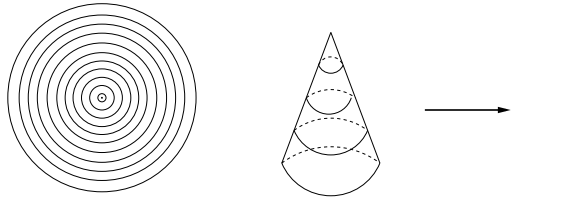, width= 3in} \caption{Collapsing of the foliated disk $D$.}
\end{figure}
\bigskip
\vspace{.2in}

\hspace{1cm} For our purpose, collapsing a disk to an edge means collapsing each circle of the foliation to a point. Cutting along $S \cup D$ means removing a regular neighborhood of it which, in this case, is homeomorphic to a thrice punctured 3-sphere. We collapse each of the three 2-spheres to edges or points, and using Lemma ~\ref{L1} we obtain the result.

\bigskip
\begin{figure}[h]
\hspace{.6in} \psfig{figure=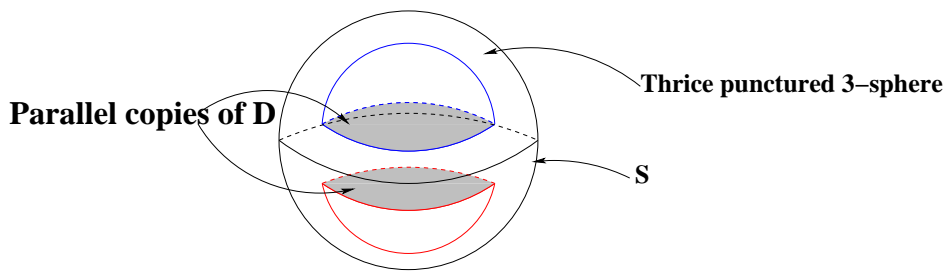, width=5in} \caption{A regular neighborhood of $S \cup D$.}
\end{figure}
\bigskip
\vspace{.2in}

\begin{definition} 

Let $S$ be an embedded 2-sphere and $A$ an embedded annulus such that $A \cap S = \partial A$. We will say that $A$ is \textbf{\textit{inessential}} if: 

\hspace{1cm} a) If A is embedded as in case 1 of Lemma~\ref{L2}, then either the sphere $A \cup D_{1} \cup D_{2}$ or $A \cup A' \cup D_1 \cup D_1$ or $A \cup A' \cup D_2 \cup D_2$  bounds a ball $B^{3}$ such that $int(B^{3})$$\cap S = \emptyset$.

\hspace{1cm} b) If A is embedded as in case 3 of Lemma~\ref{L2}, one of the four 2-spheres $D_1 \cup A \cup D_2$, $D'_1 \cup A \cup D'_2$, $D'_1 \cup A \cup D_2$, or $D_1 \cup A \cup D'_2$  bounds a ball.

\hspace{1cm} Let D be an embedded disk such that $D \cap S = \partial D$. We will say D is \textbf{\textit{inessential}} if one of the two 2-spheres $D \cup D'_{1}$ or $D' \cup D'_{2}$, defined in Lemma~\ref{L3}, bounds a ball $B^{3}$ such that $int(B^{3})$$\cap S = \emptyset$.

\hspace{1cm} Let S be an embedded 2-sphere in M and let B be an embedded Mobius band such that $B \cap S = \partial B = \beta$. Let $S_{1}$ and $S_{2}$ be the boundary of a regular neighborhood of $S \cap D$. Here $S_{2} = A \cap D \cap D'$ for some embedded annulus A and disks D and D'. We will say that B is \textbf{\textit{inessential}} if $D_{1} \cup D_{2} \cup A$ bounds a ball $B^{3}$ such that $int(B^{3}) \cap S = \emptyset$. We will say that a surgery surface is essential if it is not inessential.
\end{definition}
\hspace{1cm} For instance, let $S$ be a vertex 2-sphere. Jaco and Tollefson ~\cite{JT:gnus} proved that if $A$ is a surgery annulus, then the disks $D$ and $D'$ of Lemma~\ref{L2} are normal isotopic. Hence any surgery annulus on a vertex 2-sphere is inessential.

\bigskip
\begin{figure}[h]
\hspace{1.2in} \psfig{figure=fig3.7.eps, width= 3in} \caption{A vertex sphere with inessential surgery annulus.}
\end{figure}
\bigskip
\vspace{.2in}

\section{Proof of the Main Theorem}

\hspace{1cm} We are now ready to prove the theorem. The proof of this theorem is presented in the form of a procedure which, given a normal 2-sphere, finds triangulations of some connected summands of $M$. The drawback is that, given two topologically parallel non-trivial normal 2-spheres with different normal disks, we could end up with different decompositions of $M$. Indeed, given two such spheres, one may have surgery surfaces whereas the other may not. Hence, the geometrical decomposition that we get is $not$ a topological invariant.

\hspace{1cm} Let us start by cutting $M$ along the non-trivial normal 2-sphere $S$ that is given. We obtain a cell decomposition of $M \backslash S$ composed of the 7 types of pieces described earlier.

\hspace{1cm} {\bf Step 1:} First, we collapse one prism at a time. It is clear that if the prism is embedded, then it is homeomorphic to a 3-ball. It is a well known fact that collapsing a 3-ball to a point in a 3-manifold doesn't change its topology. Moreover, collapsing a 3-ball to a point is topologically equivalent to collapsing it to a disk. The proof is similar to Lemma~\ref{L1}: look at the deformation retract $r: D^{2} \times I \rightarrow D^{2} \times \{pt\}$. Hence, if $P$ is embedded, we can collapse it without changing the topology of $M$.

\hspace{1cm} On the other hand, it is not as clear if some, or all, of the leaves are not embedded. We will see that there are 4 different cases to consider. Note that when we collapse a prism, we also collapse its sides. Hence, the adjacent prisms, if any, have one or both of their sides being collapsed. This is not a problem since the $I$-bundle structures (defined in the previous section) of the two prisms coincide on their common side(s).

\hspace{1cm} Let's recapitulate what we do in this first step: we pick an arbitrary (truncated) prism, check into which of the 4 cases it falls, and then collapse it, keeping track of what we are doing to the manifold. We repeat this for all (truncated) prisms. Because there could not be more than two of them in each tetrahedron, the process is finite.

\hspace{1cm} {\bf Case 1 :} All of the leaves of at least one side of the (truncated) prism $P$ are not embedded in $P$, and all the other leaves are. This implies that one, or both, of the sides represents a surgery annulus or a surgery disk depending on whether $P$ is truncated or not.

\bigskip
\begin{figure}[h]
\hspace{1in} \psfig{figure=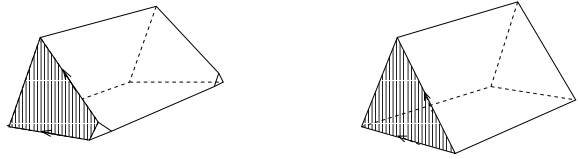, width= 3in} \caption{A surgery annulus and a surgery disk.}
\end{figure}
\bigskip

\hspace{1cm} We first cut along the surgery (annulus) disk and then collapse $P$. Using Lemmas~\ref{L2} and ~\ref{L3}, we see that this is topologically equivalent to cutting and collapsing along other 2-spheres. Note, the union of the two sides of $P$ may represent a surgery (annulus) disk. If this is the case, we cut along both sides and we collapse $P$.

\hspace{1cm} {\bf Case 2:} Exactly one leaf in $P$ from one (or both) of its sides is not embedded. This implies that one (or both) sides represents a surgery Mobius band $B$. This can only happen if the prism is truncated. Indeed, having a regular prism would imply that a vertex is identified to a point in the interior of an edge which clearly violates our definition of a triangulation.  

\bigskip
\begin{figure}[h]
\hspace{1.7in} \psfig{figure=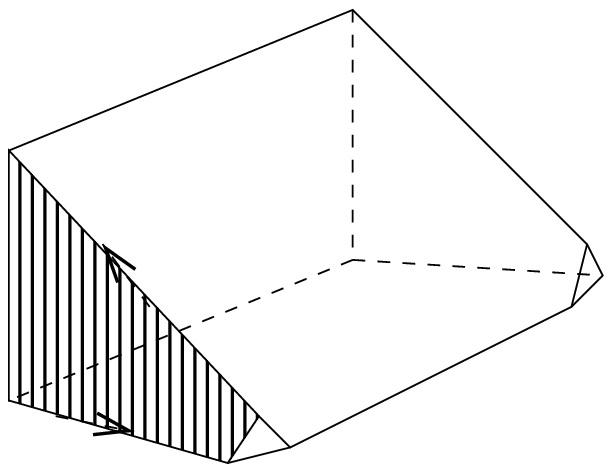, width= 1.5in} \caption{Surgery Mobius band.}
\end{figure}
\bigskip

\hspace{1cm} Let $\beta = \partial B$. Let $E$ be a twisted $I$-bundle over $B$ such that $ \partial (B)$ \~x I $\subset S$. Let $A'$ be the annulus $ \partial (B)$ \~x I and let $\beta_{1}$ and $\beta_{2}$ be its boundary. Since $\beta \subset S$, there exist two embedded disks $D_{1}$ and $D_{2}$ in $S$ such that $\partial D_{1} =\partial D_{2} = \beta$ and $D_{1} \cap D_{2} = \emptyset$. Let $D'_{1}$ (resp. $D'_{2}$) be the disk on $S$ such that $D'_{1} \subset S$ and $\partial D'_{1} = \beta_{1}$ (resp. $D'_{2} \subset S$ and $\partial D'_{2} = \beta_{2}$) and $D'_{1} \cap D'_{2} = \emptyset$. It is known that the boundary of a twisted $I$-bundle over a Mobius band is a torus. Indeed, $ \partial (B)$ \~x I and $B$ \~x $\partial (I)$ are both annuli. Hence, there is an embedded annulus $A$ such that $\partial A = \beta_{1} \cup \beta_{2}$ and $A \cap S = \partial A$. 

\bigskip
\begin{figure}[h]
\hspace{2in} \psfig{figure=fig3.9.eps, width= 3in} \caption{A solid torus $T$ such that $\partial D'_1$ and $\partial D'_2$ represent (2, 1) curves on $\partial T$.}
\end{figure}
\bigskip

\hspace{1cm} Consider the solid torus $T$ with $\partial T = A \cup A'$. The disks $D'_{1}$ and $D'_{2}$ represents 2-handles attached on $\partial T$. Because $B$ is a Mobius band, $\beta$ represents a (2, 1) curve on $\partial T$, and because $\beta_{1}$ and $\beta_{2}$ are parallel to $\beta$ they also represent a (2, 1) curve. Since $\partial D'_{1} = \beta_{1}$ and $\partial D'_{2} = \beta_{2}$, we have a solid torus $T$ on which there are two 2-handles attached on its boundary along two parallel (2, 1) curves.

\hspace{1cm} We conclude that  $T \cup D'_{1} \cup D'_{2}$ is homeomorphic to a twice punctured $\mathbf{RP}^{3}$ ($\cong  L(2, 1)$). The boundary components of a regular neighborhood of it are parallel to $D'_{1} \cup D'_{2} \cup A'$  and $D'_{1} \cup D'_{2} \cup A$. By Lemma~\ref{L2}, we cut along $B$ and we obtain a $L(2, 1)$ summand. 

\bigskip
\begin{figure}[h]
\hspace{2in} \psfig{figure=fig3.8.eps, width= 3in} \caption{A twice punctured $\mathbf{RP}^3$.}
\end{figure}
\bigskip

\hspace{1cm} Hence, if such a $B$ exists, we take an annulus $A$ (as described above) arbitrarily close to the surgery Mobius band, cut along it and collapse it to an edge. Note: the union of the two sides of the truncated prism $P$ may represents a Mobius band. In this case, we cut along both sides and collapse $P$. 

\bigskip
\begin{figure}[h]
\hspace{2in} \psfig{figure=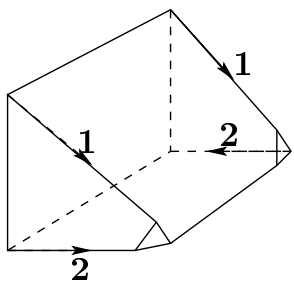, width= 1.5in} \caption{A surgery Mobius band as the union of the sides of a truncated prism.}
\end{figure}
\bigskip

\hspace{1cm} We would like to make the following remark. Because $M$ is orientable, a regular neighborhood of $B$ is homeomorphic to a twisted $I$-bundle over $B$. Such a neighborhood can be foliated by leaves with codimension one, where all the leaves are homeomorphic to an annulus except one which is homeomorphic to $B$. So if we cut along $B$, we obtain a trivial $I$-bundle over one of the leaves. All this to say that instead of constructing an annulus arbitrarily close to $B$, we can simply cut along $B$ to obtain our above result.

\medskip

\hspace{1cm} {\bf Case 3:} None of the leaves of the (truncated) prism are embedded. This implies that the top is identified to the bottom without any twist, and hence, the (truncated) prism is homeomorphic to a 3-ball $B$ in $M$.

\bigskip
\begin{figure}[h]
\hspace{2in} \psfig{figure=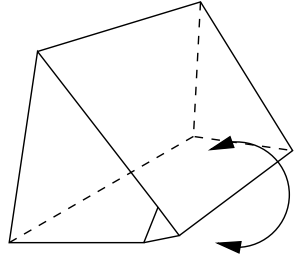, width= 1.5in} \caption{A prism with no embedded leaves.}
\end{figure}
\bigskip

\hspace{1cm} If the prism is regular, then $\partial B$ consists of two surgery disks $D_{1}$ and $D_{2}$, and one annulus $A$ (consisting of one quadrilateral) embedded in $S$. We cut along $D_{1}$ and $D_{2}$, collapse them to edges, and collapse $A$ to a point (since $A \subset S$), i.e. collapse $\partial B$ to an edge with 3 points. By Lemma~\ref{L1}, this is equivalent to collapsing $\partial B$ to a point, and hence we obtain a $S^{3}$ summand.

\bigskip
\begin{figure}[h]
\hspace{2in} \psfig{figure=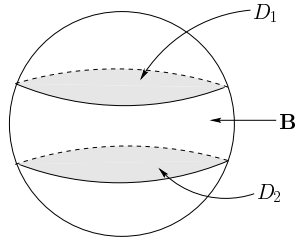, width= 1.5in} \caption{Embedding of a regular prism with no embedded leaves.}
\end{figure}
\bigskip

\hspace{1cm} If the prism is truncated, then $\partial B$ consists of two surgery annuli $A_{1}$ and $A_{2}$, one annulus $A$ and two disks $D_{1}$ and $D_{2}$ embedded in $M$. We cut along $A_{1}$ and $A_{2}$, collapse them to edges, and collapse $A$, $D_{1}$, and $D_{2}$ to points, i.e. collapse $\partial B$ to an edge with 3 points. Again, by Lemma~\ref{L1}, this is equivalent to collapsing $\partial B$ to a point, and hence we obtain a $S^{3}$ summand. 

\bigskip
\begin{figure}[h]
\hspace{2in} \psfig{figure=fig3.12.eps, width= 2.3in} \caption{Embedding of a truncated prism with no embedded leaves.}
\end{figure}
\bigskip

\hspace{1cm} {\bf Case 4:} Exactly one leaf(which does not belong to a face) in the truncated prism is not embedded (the top is identified with the bottom by a 1/3 or 2/3 twist). 

\bigskip
\begin{figure}[h]
\hspace{2in} \psfig{figure=fig3.13.eps, width= 1.6in} \caption{A truncated prism whose boundary is homeomorphic to a solid torus.}
\end{figure}
\bigskip

\hspace{1cm} We calculate the Euler characteristic of the boundary of the truncated prism $P$, and we see that $P$ is homeomorphic to a solid torus $T$. $\partial T$ is decomposed into the union of a surgery annulus $A_{1}$, composed of the two sides of $P$,  and an another annulus $A_{2}$, composed of one quadrilateral and two triangles, embedded in $S$. 

\hspace{1cm} Since $A_{2}$ is lying on $S$, both of its boundary components bound disjoint disk $D_{1}$ and $D_{2}$ such that $int(D_{1}) \cap A_{2} = int(D_{2}) \cap A_{2} =\emptyset$. Moreover, both $\partial D_{1}$ and $\partial D_{2}$ represent a (3, 1) (or a (3,2) curve) curve on $\partial T$. Hence, $T \cup D_{1} \cup D_{2}$ is homeomorphic to a twice punctured $L(3, 1)$ (or a twice punctured $L(3,2$), but $L(3,2)$ and $L(3,1)$ represent homeomorphic manifolds). The boundary of its neighborhood is parallel to $A_{1} \cup D_{1} \cup D_{2}$ and $A_{2} \cup D_{1} \cup D_{2}$. The first component is collapsed to an edge and the second one to a point. The results in a $L(3, 1)$ summand. 

\bigskip
\begin{figure}
\hspace{2in} \psfig{figure=fig3.17.eps, width= 2.3in} \caption{A twice punctured $L(3,1)$.}
\end{figure}
\bigskip

\medskip

\hspace{1cm} {\bf Step 2:} Now that we have eliminated all the (truncated) prisms, let's take a look at the tips. Remark, the side of a tip can only be identified to the side of another tip (it may be the same one) or to the side of a prism. This can be seen from the geometrical structure of each polyhedra. Indeed, the sides of the tip are triangles and the only polyhedra which have triangles as sides are tips and prisms.

\hspace{1cm} Note, if we remove the vertex of a tip, we can give it an $I$-bundle structure.

\bigskip
\begin{figure}
\hspace{2in} \psfig{figure=fig3.15.eps, width= 1in} \caption{The $I$-bundle structure of a tip minus a point.}
\end{figure}
\bigskip

\hspace{1cm} Let's go back before step 1 and look at $M \backslash S$. We take a maximal collection of tips which are glued among themselves along their sides. Their faces represent a connected subsurface of $S$ and hence, they represent a $n$-punctured disk $D_{n}$ with boundary components $S_{1}$, ..., $S_{n}$. 

\hspace{1cm} If we remove the vertex from each tip, we obtain a trivial $I$-bundle over $D_{n}$: $D_{n} \times I$. (We can think of this maximal collection of tips as the cone of $D_{n}$). Each annulus of ($\partial D_{n}) \times I$ corresponds to the union of sides of tips with one vertex removed. Because our collection of tips is maximal and from the remark above, each annulus of ($\partial D_{n}) \times I$ is identified to the side of a prism. In step 1, we collapsed all prisms and in particular we collapsed all of their sides to edges. Hence, we need to collapse each of these annuli to an edge. What is important to notice here is that the collapsing of each side of each prism preserves the $I$-bundle structure of $D_{n} \times I$.

\bigskip
\begin{figure}[h]
\hspace{2in} \psfig{figure=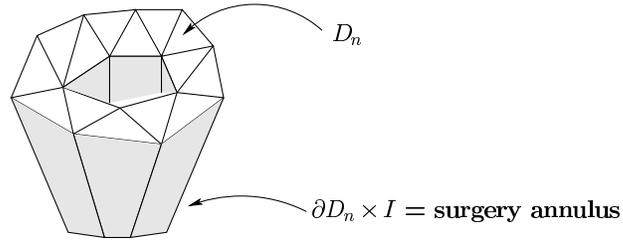, width= 3.2in} \caption{A maximal collection of tips with the vertex removed.}
\end{figure}
\bigskip

\hspace{1cm} When we collapse each annulus to an edge, we also collapse each edge of $\partial D_{n} \times$ \{pt\} to a point. Since each boundary component of $D_{n} \times$ \{pt\} is collapsed to a point, the boundary of our new $I$-bundle is homeomorphic to two disjoint 2-spheres. We get a trivial $I$-bundle over a 2-sphere. Therefore, after the collapsing of all the sides of all the prisms, our maximal collection of tips is homeomorphic to the cone of a 2-sphere, which is in turn homeomorphic to a 3-ball $B$. $\partial B$ comes from the triangulation of $D_{n}$ which was a subset of $S$. Hence, we need to collapse $\partial B$ to a point. This clearly gives us a $S^{3}$ summand.

\hspace{1cm} We now recapitulate what we do in step 2: we look for maximal connected components of tips glued along their sides. Each such component is homeomorphic to the suspension of a 2-sphere, and hence, contributes to an $S^{3}$ summand. 

\bigskip

\hspace{1cm} {\bf Step 3:} We now look at the $I$-bundles. Note, in $M \backslash S$, the side of an $I$-bundle can only be identified to the side of another $I$-bundle (it may be the same one) or to the side of a truncated prism. In fact, the sides of an $I$-bundle are rectangles and the only polyhedra having rectangles for their sides are $I$-bundles and truncated prisms. Let us go back before step 1 and look at the cell decomposition of $M \backslash S$. We take a maximal collection of $I$-bundles which are glued among themselves along their sides. This gives us an $I$-bundle $R$ which may or may not be trivial. $\partial R$, restricted to $\partial I$, consist of faces of the $I$-bundles and hence, it represents a, possibly disconnected, subsurface of $S$. There are two possibilities: $\partial R \vert_{\partial I}$ is either connected or not. By the classification of $I$-bundles (see chapter 2), if it is disconnected, then it represents the trivial $I$-bundle. If it is connected, then we have a twisted $I$-bundle over a projective plane (the base surface could not have higher genus since $\partial R \vert_{\partial I}$ is planar). 

\hspace{1cm} Note, because our collection of $I$-bundles is maximal, each face of ($\partial D_{n}$) \~{x} $I$ (or ($\partial D_{n}) \times  I$) corresponds to a side of an $I$-bundle being adjacent to a truncated prism.

\bigskip
\begin{figure}[h]
\hspace{2in} \psfig{figure=fig3.1.eps, width= 3in} \caption{The collapsing of a surgery annulus on the side of an $I$-bundle preserves the $I$-bundle structure.}
\end{figure}
\vspace{.2in}

\hspace{1cm} As we collapse each side of each truncated prism to an edge (as in case 1, 2, and 4 of step 1) , we also need to collapse each face of ($\partial D_{n}$) \~{x} $I$ (or ($\partial D_{n}) \times I)$ to an edge. Note, the collapsing of each of these faces preserves the $I$-bundle structure. Hence, collapsing each face of ($\partial D_{n}$) \~{x} $I$ (or ($\partial D_{n}) \times I$) corresponds to collapsing each edge of $\partial D_{n}$ (or the disjoint union of two copies of $\partial D_{n}$) to a point. Suppose one of the sides of $R$ is homeomorphic to a Mobius band $B$. This implies that the $I$-bundle $R$ contains a two sided Mobius band which implies that $R$ is non-orientable. Clearly, this cannot happen in an orientable 3-manifold.

\hspace{1cm} After the collapsing of all the faces of $(\partial D_{n})$ \~{x} $I$ (or ($\partial D_{n}$) $\times$ $I$), $D_{n}$ becomes a surface homeomorphic to a 2-sphere, or to two copies of a 2-sphere. If $R$ is trivial (i.e. its boundary component is disconnected), we obtain $S^{2}$ $\times$ $I$ which is homeomorphic to a twice punctured $S^{3}$. If $R$ is non-trivial, i.e. its boundary component is connected, we obtain a punctured $\mathbf{RP}^{3}$. 

\bigskip
\begin{figure}[h]
\hspace{2in} \psfig{figure=fig3.18.eps, width= 3.2in} \caption{Maximal connected collection of $I$-bundles.}
\end{figure}
\vspace{.2in}

\hspace{1cm} We recapitulate what we do in step 3: we look for maximal connected components of $I$-bundles. For each of them, we check if its boundary is connected or not. If it is connected, then the component contributes to an $\mathbf {RP}^{3}$ summand. If it is disconnected, then it contributes to an $S^{3}$ summand. 

\medskip

\hspace{1cm} {\bf Step 4:} We are now left with tetrahedra and truncated tetrahedra only. We collapse each triangle of each truncated tetrahedra to a point, and we obtain the desired triangulation.
\medskip

\hspace{1cm} {\bf Step 5:} In this step, we calculated the number of ($S^{1} \times S^{2}$) connected summands for $M$ found after the collapsing. We first need the following result:

\begin{claim} ~\label{C1}
Let $M$ be an orientable closed 3-manifold and $S$ an embedded non-separating 2-sphere. Then the 3-manifold obtained by cutting $M$ along $S$ and collapsing the two copies of $S$ to points is homeomorphic to $M'$, where $M'$ is an orientable closed 3-manifold such that $M \cong M' \# (S^{1} \times S^{2}$).
\end{claim}

\hspace{1cm} \underline {{\bf Proof :}}
Since $S$ is a non-separating sphere, there is a closed curve $\gamma$ intersecting $S$ in one point only. A regular neighborhood of $S \cup \gamma$ is homeomorphic 
to a punctured $S^{1} \times S^{2}$. So to obtain $M' \# (S^{1} \times S^{2}$), we cut along $\partial (Nbhd(S \cup \gamma))$ and collapse each copy to a point.

\hspace{1cm} How is this construction topologically equivalent to cutting along $S$ only and collapsing its two copies to points? Consider $Nbhd(S \cup \gamma)$, and let $S_{1}$ and $S_{2}$ be the boundary of a regular neighborhood of $S$. The boundary of $Nbhd(S \cup \gamma)$ is a 2-sphere composed of $S_{1}$ minus a disk $D_{1}$, $S_{2}$ minus a disk $D_{2}$, and an annulus $A$ which comes from the boundary of a regular neighborhood of $\gamma$. We now have a 2-sphere $(S_{1} / D_{1}) \cup (S_{2} / D_{2}) \cup A$ such that $D_{1} \cup D_{2} \cup A$ bounds a 3-ball. It is now clear that cutting along this 2-sphere and collapsing it to a point is topologically equivalent to collapsing $S_{1}$ and $S_{2}$ to points. Q.E.D. 

\bigskip
\begin{figure}[h]
\hspace{1.6in} \psfig{figure=fig3.20.eps, width= 3.2in} \caption{The boundary of a regular neighborhood of $S \cup \gamma$.}
\end{figure}
\bigskip

\hspace{1cm} Using this result, we know that each non-separating 2-sphere that we cut along gives us an extra $S^{1} \times S^{2}$ summand. So what we first do is count the number $n$ of summands we have, including the $L(3, 1)$, $\mathbf{RP}^{3}$, and $S^{3}$'s which were found in steps 1, 2, and 3.
We then count the number $m$ of surgery annuli and disks we found in step 1. Suppose that, out of the $m$ surgery surfaces, there are $n_{1}$ disks, $n_{2}$ annuli. If $M$ did not have any connected summands homeomorphic to $S^{1} \times S^{2}$, it would have $n_{1} + n_{2} + 2$ connected summands. Now suppose that, after the collapsing, we obtain $k$ connected summands for $M$. This means that there are $n_1 + n_2 + 2 - k$ non-parallel  disjoint embedded 2-spheres in $M$. By the claim above, $M$ contains $n_1 + n_2 + 2 - k$ connected summands homeomorphic to $S^{1} \times S^{2}$. For each triangulated summand that we obtain, we can apply the Rubinstein-Thompson algorithm to check if it is homeomorphic to $S^{3}$. This can be done by looking for an almost normal 2-sphere in each connected summands (see ~\cite{Tho:gnus} and ~\cite{Ru:gnus}). The resulting decomposition of $M$ will not have any trivial connected summands. We would like to make the following remark: the positive side to this algorithm is that we only need one non-trivial normal 2-sphere to obtain, perhaps, several summands for $M$. Of course, we cannot choose in any way what this number of summands is. It is completely determined by the 2-sphere and the triangulation. Perhaps, one could try to construct a normal 2-sphere with a maximal number of surgery surfaces. In this way, the complete 2-sphere decomposition of $M$ would be obtain somehow rather ``quickly''.

\bigskip

\hspace{1cm} We now need to verify the strict inequality $|M_{1}| +...+ |M_{k}| < |M|$. Because every (truncated) prism is collapsed to a triangle, every tetrahedron in $M$ intersected by $S$ in a quadrilateral will disappear. Since the procedure does not create new tetrahedra, the result follows. In particular, if $S$ intersects n tetrahedra in quadrilaterals, then $|M_{1}| +...+ |M_{k}| + n \leq |M|$. We do not have an equality since some of the summands may be homeomorphic to $S^{3}$. 

\bigskip

\hspace{1cm} To finish the proof of the theorem, we need to check that this construction does give a well-defined triangulation for the connected summands of $M$. First of all, for every $L(3, 1)$, $\mathbf{RP}^{3}$, $S^{3}$, and $S^{1} \times S^{2}$ summand found in the previous steps, no triangulations will be needed (one could always construct a 2-tetrahedra triangulation of $L(3, 1)$, $\mathbf{RP}^{3}$, or $S^{1} \times S^{2}$, and a 1-tetrahedron triangulation of $S^{3}$). But what about the other summands? 

We need to check that:

\hspace{1cm} 1) Each face of each tetrahedron is glued to a unique other one.

\hspace{1cm} 2) The boundary of a regular neighborhood of each vertex is homeomorphic to a 2-sphere.

\hspace{1cm} Let $f$ be a face of a tetrahedron in $M_{1}$. If we think of the collapsing of a (truncated) prism as a continuous map $p$ from a prism to its top (or bottom) face, it makes sense to talk about the preimage $p^{-1}(f)$ of $f$. If $p^{-1}(f)$ represents a face of a tetrahedron in the original triangulation, there is nothing to prove. Suppose that $p^{-1}(f)$ represents a collection of (truncated) prisms in the original triangulation: $P_{1}$, ..., $P_{k}$ where the top of $P_{i}$ is glued to the bottom of $P_{i + 1}$. We assume that the bottom of $P_{1}$ is not glued to the top of $P_{k}$ (the contrary would imply the existence of a $S^{3}$ or $L(3, 1)$ summand). 

\hspace{1cm} Let $t_{1}$ and $t_{2}$ be the two truncated tetrahedra which are glued to the bottom of $P_{1}$ and to the top of $P_{k}$. Since truncated tetrahedra are collapsed to tetrahedra, $t_{1}$ and $t_{2}$ are the only two tetrahedra in $M_{1}$ having $f$ for one of their faces.

\hspace{1cm} By Lemma~\ref{L1} and ~\ref{L2}, we know that each piece $M_{i}$ is a closed orientable 3-manifold. In particular, the boundary of a regular neighborhood of each vertex is homeomorphic to a 2-sphere which proves 2-.

\bigskip

\hspace{1cm} Since the proof of the theorem was given in the form of an algorithm, there is a natural question which one can ask: {\bf What is the running time of the algorithm?}

\hspace{1cm} Step 1 : Because there are $t$ tetrahedra, there are no more than $2t$ (truncated) prisms, and hence, no more than $2t$ collapsings.

\hspace{1cm} Step 2: We only get $S^{3}$ summands so there is nothing to do here.

\hspace{1cm} Step 3: We need to check if the boundary of each $I$-bundle is connected or not. So in each tetrahedron, we collapse each $I$-bundle to an edge, i.e. quad$\times I \rightarrow \{pt\} \times I$. In this way, we are collapsing each leaf of the $I$-bundle to a point. Hence, at the end of this process we end up with one edge for each maximal connected collection of $I$-bundles (see step 3). If one end point of an edge lies in the interior of $M \backslash S$, then the $I$-bundle represents an $\mathbf{RP}^{3}$ summand. If both end points lies in the boundary of $M \backslash S$, then the $I$-bundle represents an $S^{3}$ summand. See chapter 2 to verify these assertions. Unfortunately, the running time of this collapsing process depends on the weight of the 2-sphere (more precisely, on the number of $I$-bundles in each tetrahedron) which may be arbitrarily large. So let us describe another way to calculate the number of non-trivial $I$-bundles. 

\hspace{1cm} First, we calculate $H_{2}(M; \mathbf{Z})$. This can be done in polynomial time in the number of tetrahedra. Indeed, calculating $H_{2}$ is equivalent to transforming a $(2t) \times t$ matrix to its row reduced echelon form, where $t$ is the number of tetrahedra (see ~\cite{Sch:gnus}). The second homology of a compact 3-manifold is a finitely generated abelian group (this can be deduce from the fact that there are finitely many faces in the triangulation). Hence, $H_{2}(M; \mathbf{Z}) \cong \mathbf{Z} \oplus$ ... $\oplus \mathbf{Z}_{p_{1}} \oplus$ ... $\oplus \mathbf{Z} \oplus \mathbf{Z}_{p_{n}}$ by the classification of finitely generated abelian groups (see ~\cite{Ar:gnus}). In particular, every $\mathbf{RP}^{3}$ summand corresponds to a $\mathbf{Z}_{2}$ factor. Let's now calculate the second homology of the resulting manifold after the collapsing. In order to know how many non-trivial $I$-bundles there are, all we need to do is to count how many $\mathbf{Z}_{2}$ factors have disappeared.

\hspace{1cm} Step 4 clearly has a running time linear in $t$ (there are at most $4t$ triangles to collapse to points). Finally, in step 5, we need to count the number of $S^{1} \times S^{2}$ summands.  This was explained earlier: we need to count the number of summands obtained in our decomposition of $M$ and count the number of surgery surfaces. This can be done in linear time.
Therefore, the running time of this collapsing process is polynomial in the number of tetrahedra.

\hspace{1cm} Note, we do not include here the running time to verify if some of the summands are homeomorphic to $S^{3}$. This verification takes exponential time in the number of tetrahedra.

\bigskip

\hspace{1cm} In the decomposition of $M$ obtained through the collapsing process, let's try to calculate the number of vertices in each triangulated summand. Let $S$ be the non-trivial normal 2-sphere and let $A$ be a surgery annulus. When we apply Lemma~\ref{L1}, we see that one or more of the boundary components of $M \backslash (S \cup A)$ will be collapsed to one or more edges. For instance, the 2-sphere $D_{2} \cup A_{2} \cup D'_{2}$ in case 1 of Lemma~\ref{L2} is collapse to an edge with disjoint vertices at its endpoints. Hence, if $M$ has, say $k$ vertices, and $A$ is embedded as in case 1 of Lemma~\ref{L2}, and if we cut along and collapse $S \cup A$, we obtain triangulated summands such that the sum of their number of vertices is $k+4$.

\medskip

\hspace{1cm} Similarly, cutting along a surgery disk and collapsing it to an edge, increases the number of vertices in the triangulation of some summands of $M$. Suppose that S contains $k_{1}$ surgery disks and $k_{2}$ surgery annuli embedded as in case 1 of Lemma~\ref{L2} (we assume this particular embedding in order to simplify our argument but a similar argument would hold for the other case). It is clear that each surgery annulus (or disk) corresponds to adding one extra vertex in one summand and two vertices in another. Hence, suppose we have $k_{1}$ surgery disks, $k_{2}$ surgery annuli, and $k_{3}$ surgery Mobius bands. Then, the sum of the vertices of the summands in the decomposition of $M$ will be $(3k_{1} - 2) + (2k_{2} + 1 ) + (k_{3} + 1)$, where $k_{1} \geq 2$. What needs to be remembered is that, if $M$ has, say 1 vertex, and contains a non-trivial 2-sphere $S$, then after collapsing $S$, it is possible to end up with summands of $M$ with a relatively large number of vertices.

\bigskip

\begin{proposition} \label{Pro1} (~\cite{JR:gnus})
Let $M$ be a closed orientable irreducible 3-manifold triangulated with $t$ tetrahedra and more than one vertex. Then $M$ can be made into a 1-vertex triangulation with strictly less tetrahedra, unless it is homeomorphic to $L(3, 1)$, $\mathbf{RP}^{3}$, or $S^{3}$.
\end{proposition}

\hspace{1cm} \underline{{\bf Proof :}} We would like to mention that this proof is independent of Jaco and Rubinstein's one. Let $T$ be a maximal tree of the 1-skeleton of $M$. Let $S$ be the boundary of a regular neighborhood of $T$. Because $T$ is simply connected, $S$ is homeomorphic to a 2-sphere. We normalize $S$. In the normalization process, we may have to cut $S$ along disks (which are embedded in the 2-skeleton).

\bigskip
\begin{figure}[h]
\hspace{1in} \psfig{figure=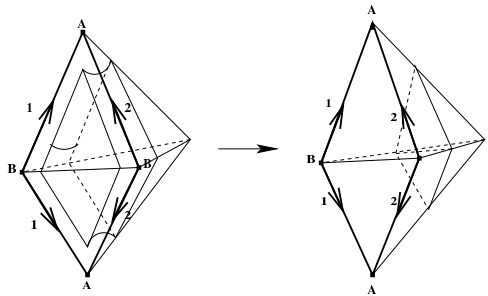, width= 2.5in} \caption{Compressing the boundary of a regular neighborhood of $T^{(1)}$.}
\end{figure}
\vspace{.2in}

\hspace{1cm} Because normalizing an embedded surface which is in general position with the 1-skeleton decreases its weight, the normalization of $S$ does not intersect $T$. In ~\cite{JR:gnus}, $S$ is called a barrier surface. At the end of the normalization process, we obtain $k$ disjoint normal 2-spheres $S_{1}$, ..., $S_{k}$. We discard the ones which have 0-weight and the ones which are composed of triangles only (they bound balls). Because $M$ is irreducible, each of these spheres is separating and hence, it makes sense to talk about the inside and outside of each of them. We now cut along all the $S_{i}$'s and obtain $(k+1)$ pieces. The piece containing $T$ is homeomorphic to a $n$-punctured $S \sp 3$. We will call the outside of $S_{i}$ the side containing the $n$-punctured $S^{3}$. Each other piece is homeomorphic to a 3-manifold with boundary $S_{i}$. Suppose that the inside of $S_{1}$ is not homeomorphic to a 3-ball. Then the inside of all the other spheres must be homeomorphic to a 3-ball because $M$ is irreducible. Hence, there can only be at most one $S_{i}$ such that its inside is not homeomorphic to a 3-ball. We now collapse each of these 2-spheres to points, following the proof of Theorem~\ref{MT}, and we obtain $n$ 3-manifolds $M_{1}$, ..., $M_{n}$ (Note that $n \geq k$ since the normal 2-spheres may have surgery surfaces, giving rises to more than $k$ summands). If one of the $M_{i}$'s is homeomorphic to $L(3,1)$ or $\mathbf{RP}^{3}$, we stop here and give $M$ its 1-vertex 2-tetrahedra triangulation. We see here that if $M$ had originally 2 tetrahedra, we do not decrease the number of tetrahedra.

\hspace{1cm} If none of the summands are homeomorphic to neither $L(3,1)$ nor $\mathbf{RP}^{3}$, we apply the Thompson-Rubinstein algorithm to check which of them is not homeomorphic to $S^{3}$. If all of the $M_{i}$'s are homeomorphic to a 3-sphere, then $M$ is homeomorphic to $S \sp 3$. In this case, there exists a 1-vertex 1-tetrahedron triangulation of $M$. Without loss of generality, suppose $M_{1}$ was found not to be homeomorphic to $S \sp 3$ (in this case $M_{1} \cong M$). If $M_{1}$ has one vertex then we are done. If $M_{1}$ has more than one vertex (recall that in the Theorem~\ref{MT}, the number of vertices of a component $M_{i}$ may have increased and may be bigger than the number of vertices of the original manifold), we apply this construction again to $M_{1}$. This process must terminate after a finite number of steps since $|M_{1}| < |M|$.

\medskip

\hspace{1cm} If $M$ is homeomorphic to either $\mathbf{RP}^{3}$ or $L(3, 1)$, then $M$ could be equipped with the minimal 2-vertex 2-tetrahedron triangulation and in this case we could not get a 1-vertex triangulation with fewer tetrahedra. Later on, we will see that by adding tetrahedra in the triangulation of any closed orientable manifold, we can always make the triangulation a 1-vertex one. The point in this proposition is that the 1-vertex triangulation of $M$ which we obtain has strictly less tetrahedra than the original one.  Q.E.D.

\hspace{1cm} It is not hard to see that any compact orientable irreducible 3-manifold $M$ (except $\mathbf{RP}^{3}$) admits a triangulation in which there are only trivial 2-spheres (these are called 0-efficient triangulations in ~\cite{JR:gnus}): suppose $M$ contains a non-trivial 2-sphere. Cut along it and then collapse it to a point, following the proof of Theorem~\ref{MT}. We may have the following decompositions: $M \cong M_{1}$ with $|M_{1}| < |M|$, $M \cong M_{1} \#$ ...$\# M_{r}$, or $M \cong L(3, 1)$. In the first case we strictly decrease the number of tetrahedra so eventually if we repeat this procedure, we will get rid of all non-trivial 2-spheres. In the second case, we use the Thompson-Rubinstein algorithm to check which summand is not homeomorphic to $S^{3}$. No matter which summand it is, it will have strictly less tetrahedra then $M$. By repeating this procedure we will get rid of all non-trivial 2-spheres. Finally, in the last case, take the 1-vertex 2-tetrahedra triangulation of $L(3, 1)$: it does not contain any non-trivial 2-sphere. The reason we do not consider $\mathbf{RP}^{3}$ is that the two 2-vertex 2-tetrahedra triangulation of it, does contain a non-trivial 2-sphere. 

\bigskip
\begin{figure}[h]
\hspace{1.8in} \psfig{figure=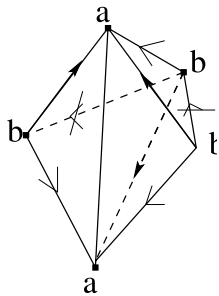, width= 1.1in} \caption{Triangulation of $\mathbf{RP}^3$ with a non-trivial normal 2-sphere isotopic to the boundary of a normal neighborhood of the edge $ab$.}
\end{figure}
\vspace{.2in}

\hspace{1cm} What we hope to do is, given a triangulation of a to compact reducible orientable 3-manifold $M$ with more than one vertex, to construct a 1-vertex triangulation of $M$ with less tetrahedra. On the other hand, we will soon show a construction to change any triangulation into a 1-vertex one by increasing the number of tetrahedra. From an algorithmic point of view, this construction is not very interesting since the number of tetrahedra that are added depends on the number of vertices in the original triangulation.

\begin{proposition} \label{Pro2} (~\cite{JR:gnus})
Let M be a closed orientable irreducible 3-manifold with a minimal triangulation. Then $M$ has one vertex or M is homeomorphic to either $S^{3}$, L(3,1), or $\mathbf{RP}^{3}$.
\end{proposition}

\hspace{1cm} \underline{{\bf Proof :}} Without loss of generality, suppose $M$ has 2 vertices. Following the proof of Proposition~\ref{Pro1}, we normalize the boundary of a regular neighborhood of an embedded edge. We then obtain, say $k$ normal 2-spheres, and we collapse them all. We obtain either $M \approx M_{1}$ with $|M_{1}| < |M|$, $M \cong M_{1} \#$ ...$\# M_{r}$, $M \cong S^{3}$, $M \cong  \mathbf{RP}^{3}$, or $M \cong L(3, 1)$. The first two cases contradict the minimality of $M$. The reason we are not considering $\mathbf{RP}^{3}$, $S^{3}$, or $L(3,1)$ is because $\mathbf{RP}^{3}$ and $L(3,1)$ both have a 2-tetrahedra minimal triangulation with 2 vertices and $S^{3}$ has a 1-tetrahedra minimal triangulation with 2 vertices. Q.E.D.

\newpage
\pagestyle{myheadings} 
\markright{  \rm \normalsize CHAPTER 4. \hspace{0.5cm}
 Minimal Triangulations}
\chapter{Minimal Triangulations}
\thispagestyle{myheadings} 
\section{Connected sums of Triangulated Orientable Closed 3-Manifolds}

\hspace{1cm} We first describe a geometrical construction to obtain a triangulation for the connected sum of two closed orientable triangulated 3-manifolds. Intuitively, one could think of removing the interior of a tetrahedron in each manifold and glue them together along their boundaries. The problem which arises is that the boundary of a tetrahedron may not be homeomorphic to a 2-sphere. One the other hand, one could take the second barycentric subdivision of the triangulations to obtain tetrahedra with embedded boundaries, but this would substantially increase the number of tetrahedra, and, from an algorithmic point of view, this is not interesting.

\smallskip

\begin{construction} \label{cons1}
Let $P$ and $N$ be two triangulated closed orientable 3-manifolds with $t_{1}$ and $t_{2}$ tetrahedra, and $v_{1}$ and $v_{2}$ vertices, respectively. Assume that not both of $v_{1}$ and $v_{2}$ are equal to 1. Then there is a triangulation for $P \# N$ with $t_{1} + t_{2} + 2$ tetrahedra and $v_{1} + v_{2} - 2$ vertices. There are some cases, depending on the triangulations of P and N, where the connected sum will have $t_{1} + t_{2} + 1$ tetrahedra only. If $v_{1} = v_{2} = 1$, then $P \# N$ has $t_{1} + t_{2} + 4$ tetrahedra and one vertex. This triangulation of $P \# N$ does not have to be minimal, even if $P$ and $N$ are minimal.
\end{construction}

\medskip

\hspace{1cm} One reason we are assuming that either $P$ or $N$ has more than one vertex is because this construction is seen as the inverse construction of Theorem~\ref{MT}: if $M$ has one vertex and if $M \cong P \# N$, then $P$ must have at least one vertex (one copy of $S$) and $N$ must have at least 2 vertices (the other copy of $S$ and the original vertex of $M$). Hence, it is somewhat natural to make this assumption.

\hspace{1cm} We first assume that $v_{1}$ and $v_{2}$ are both strictly greater than 1. Let $A$ be a vertex of $P$ and $B$ a vertex of $N$. We remove a normal neighborhood of $A$ and $B$. $P/S_{1}$ and $N/S_{2}$ are now 3-manifolds composed of tetrahedra and truncated tetrahedra with a boundary component being a triangulated 2-sphere, $S_{1}$ or $S_{2}$. Our goal is to change the cell decompositions of $P/S_{1}$ and $N/S_{2}$ in order to glue the two manifolds along their boundary and obtain a well-defined triangulation for their connected sum. To simplify the notation, we will denote $P/S_{1}$ and $N/S_{2}$ by $P'$ and $N'$, respectively. 

\hspace{1cm} Because $P$ has more than one vertex, we can assume, without loss of generality, that $A$ (resp. $B$) was chosen so that $P'$ (resp. $N'$) contains a truncated face as in the figure below:

\bigskip
\begin{figure}[h] \label{fig4.1}
\centering{\psfig{figure=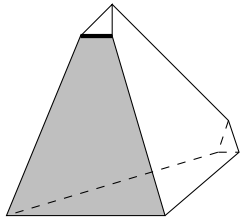, width= 1.5in} \caption{A truncated tetrahedron in $P'$ and $N'$.}}
\end{figure}
\bigskip
\vspace{.2in}

\hspace{1cm} For future reference, we keep track of the above thickened edges on the triangulations of $S_{1}$ and $S_{2}$.

\begin{definition}

We will say that a face of a tetrahedron in a triangulation is a \textbf{\textit{cone}} if two of its edges are identified as follows:

\bigskip
\begin{figure}[h]
\centering{\psfig{figure=fig4.2.eps, width= 1.3in} \caption{A cone in a tetrahedron.}}
\end{figure}
\bigskip
\vspace{.2in}

Let $\tau$ be the triangulation of a surface $S$. The $dual$ $\tau'$ of $\tau$ is a cell decomposition such that :

\hspace{1cm} - There is a 1-1 correspondence between the i-simplices of $\tau$ and the $(2 - i)$-simplices of $\tau'$.

\hspace{1cm} -Every triangle of $\tau$ (resp. every polygon of $\tau'$) contains exactly one vertex of $\tau'$ (resp. $\tau$).

\hspace{1cm} -Every edge of $\tau$ (reps. $\tau'$) intersects exactly one edge of $\tau'$ (resp. $\tau$).

\end{definition}

\bigskip
\begin{figure}[h]
\centering{\psfig{figure=fig4.3.eps, width= 3.8in} \caption{A triangulation of a 2-sphere and its dual.}}
\end{figure}
\bigskip
\vspace{.2in}

\hspace{1cm} It is not hard to check that for any triangulation there is exactly one dual, up to isomorphism of cell decompositions. We note that the union of the vertices and the edges of $\tau'$ is a 3-valent planar graph, which we denote by \textbf{\textit{G}}:

\hspace{1cm} Indeed, the 3-valency comes from the fact that $\tau$ is decomposed into triangles only. Because $G$ is embedded on a 2-sphere, it must miss a point on that 2-sphere. We remove this point from the 2-sphere, and we observe that $G$ is embedded in a punctured 2-sphere which is homeomorphic to the plane. Hence, any graph embedded in $S^{2}$ is planar. Conversely, any planar graph can be embedded in $S^{2}$ (clearly, since $\mathbf{R}^{2} \subset S^{2}$).

\hspace{1cm} Let $T_{1}$ and $T_{2}$ be the triangulations corresponding to the 2-spheres $S_{1}$ and $S_{2}$, respectively.
Let $G_{1}$ and $G_{2}$ be the graphs corresponding to the dual of $T_{1}$ and $T_{2}$, respectively.

\hspace{1cm} We will first describe a construction to transform the graphs $G_{1}$ and $G_{2}$ into isomorphic graphs $G'_1$ and $G'_2$. We will then describe how the triangulations $T_{1}$ and $T_{2}$ are transformed, and finally we will show what happens to the topology of $P'$ and $N'$.

\hspace{1cm} Let $G_{1}$ and $G_{2}$ be given.
We first color the vertices of $G_{1}$ in white and the vertices of $G_{2}$ in black.
Here are the rules to transform $G_{1}$ and $G_{2}$. We can add 4-valent red vertices on old edges, and 3-valent black and white vertices and edges using the following rules:

\hspace{1cm} 1- Black vertices cannot be joined by an edge to white vertices.

\hspace{1cm} 2- Black (resp. white) vertices cannot be added to $G_{2}$ (resp. $G_{1}$).

\hspace{1cm} 3- Let $r$ be a red vertex. Let $e_{1}$ and $e_{2}$ be opposite edges with one common end point $r$ (since r is 4-valent and $G_{1}$ is planar, the notion of opposite edges makes sense). Let $a_{1}$ and $a_{2}$ be the two other end points of $e_{1}$ and $e_{2}$, respectively. If $a_{1}$ is black (resp. white), then $a_{2}$ cannot be white (resp. black).

\hspace{1cm} 4- The resulting graphs $G'_{1}$ and $G'_{2}$ have to be connected and planar.

\bigskip
\begin{figure}[h]
\centering{\psfig{figure=fig4.4.eps, width= 3in} \caption{Example of a 3-valent graph $G$ and another graph $G'$ obtained by adding edges and vertices according to the above rules.}}
\end{figure}
\vspace{.2in}

\hspace{1cm} How do we transform $G_{1}$ and $G_{2}$ to obtain isomorphic 3-valent planar colored graphs? Let $S_{1}$ and $S_{2}$ be given. To the thickened edge in Figure~\ref{fig4.1} corresponds a unique edge in $G_{1}$ and in $G_{2}$, $e_{1}$ and $e_{2}$ respectively. We insert two red vertices on each of the edges: $r_{1}$ and $r'_{1}$ on $G_{1}$ and $r_{2}$ and $r'_{2}$ on $G_{2}$. We then connect $r_{i}$ and $r'_{i}$ by an edge. Consider now a copy of $G_{2}$ with the edge $e_{2}$ removed. We call this new graph $G_{2, e_{2}}$. Color all the vertices of $G_{2, e_{2}}$ in black. Draw an edge emanating from $r_{1}$ and one emanating from $r'_{1}$. On these two edges, draw the graph corresponding to $G_{2, e_{2}}$. The same procedure can be done by adding $G_{1, e_1}$ to $G_{2}$. It is clear from the picture below that $G'_{1}$ and $G'_{2}$ are isomorphic graphs. 

\bigskip
\begin{figure}[h]
\centering{\psfig{figure=fig4.5.eps, width= 3.3in} \caption{Isomorphism between $G'_1$ and $G'_2$.}}
\end{figure}
\vspace{.2in}

\hspace{1cm} What happens to the triangulations $S_{1}$ and $S_{2}$? Adding 4-valent vertices on $G_1$ can be thought as adding quadrilaterals on $S_1$.  Adding 3-valent black (white) vertices is thought as adding black (white) triangles.

\hspace{1cm} Note, because $G'_{1}$ and $G'_{2}$ are planar and connected, they represent the cell decompositions of 2-spheres. This implies that $S'_{1}$ and $S'_{2}$ also represent the cell decompositions of 2-spheres. From the definition of a dual triangulation, each 3-valent (resp. 4-valent) vertex added on $G_{1}$ corresponds to a triangle (resp. quadrilateral) added on $S_{1}$. Each edge added on $G_{1}$ corresponds to the gluing of two polygons along one of their edges on $S_{1}$. Hence, we changed the triangulation of $S_{1}$ (and $S_2$) into a cell decomposition of a 2-sphere with triangles and quadrilaterals. 

\bigskip
\begin{figure}[h]
\centering{\psfig{figure=fig4.6.eps, width= 3.3in} \caption{View of the new triangulation of $\partial P'$.}}
\end{figure}
\vspace{.2in}

\hspace{1cm} What about the cell decompositions of $P'$ and $N'$? We saw that each 4-valent vertices added on $G_{1}$ corresponds to adding a quadrilateral on $S_{1}$. Here, adding a quadrilateral on $S_{1}$ is seen as inserting a prism in $P'$. See Figure below. Similarly, adding a triangle on $S_{1}$ corresponds to inserting the tip of a tetrahedron in $P'$. We would like to explain why the topology of $P'$ is unchanged at the end of this cnstruction.

\hspace{1cm} As described above, the black vertices that we added on $G_{1}$, which correspond to tips of tetrahedra inserted in $P'$, represent a graph isomorphic to $G_{2, e_{2}}$. Since $G_{2}$ is the dual of the triangulation of a 2-sphere, $G_{2, e_{2}}$ is the dual of the triangulation of a disk. Hence, the union of all the tips inserted in $P'$ is homeomorphic to the cone of a disk, which is homeomorphic to a 3-ball. Hence, inserting the quadrilaterals and tips of tetrahedra can be seen as thickening a face of one of the truncated tetrahedra. This does not change the topology of $P'$.

\bigskip
\begin{figure}[h] ~\label{fig4.7}
\centering{\psfig{figure=fig4.7.eps, width= 2.8in} \caption{The topology of $P$ is unchanged.}}
\end{figure}
\vspace{.2in}

\hspace{1cm} Note, because $G'_{1}$ is isomorphic to $G'_{2}$, $S'_{1}$ is isomorphic to $S'_{2}$, so there exists a cell-preserving homeomorphism $\phi$ sending $S'_{1}$ to $S'_{2}$. From the construction above, every white triangle in $S'_{1}$ corresponds to a truncated tetrahedron, every white triangle in $S'_{2}$ corresponds to the tip of a tetrahedra, and every quadrilateral in $S'_{1}$ corresponds to a quadrilateral in $S'_{2}$. Hence, we apply $\phi$, and we get a well-defined triangulation for $P\# N$. 

\bigskip
\begin{figure}[h]
\centering{\psfig{figure=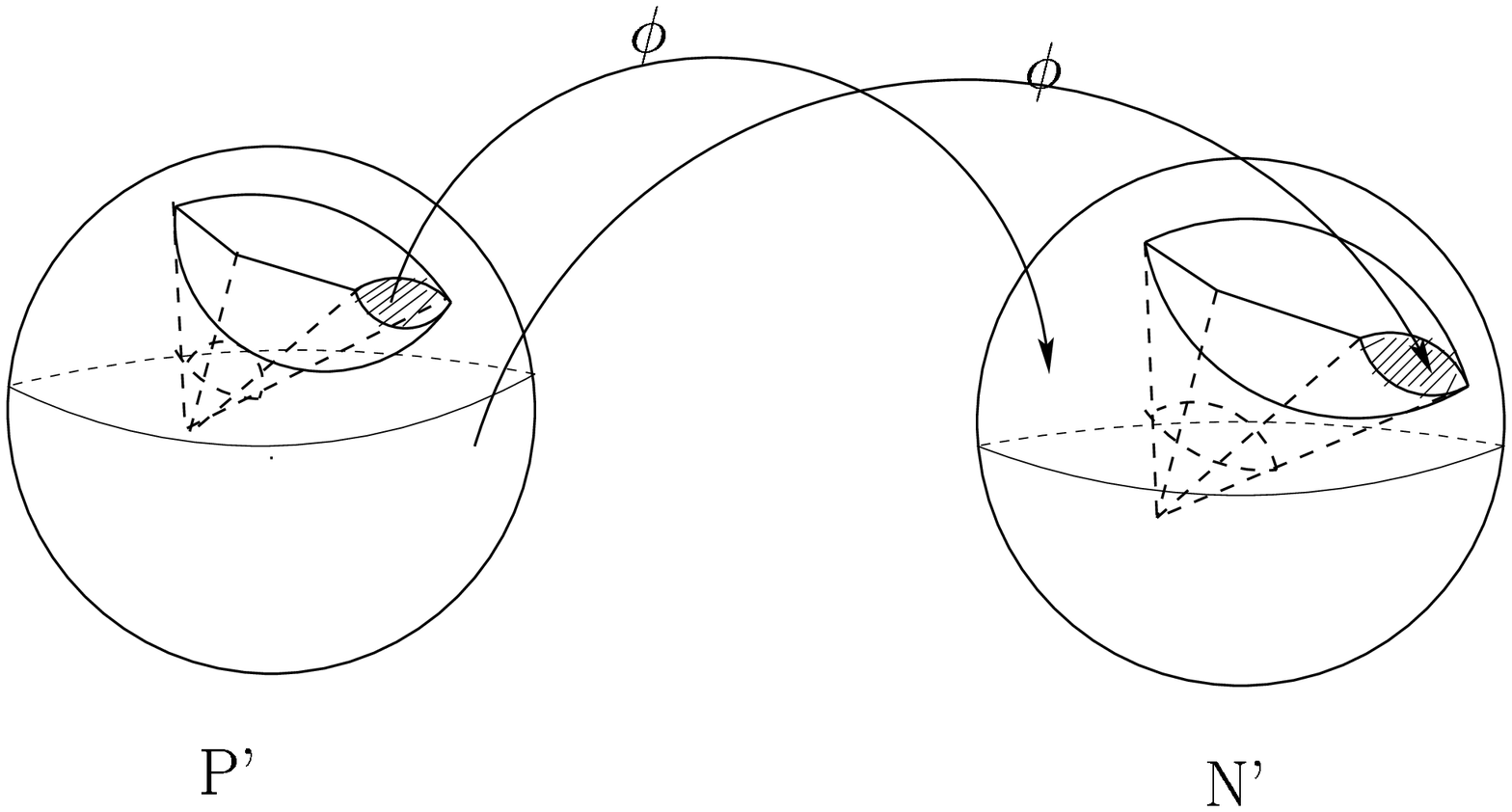, width= 2.8in} \caption{The connected sum of $P$ and $N$.}}
\end{figure}
\vspace{.2in}

\hspace{1cm} Suppose now that $N$ has more than one vertex and that $P$ has exactly one vertex. Because $P'$ does not contain any vertices from $P$, only truncated prisms and truncated tips can be inserted. Let $G_{1}$ and $G_{2}$ be the induced graphs of $P$ and $N$ respectively. Consider $G_{1}$. As before, we add two 4-valent vertices which correspond to inserting two truncated prisms. Because the prisms are truncated, each time we add a 4-valent vertex on $G_{1}$ we need to add two pair of black vertices too.

\hspace{1cm} Adding black vertices on $G_{1}$ goes against our rules defined earlier, but this is not a problem since we are not given a choice on how and where to place these two pairs. Moreover, when we add white vertices on $G_{1}$ corresponding to a graph isomorphic to $G_{2, e_{2}}$, these whites vertices correspond to truncated tips of tetrahedra and so, each time we add a white vertex we need to add a black vertex. Hence, after adding all the necessary vertices on $G_{1}$, we obtain a graph as below:

\bigskip
\begin{figure}[h]
\centering{\psfig{figure=fig4.9.eps, width= 4.1in} \caption{Isomorphism between $G'_1$ and $G'_2$.}}
\end{figure}
\vspace{.2in}

\hspace{1cm} Here again, the black and white vertices added on $G_{1}$ correspond to two copies of $G_{2, e_2}$ which correspond to the triangulations of two disks. Hence, $G'_{1}$ is the dual of the cell decomposition of a 2-sphere.

\hspace{1cm} Suppose that both $P$ and $N$ have exactly one vertex in their triangulations. If we remove, as above, the normal neighborhood of each vertex, we end up with a cell decomposition with no vertices which, of course, is impossible. What we can do, on the other hand, is take a subdivision of a tetrahedron of, say, $P$, to obtain a new triangulation of $P$ with two vertices (see figure below). We then apply the construction above to the new triangulation of $P$ and to $N$. Note that the triangulation of $P \# N$ has $|P_{new}| + |N| + 2 = (|P| + 3) + |N| + 2 = |P| + |N| + 5$ tetrahedra.  

\bigskip
\begin{figure}[h] \label{fig4.10}
\centering{\psfig{figure=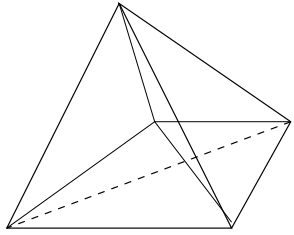, width= 1.5in} \caption{Barycentric subdivision of a tetrahedron.}}
\end{figure}
\vspace{.2in}

\hspace{1cm} Let us describe a construction which requires four more tetrahedra only. Consider the 2-tetrahedron 3-vertex triangulation of $S^{3}$ in the figure below. Take the connected sum of $P$ with $S^{3}$ by removing normal neighborhoods of the vertex $v_1$ in $S^{3}$, and $v$ in $P$ . This vertex $v_1$ as the nice property that its induced 3-valent graph $G$ is 1-edge-connected. For future reference, we will say that $v_1$ is a \textbf{\textit{good vertex}}. Indeed, it is not hard to see that a vertex $A$ is a good vertex if it is the end point of two edges forming a cone.

\bigskip
\begin{figure}[h]  \label{fig4.11}
\centering{\psfig{figure=fig4.11.eps, width= 1.8in}  \caption{A 3-vertex 2-tetrahedron triangulation of $S^3$.}}
\end{figure}
\vspace{.2in}

\hspace{1cm} Looking at figure 4.12 below, we see that $P \# S^{3}$ is made of the tetrahedra from $P$, the tetrahedra from $S^{3}$, plus one extra tetrahedron. In fact, when there is a good vertex in the triangulation of one of two 3-manifolds, there is one red vertex only that needs to be added to $G_1$ and $G_2$ to obtain an isomorphism between $G'_1$ and $G'_2$. Hence, we obtain a triangulation of $P$ with 3 more tetrahedra and 2 vertices. What we need to notice here is that one of the two vertices of the new triangulation of $P$ is also a good vertex (vertex $v_3$ in figure 4.11). Hence, $|P_{new} \# N| = |P_{new}| + |N| + 1 = |P| + 3 + |N| + 1 = |P| + |N| + 4$. 

\bigskip
\begin{figure}[h]  \label{fig4.12}
\centering{\psfig{figure=fig4.12.eps, width= 3.8in} \caption{Only one red vertex needed to obtain isomorphic graphs $G'_1$ and $G'_2$.}}
\end{figure}
\vspace{.2in}

\hspace{1cm} We summarize the above constructions. If $P$ has two vertices, then we can take the connected sum of $P$ and $N$ by inserting 2  tetrahedra unless one of the vertices of $P$ or $N$ is good which would require inserting only 1 tetrahedron. If both $P$ and $N$ have one vertex, then $P \# N$ requires inserting 4 tetrahedra (we will see later that if $P$ or $N$ has a good vertex then this construction requires the insertion of 2 tetrahedra only).

\hspace{1cm} We describe the following constructions:

\begin{construction} \label{cons2}
Let $M_{1}$, ..., $M_{k}$ be closed orientable 3-manifolds with at least 2 vertices in each of their triangulations. Then there exists a triangulation for $M_{1} \#$ ...$M_{k}$ with $\sum |M_{i}| + k + 2$ tetrahedra and at least 2 vertices.
\end{construction}

\bigskip
\begin{figure}[h] \label{fig4.13}
\centering{\psfig{figure=fig4.13.eps, width= 2.5in} \caption{The cell decomposition of a punctured $M_i$.}}
\end{figure}
\vspace{.2in}

\hspace{1cm} Since the $M_{i}$'s have two vertices, $v_{i}$ and $v'_{i}$, it is possible to remove a regular neighborhood of, say, $v'_{i}$, so that we have a truncated tetrahedron as in Figure~\ref{fig4.1}. We insert a prism, and we observe that the resulting cell decomposition represents $M_{i}$ with two balls removed such that the boundary of $M_i /(B \cup B')$ is composed of an annulus and two disks from the sides of the prism.

\bigskip
\begin{figure}[h] \label{fig4.14}
\centering{\psfig{figure=fig4.14.eps, width= 3.2in} \caption{Graph corresponding to the dual of the cell decomposition of the boundary of the punctured manifold $M_2 \# M_3 \#$... $\# M_k$.}}
\end{figure}
\vspace{.2in}

\hspace{1cm} We first construct the connected sum of $M_{1}$ with the 3-sphere described in Figure 4.11 to obtain a new triangulation of $M_{1}$ with at least 3 vertices and exactly 3 more tetrahedra. Note again that one of the vertices, say $A$, is a good vertex. Let $G_1$ be the graph corresponding to the dual of the link of $A$. We now construct a cell decomposition of the connected sum of the $M_i$'s, for $2 \leq i \leq k$.

\begin{construction} \label{cons3}
Let $M$ be a triangulated closed orientable 3-manifold with $|M| = n$ and with exactly one vertex. Then, for $k \geq 3$, there exists a $k$-vertex triangulation for $M$ with $(n + k + 1)$ tetrahedra.
\end{construction}

\hspace{1cm} We first take the connected sum of $M$ with $S^{3}$ as above. We obtain a new triangulation for $M$ with 3 more tetrahedra and 2 vertices, and one of the vertices, say $v$, is a good vertex. We remove a normal neighborhood of $v$ and we change the induced graph $G$ corresponding to $\partial Nbhd(v)$ in the following way:

\bigskip
\begin{figure}[h] \label{fig4.15}
\centering{\psfig{figure=fig4.15.eps, width= 2.4in} \caption{Induced graph $G'$ after inserting $(k - 2)$ (truncated) prisms.}}
\end{figure}
\vspace{.2in}

Note, it is possible to construct such a graph because $G$ is one-connected.

\bigskip
\begin{figure}[h] \label{fig4.16}
\centering{\psfig{figure=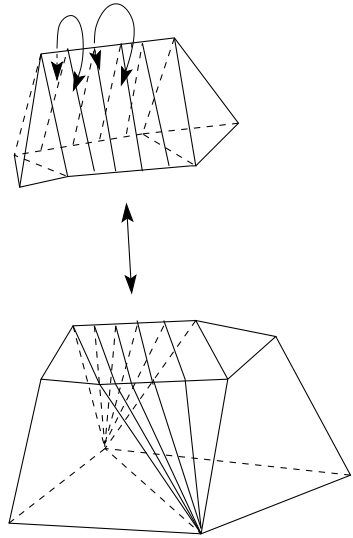, width= 1.8in} \caption{We insert $(k - 2)$ prisms in the new cell decomposition of $M$.}}
\end{figure}
\vspace{.2in}

\begin{construction} \label{cons4}
Let $M$ be a triangulated closed orientable 3-manifold with $|M| = n$ and with exactly $k$ vertices. Then, there exists a $(k+1)$-vertex triangulation of $M$ with $(n + 2)$ tetrahedra. If $M$ has a good vertex, this construction can be done by adding a single tetrahedron, and hence, obtaining a $(k+1)$-vertex triangulation of $M$ with (n + 1) tetrahedra.
\end{construction}

\hspace{1cm} To simplify the notation, suppose $M$ has two vertices. We remove a normal neighborhood of a vertex as in Figure 4.1, and we construct the following graph from $G$. 

\bigskip
\begin{figure}[h] \label{fig4.17}
\centering{\psfig{figure=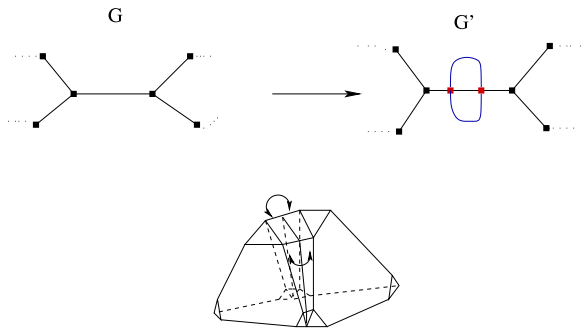, width= 3in} \caption{Insertion of 2 truncated prisms to create a new vertex.}}
\end{figure}
\vspace{.2in}

\begin{construction}
Let $M_{1}$, ..., $M_{k}$ be closed orientable 3-manifolds with 1 vertex in each of their triangulations, and $k \geq 3$. Then there exists a triangulation for $M_{1} \#$ ...$M_{k}$ with $\sum |M_{i}| + 2k$ tetrahedra and 1 vertex.
\end{construction}

\hspace{1cm} To $M_1$, we apply Construction~\ref{cons3} to obtain a 3-manifold homeomorphic to $M$ whose triangulation has $(|M_1| + k + 1)$  tetrahedra and $k$ vertices. Note, all vertices except one are good vertices. Hence, the connected sum of all the $M_i$'s, $2 \leq i \leq k$, requires $(k-1)$ more tetrahedra.  

\bigskip

\begin{construction} \label{cons6}
Let M be a triangulated closed orientable 3-manifold with $|M| = n$ and at least two vertices. Then there exists a triangulation for $M \# L(3, 1)$ with (n + 2) tetrahedra.
\end{construction}

\hspace{1cm} Without loss of generality, we assume we can find a tetrahedron as in Figure~\ref{fig4.1}. Let $M'$ be the manifold obtained from $M$ after removing the link of a vertex. Let $T$ be the triangulation of $\partial M'$, and $G$ its dual. Let $e$ be the thickened edge from Figure~\ref{fig4.1}.  We insert 2 prisms along the shaded face of Figure~\ref{fig4.1}, and we insert some tips of tetrahedra so that the new triangulation $T'$ of $\partial M'$ and its dual $G'$ correspond to the following:

\bigskip
\begin{figure}[h] \label{fig4.27}
\centering{\psfig{figure=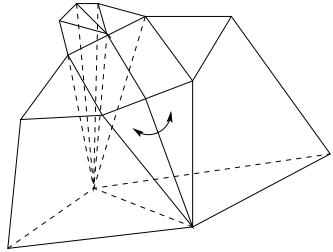, width= 1.7in} \caption{The new cell decomposition of $M'$.}}
\end{figure}
\vspace{.2in}

\bigskip
\begin{figure}[h] \label{fig4.28}
\centering{\psfig{figure=fig4.28.eps, width= 3in} \caption{The new triangulation $G'$ of $G$.}}
\end{figure}
\vspace{.2in}

\hspace{1cm} Since $G'$ is planar, it is the dual of a cell decomposition of a 2-sphere, and hence, we did not change the topology of $M'$. Consider the following cell decomposition of the punctured $L(3,1)$ (which we will call $N'$):

\bigskip
\begin{figure}[h] \label{fig4.29}
\centering{\psfig{figure=fig4.29.eps, width= 2.9in} \caption{A cell decomposition of $N'$ and the dual $G_1$ of $\partial N'$.}}
\end{figure}
\vspace{.2in}

\hspace{1cm} We now cut along face 4 and insert truncated prisms so that the resulting graph $G'_1$ of $G_1$ is isomorphic to $G'$. It is now clear that the triangulations of $\partial M'$ and $\partial N'$ matches. We glue the two boundaries together to obtain a well-defined triangulation of $M \# L(3,1)$ with $(n+2)$ tetrahedra.

\begin{construction} \label{cons7}
Let M be a triangulated closed orientable 3-manifold with $|M| = n$ and at least two vertices. Then there exists a triangulation for $M \# \mathbf{RP}^{3}$ with (n + 2) tetrahedra.
\end{construction}

\hspace{1cm} Without loss of generality, we assume we can find a tetrahedron as in Figure~\ref{fig4.1}. Let $M'$ be the manifold obtained from $M$ after removing the link of a vertex. Let $T$ be the triangulation of $\partial M'$, and $G$ its dual. Let $e$ be the thickened edge from Figure~\ref{fig4.1}. We insert 4 prisms along the shaded face of Figure~\ref{fig4.1}, and we insert some tips of tetrahedra so that the new triangulation $T'$ of $\partial M'$ and its dual $G'$ correspond to the following:

\bigskip
\begin{figure}[h] \label{fig4.24}
\centering{\psfig{figure=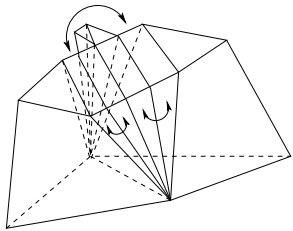, width= 1.5in} \caption{The new cell decomposition of $M'$.}}
\end{figure}
\vspace{.2in}

\bigskip
\begin{figure}[h] \label{fig4.25}
\centering{\psfig{figure=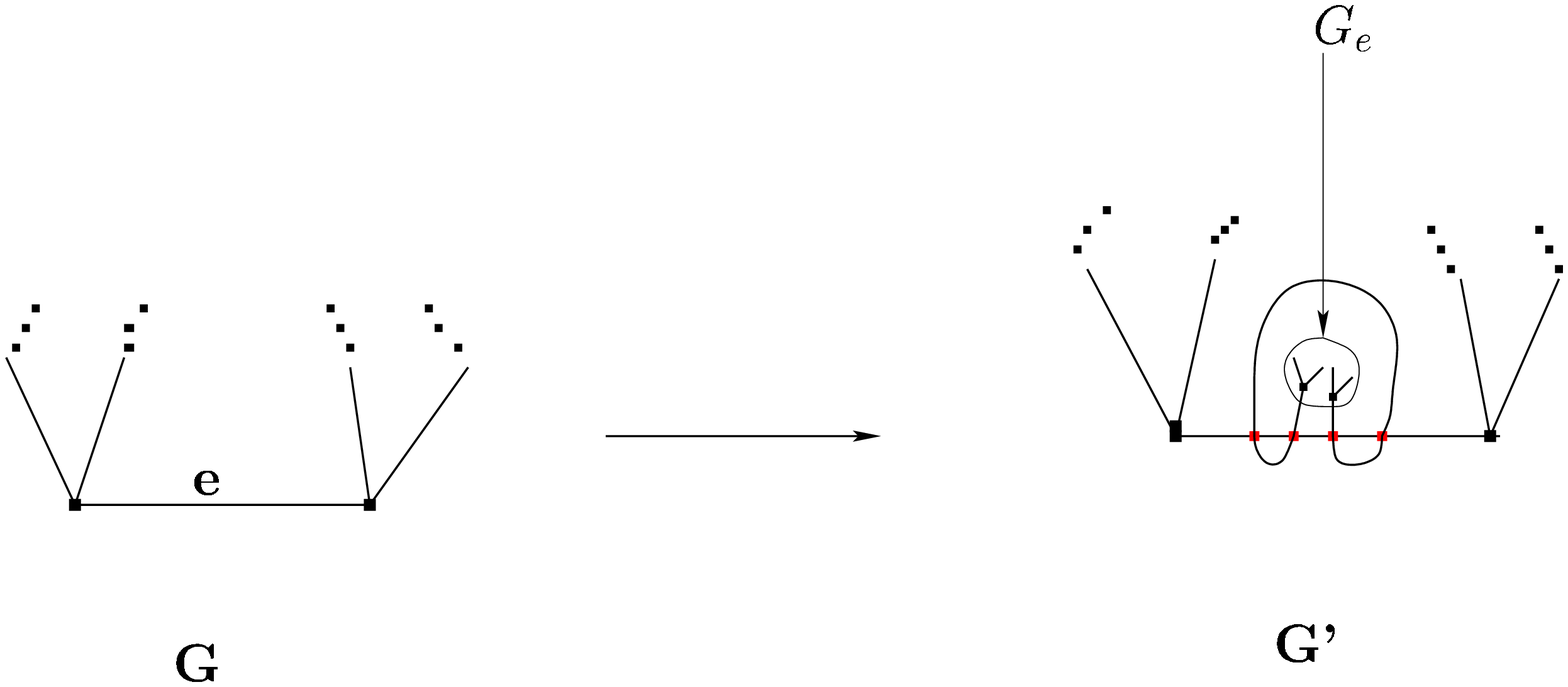, width= 3in} \caption{The new triangulation $G'$ of $G$.}}
\end{figure}
\vspace{.2in}

\bigskip
\begin{figure} \label{fig4.26}
\centering{\psfig{figure=fig4.26.eps, width= 2.9in} \caption{A cell decomposition of $N'$ and the dual $G_1$ of $\partial N'$.}}
\end{figure}
\vspace{.2in}

\hspace{1cm} Note, $G'$ is planar and so we did not changed the topology of $M'$. Consider now the following cell decomposition of the punctured $\mathbf{RP}^3$ (which we will call $N'$):

\hspace{1cm} We now cut along face 4 and insert thickened triangles (which we think as truncated tips) so that the resulting graph $G'_1$ of $G_1$ is isomorphic to $G'$. At this point, it is clear that the triangulations of $\partial M'$ and $\partial N'$ matches. We glue the two boundaries together to obtain a well-defined triangulation of $M \# \mathbf{RP}^3$ with $(n+2)$ tetrahedra.

\vspace{1cm}

\begin{construction} ~\label{cons8}
Let M be a triangulated closed orientable 3-manifold with $|M| = n$ and at least three vertices. Then there exists a triangulation for $M \# (S^{1} \times S^{2}$) with (n + 2) tetrahedra. If M contains a good vertex, this construction can be done with one extra tetrahedron only.
\end{construction}

\hspace{1cm} This construction is almost identical to taking the connected sum of two triangulated 3-manifolds, where one triangulation has at least 2 vertices. As it was described in the claim at the end of Theorem~\ref{MT}, in order to construct $M \# (S^{1} \times S^{2}$), we simply need to remove two disjoint 3-balls from $M$ and identify their boundaries, $S_1$ and $S_2$, through a homeomorphism. Let us remove the normal neighborhoods of two vertices. Because $M$ has three vertices and is connected, we can find the following truncated tetrahedra:

\bigskip
\begin{figure} ~\label{fig4.18}
\centering{\psfig{figure=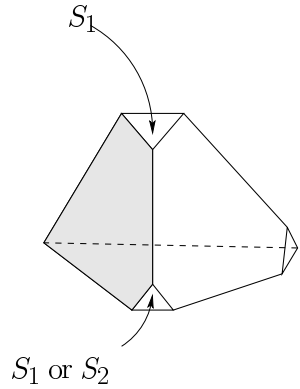, width= 1.5in} \caption{A truncated tetrahedron in the cell decomposition of the twice punctured manifold with boundary components $S_1$ and $S_2$.}}
\end{figure}
\vspace{.2in}

\begin{proposition}
There exists a 1-vertex triangulation of any closed orientable 3-manifold.
\end{proposition}

\hspace{1cm} \underline{{\bf Proof :}} Let $M$ be a closed orientable reducible 3-manifold equipped with a triangulation. Suppose $M$ has more than 1 vertex. Contrary to Proposition~\ref{Pro1}, cutting along the boundary of the neighborhood of a maximal tree does not solve our problem. The reason is that there may be more than one $S_{i}$ not bounding a ball on the inside (in fact, the $S_{i}$'s may not even be separating).

\hspace{1cm} We consider the 1-vertex 1-tetrahedron triangulation of $S \sp 3$. Then $M \# S^{3}$ has $(t + 3)$ tetrahedra and $(v-1)$ vertices by our Construction~\ref{cons1}. We can repeat this construction $(v-2)$ times and obtain a manifold homeomorphic to $M$ with $(t+3(v-1))$ tetrahedra and 1-vertex only. Q.E.D.

\medskip

\section{Small normal 2-Spheres in Minimal Triangulations}

\hspace{1cm} Using the constructions in the previous section and Theorem~\ref{MT}, we show that, in some sense, minimal triangulations of reducible 3-manifolds contain ``small" non-trivial normal 2-spheres. By small, we mean normal 2-spheres which have quadrilateral types in ``few" tetrahedra. This will be extremely useful in Andrew Casson's algorithm to check if a minimal triangulation is reducible or not.  

\smallskip

\begin{definition} A 3-manifold $M$ is said to have a \textbf{\textit{minimal triangulation}} $\tau$ if $\tau$ contains the smallest number of tetrahedra over all possible triangulations of $M$ (by abuse of language, we say that $M$ is minimal). 

\smallskip

\hspace{1cm} The \textbf{\textit{weight}} of a surface $S$, denoted by $wt(S)$, is the number of points in the intersection of $S$ and the 1-skeleton, i.e. $wt(S) = card (\{ S \cap T^{(1)} \})$. 

\smallskip

\hspace{1cm} We will say that a normal surface $S$ has $n$ \textbf{\textit{tetraquads}} (denoted by $<S>$ = n) if there exist exactly n tetrahedra with the property that S intersects each of these tetrahedra in at least one quadrilateral type. Note, if S and T are two compatible normal surfaces with $<$S$>$ = n and $<$T$>$ = m, then $(n+m) \geq$ $<$$S + T$$>$ $\geq $max(n, m). Moreover, normalizing an embedded surface decreases its weight, but it may increase its number of tetra-quads.

\smallskip

\hspace{1cm} We will denote by $\mathbf{\#_q (F)}$, the number of quadrilaterals representing the normal surface $F$.

\end{definition}

\smallskip

\hspace{1cm} To prove the existence of small normal 2-spheres in minimal triangulations, we first want to show the existence of a normal 2-sphere with the property that all the surgery surfaces are inessential. If we can find such a 2-sphere, then cutting along it and collapsing it to a point will result in a decomposition of M in exactly two summands. Because $M$ is minimal, we use the constructions in the previous section to conclude that $S$ cannot have too many tetraquads.

\medskip

\hspace{1cm} A naive approach would be to first look at a non-trivial normal 2-sphere with minimal weight over all non-trivial normal ones. There is a problem though. Suppose $S$ is a non-trivial normal 2-sphere and has minimal weight. Let $A$ be a surgery annulus as in Lemma~\ref{L2}. If $A$ is essential then we can find a non-trivial normal 2-sphere, parallel to $T \cup D_1 \cup D_2$ (see Lemma~\ref{L2} for notation), with smaller weight, and we reach a contradiction. On the other hand, let $D$ be an essential surgery disk. The problem here is that the normalization of the two 2-spheres $D \cup D_1$ and $D \cup D_2$ (see Lemma~\ref{L3} for notation) may give us 2-spheres with bigger weights. Hence, having a minimal weight 2-sphere is not enough to conclude that all surgery surfaces are inessential. 

\medskip

\hspace{1cm} Another approach would be to look at a non-trivial normal 2-sphere which has the smallest number of quadrilaterals. So suppose that $D$ is an essential disk. We construct in Lemma~\ref{L4} a new non-trivial normal 2-sphere which has a strictly smaller number of quadrilaterals. This is a contradiction and hence, $D$ had to be inessential. Let $A$ be an essential annulus. The problem here is that we haven't been able to construct a non-trivial normal 2-sphere with a strictly smaller number of quadrilaterals. In particular, we construct a 2-sphere in Lemma~\ref{L4} which may have the same number of quadrilaterals as $S$.

\medskip

\hspace{1cm} Here is another idea, without a proof, to go around these obstructions. Let $S$ be a non-trivial normal 2-sphere with the smallest number of essential surgery surfaces over all possible non-trivial normal 2-spheres. Let $A$ be an essential surgery annulus or disk. We construct a new 2-sphere $S'$ which is described in Lemma~\ref{L4}. One can show that $S'$ has strictly less essential surgery surfaces than $S$. This is a contradiction and hence, $A$ had to be inessential.

\bigskip

\hspace{1cm} To show the existence of a non-trivial normal 2-sphere with no essential surgery surfaces, we are going to use a result due to W. Jaco and J. Tollefson ~\cite{JT:gnus}. Here is the theorem:

\medskip

\hspace{1cm} \underline{Theorem 4.1} (~\cite{JT:gnus}): A normal two-sphere $F$ is a vertex surface if and only if $F$ has the property that whenever there exists an annulus $A$ which is an exchange surface for $F$ then the two disjoint disks in $F$ bounded by $\partial A$ are normal isotopic.

\medskip

\hspace{1cm} Here, an \textbf{\textit{exchange surface}} $A$ for $F$ is a surface with the following properties: 1) $fr(A) = A \cap F$, 2) $A$ has an orientable regular neighborhood $N(A)$, and 3) for every tetrahedron $\Delta$, each component of $\Delta \cap A$ is a 0-weight disk $L$ spanning two distinct elementary disks $E_1$, $E_2$ of $F$ such that $\partial L = L \cap (E_1 \cup E_2 \cup \partial \Delta)$ and $L \cap E_i$ is an arc joining the interiors of two distinct 2-faces of $\Delta$.

\hspace{1cm} It is clear from this definition that if $A$ is a surgery annulus or a surgery Mobius band, then we can either push $A$ off the 2-skeleton or look at $\partial Nbhd(A)$ to find an annulus which is an exchange surface. Hence, if $F$ is a vertex 2-sphere, we now know that any surgery annulus or Mobius band is inessential. On the other hand, Jaco and Tollefson's theorem does not say anything about surgery disks. So, to go around this problem, we define a complexity $Q$. Let $S$ be a normal surface. Then \textbf{\textit{Q(S)}} represents a pair, ordered lexicographically, whose first and second entries are $< S >$ and $\#_q (S)$ respectively: Q(S) = ($<S>, \#_q (S)$). 

We are now ready to prove the main Lemma.

\bigskip

\begin{lemma} \label{L4}
Let M be a triangulated closed orientable 3-manifold which contains a non-trivial normal 2-sphere. Then M contains a non-trivial normal 2-sphere whose surgery surfaces are all inessential.
\end{lemma}

\smallskip

\hspace{1cm} \underline {{\bf Proof :}} Consider the set of non-trivial normal 2-spheres. This set is non-empty by assumption. In this set, choose the 2-sphere $F$ which is minimal with respect to $Q$. We want to show that such a 2-sphere does not have any essential surgery surfaces, but first, we want to show that it is a vertex surface. See Chapter 1 for the definitions of vertex surfaces and vertex solutions.

\hspace{1cm} Suppose $F$ is not a vertex surface, i.e. suppose $x_F$ is not a vertex of \textit{P(M, T)}, i.e. suppose for all positive integer $k$, $k \cdot x_F$ is not a vertex solution. Let $k$ be the smallest positive integer such that $k \cdot x_F = S$ is an integral solution (and, by assumption, not a vertex solution). Then $\chi(S)$ must divide $\chi(F)$, and so $\chi(S) = 1$ or 2.

\hspace{1cm} \underline{Case 1:} $\chi(S) = 2$. Since $F$ is an integral multiple of $S$, $S = F$. Suppose there exists a positive integer $n$ such that $nF = V_1 +$... $+ V_k + W_1 +$... $+W_r + X$, where at least one of the summands is not an integral multiple of $F$. Without loss of generality, we can assume that the $V'_i$s are 2-spheres, the $W'_i$s are real projective planes, and $X$ is a (possibly empty and possibly disconnected) surface with non-positive Euler characteristic. Because the Euler characteristic is preserved under surface addition, we have $2k + r \geq 2n$. 

\hspace{1cm} \underline{Subcase 1:} Suppose that $r = 0$. Then $nF = V_1 +$... $+ V_k + X$ with $k \geq n$. Consider the 2-sphere, say $V_1$, which has the smallest $\#_q$ over all the $V_i$'. Since the surface addition preserves the quadrilateral types and the number of quadrilaterals for each quadrilateral type, we have $< V_1 > \leq < F >$. But by assumption, $< F > \leq < V_1>$. Hence, $< F > = < V_1 >$. It is crucial to notice that, not only $F$ and $V_1$ have the same number of tetraquads, but they must also have their quadrilaterals in the same tetrahedra. Moreover, $n \cdot (\#_q ( F )) \geq k \cdot (\#_q ( V_1 )) + \#_q ( X )$. Since $k \geq n$ and $\#_q (V_1) \geq \#_q(F)$, the only way for the inequality to be true is if it is an equality and $k = n$ and $\#_q( X ) = 0$. We conclude that $X$ was actually empty. Also, $\#_q (V_1) \leq \#_q(F)$. By assumption, $\#_q (F) \leq \#_q(V_1)$, so $\#_q(F) = \#_q(V_1)$. Therefore, $V_1$ is a parallel copy of $F$. We remove $V_1$, X,  and a copy of $F$ from the equation to obtain $(n-1)F = V_2 +$... $+ V_n$. We repeat the argument for the next 2-sphere, say $V_2$, which as the smallest $\#_q$. We conclude that $V_i$ is a parallel copy of $F$ for $1 \leq i \leq n$. This contradicts the fact that $F$ was not an integral solution.

\hspace{1cm} \underline{Subcase 2:} Suppose $ r \neq 0$. We look at the equation $2nF = 2V_1 +$... $+ 2V_k + 2W_1 +$... $+2W_r + 2X$, where $2V_i$ represents two copies of a 2-sphere and $2W_j$ represents a 2-sphere. We repeat the same argument as in Subcase 1 to conclude that each $V_i$ is a parallel copy of $F$, that each $2W_i$ is a parallel copy of $F$, and that $X$ is empty. Since, say $2W_1$ is parallel to $F$, we conclude that $F$ is the double of the projective plane $W_1$. This contradicts the fact that $k$ was the smallest positive integer such that $S$ (= $F$) is an integral solution.

\smallskip

\hspace{1cm} \underline{Case 2:} $\chi(S) = 1$. Since $F$ is an integral multiple of $S$, $F = 2S$. Suppose there exists a positive integer such that $nS = V_1 +$... $+ V_k + W_1 +$... $+W_r + X$, where the $V_i$'s are 2-spheres, the $W_j$'s are real projective planes, and $\chi(X) \leq 0$. We look at the new equation $nF = 2nS = 2V_1 +$... $+ 2V_k + 2W_1 +$... $+ 2W_r + 2X$. We run the same argument as in Subcase 1 of Case 1 to conclude that each $V_i$ and each $2W_j$ is a parallel copy of $F$, and that $X$ is empty. Hence, each $V_i$ and $2W_j$ represents the double of $S$. Therefore, each summand of the equation is an integral multiple of $S$. We conclude that $S$ is a vertex solution, or equivalently, $F$ is a vertex surface. Contradiction.

\medskip

\hspace{1cm} Using the above theorem from Jaco and Tollefson, we now know that, if $F$ is a non-trivial normal 2-sphere which minimizes $Q$, then $F$ does not have any inessential surgery annuli or Mobius bands.

\smallskip

\hspace{1cm} Let $D$ be a surgery disk. We want to construct another non-trivial normal 2-sphere with a smaller number of quadrilaterals. Note, there are two ways of constructing the 2-spheres $D \cup D_1$ and $D \cup D_2$ as in Lemma~\ref{L3}. One way is to normalize, say $D \cup D_1$, and obtain a possibly disconnected surface whose connected components are all normal 2-spheres. Another way is to discard one of the disks, say $D_2$, and replace it by triangles. 

\bigskip
\begin{figure}[h] ~\label{fig4.20}
\centering{\psfig{figure=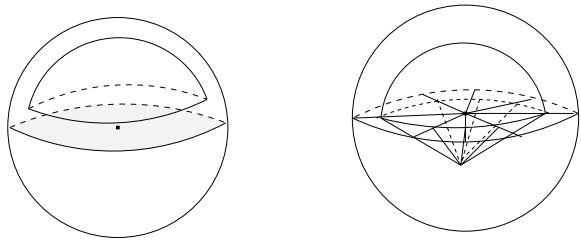, width= 3in} \caption{Construction of $D \cup D_1$.}}
\end{figure}
\vspace{.2in}

\hspace{1cm} We note that in the latter case, the number of quadrilaterals is strictly reduced, and unless the surgery disk is inessential, we obtain a new 2-sphere which is non-trivial. Moreover, this new 2-sphere $D \cup D_1$ is already in normal form. This proves the Lemma.

\bigskip

We are now ready to prove the main theorem of this Chapter. 

\begin{theorem} ~\label{MT3}
Let M be a closed orientable 3-manifold equipped with a minimal triangulation. If M is reducible, then M contains a non-trivial normal 2-sphere S' such that $<S'> \leq$ 4.
\end{theorem}

\underline {{\bf Proof :}}

\hspace{1cm} Because $M$ is reducible, Kneser finiteness theorem ~\cite{Kn:gnus} tells us that there exists a normal essential 2-sphere in $M$. So by Lemma~\ref{L4} we know there exists a non-trivial normal 2-sphere with only inessential surgery surfaces. If $<S> \leq 4$, then we are done. So we will assume that $<S>$ $\geq 5$ and we will either contradict the fact that M is minimal or show that $M$ is homeomorphic to one of the following spaces: $L(3, 1) \# L(3, 1)$, $L(3, 1) \# \mathbf{RP}^3$, $\mathbf{RP}^3 \# \mathbf{RP}^3$, $L(3, 1)$, $\mathbf{RP}^3$, or $S^3$. As in the proof of Theorem~\ref{MT}, we want to look for surgery annuli, Mobius bands, or disks. If there are none, then we follow the proof of Theorem~\ref{MT} and we may obtain 2 different decompositions: 

\hspace{1cm} 1) $M \cong M_{1} \# M_{2}$, with $|M_{1}| + |M_{2}| + 5 \leq |M|$ and either $M_{1}$ or $M_{2}$ has more than 1 vertex in its triangulation. Indeed, when we cut along $S$ and collapse the two boundary components of its regular neighborhood to points, we create two new vertices. Since $M$ has at least one vertex in its original triangulation, $M_{1}$ or $M_{2}$ must have at least two vertices. We use Construction~\ref{cons1} and obtain a triangulation of $M$ with 3 less tetrahedra. Contradiction.

\hspace{1cm} 2) $M \cong M_{1} \# (S^{1} \times S^{2})$ ,with $|M_{1}| + 5 \leq |M|$ and $M_{1}$ has at least 3 vertices. Here, the three vertices come from the collapsing of the two parallel copies of $S$ and the vertex in the original triangulation of $M$. We apply Construction~\ref{cons8} and we obtain a triangulation of $M$ with 3 less tetrahedra. Contradiction.

\smallskip

\hspace{1cm} Hence, we need to assume the existence of such surgery surfaces. Let $S$ be the non-trivial normal 2-sphere, constructed in Lemma~\ref{L4}, which does not have any essential surgery surfaces. Hence, after collapsing $S$, we can get at most two non-trivial summands for $M$. 

\hspace{1cm} Suppose first that $S$ is non-separating. Let $A$ be an inessential surgery annulus. Then, following the proof of Theorem~\ref{MT}, we obtain the following decomposition:

\hspace{1cm} $M \cong M_{1} \# (S^{1} \times S^{2})$. Note, $M_{1}$ may have only two vertices and so taking the connected sum of $M_{1}$ with $(S^{1} \times S^{2})$ may require 4 additional tetrahedra: we first add a vertex to $M$ as in Construction~\ref{cons4}. This increases the number of tetrahedra by 2. We now apply Construction~\ref{cons8}, which increases by 2 the number of tetrahedra in the new triangulation. Hence, we obtain a new triangulation of $M$ which has 4 more tetrahedra than $M_{1}$. Because $<S> \geq 5$, $|M_{1}| + 4 < |M|$ and we obtain a new triangulation of $M$ with fewer tetrahedra. Contradiction.

\hspace{1cm} Let $A$ be a surgery disk. We obtain the same decomposition but in this case, $M_{1}$ has 3 vertices and we can apply Construction~\ref{cons8} or ~\ref{cons1}. We get a new triangulation of $M$ with at least 3 fewer tetrahedra. This contradicts the minimality of $M$.

\bigskip

\setlength{\parskip}{0.03in}

\hspace{1cm} Suppose now that $S$ is separating. Note, because all surgery surfaces are inessential, it is only necessary to look at the cases where there is exactly one surgery surface on each side of $S$.

{\bf 1-} \space \space Let $A$ and $A'$ be two surgery annuli on each side of $S$.

\hspace{1cm} {\bf Case 1:} Both $A$ and $A'$ bound solid tori. What we mean by a surgery annulus bounding a solid torus is the following: let $A$ be a surgery annulus embedded as in case 1 of Lemma~\ref{L1}. We say, by abuse of language, that $A$ bounds a solid torus if $A_{1} \cup A'_{1}$ bounds a solid torus made of truncated prisms. Because $A$ and $A'$ are both inessential, we have the following decompositions: $M \cong L(3, 1) \# L(3, 1)$, $M \cong L(3, 1) \# S^{3}$, or $M \cong S^{3} \# S^{3}$.

\hspace{1cm} {\bf Case 2:} One annulus bounds a solid torus and the other one doesn't. There are two possibilities: either $M \cong M_{1} \# L(3, 1)$ or $M \cong M_{1} \# S^{3}$. The latter case clearly leads to a contradiction. In the former case, $M_{1}$ may have only one vertex. What we do here is take the connected sum of $M_{1}$ with the 2-vertex 2-tetrahedron triangulation of $L(3, 1)$. This results in a triangulation of $M$ with 4 more tetrahedra then $M_{1}$. Because $<S> \geq 5$, this new triangulation of $M$ has strictly less tetrahedra than the original one. Contradiction. 

\hspace{1cm} {\bf Case 3:} No annuli bound solid tori. We have: $M \cong M_{1} \# M_{2}$. In the worst case, both $M_{1}$ and $M_{2}$ have one-vertex triangulations. We can still apply construction 4.1 and get a new triangulation for $M$ with at least one less tetrahedron. Contradiction.

\setlength{\parskip}{0.1in}

\vspace{.2in}

\medskip

\setlength{\parskip}{0.03in}

{\bf 2-} \space \space Let $A$ be a surgery annulus on one side of $S$ and $B$ a surgery Mobius band on the other side. 

\hspace{1cm} {\bf Case 1:} $A$ bounds a solid torus. We have: $M \cong L(3, 1) \# \mathbf{RP}^{3}$ or $M \cong S^{3} \# \mathbf{R}P^{3}$.

\hspace{1cm} {\bf Case 2:} $A$ doesn't bound a solid torus. We have: $M \cong M_{1} \# \mathbf{RP}^{3}$. Since $M_{1}$ may have only one vertex, we take the connected sum of $M_{1}$ with the 2-vertex 2-tetrahedra triangulation of $\mathbf{RP}^{3}$. We then obtain a new triangulation of $M$ with 4 more tetrahedra then $M_{1}$. Since $<S> \geq 5$, we constructed a triangulation of $M$ with strictly less tetrahedra then the original one. Contradiction. 

\setlength{\parskip}{0.1in}

\bigskip
\begin{figure}[h] ~\label{fig4.22}
\centering{\psfig{figure=fig4.22.eps, width= 2.2in} \caption{An inessential surgery Mobius band on one side of $S$ and an inessential surgery annulus, on the other side, not bounding a solid torus.}}
\end{figure}
\vspace{.2in}

\medskip

{\bf 3-} \space \space Let $B_{1}$ and $B_{2}$ be two surgery Mobius bands on each side of $S$. We have: $M \cong \mathbf{RP}^{3} \# \mathbf{RP}^{3}$.

\medskip

{\bf 4-} \space \space Let $D_{1}$ and $D_{2}$ be two surgery disks on each side of $S$. We have: $M \cong M_{1} \# S^{3}$, $M \cong M_{1} \# M_2$, or $M \cong S^{3} \# S^{3}$. The former case clearly leads to a contradiction. If $M \cong M_{1} \# M_{2}$, $M_{1}$ and $M_{2}$ must each have at least two vertices (which come from the collapsing of the surgery disks). So we can apply Construction 4.1 to obtain a triangulation of $M$ with 3 less tetrahedra. Contradiction. 

\medskip

\setlength{\parskip}{0.03in}

{\bf 5-} \space \space Let $D$ be a surgery disk on one side of $S$ and $A$ a surgery annulus on the other side.

\hspace{1cm} {\bf Case 1:} $A$ bounds a solid torus. We have: $M \cong M_{1} \# L(3, 1)$, $M \cong M_{1} \# S^{3}$, $M \cong S^{3} \# S^{3}$, or $M \cong S^{3} \# L(3, 1)$. In the former case, $M_{1}$ must have at least 2 vertices. We use Construction~\ref{cons6} to obtain a new triangulation of $M$ with 3 less tetrahedra than the original one. Contradiction. The other cases clearly contradict the minimality of $M$.

\hspace{1cm} {\bf Case 2:} $A$ does not bound a solid torus. We have:  $M \cong S^{3} \# M_{2}$ or $M \cong M_{1} \# M_{2}$. The former case is a clear contradiction. In the latter case, $M_{1}$ has at least 2 vertices and so we apply Construction 4.1 to obtain a new triangulation of $M$ with 3 less tetrahedra than the original one. 

\bigskip

{\bf 6-} \space \space Let $D$ be a surgery disk on one side of $S$ and $B$ a Mobius band on the other side. We then have $M \cong M_{1} \# \mathbf{RP}^{3}$ or $M \cong S^{3} \# \mathbf{RP}^{3}$. In the former case, $M_{1}$ has at least two tetrahedra so we can apply Construction 4.7 and contradict the minimality of $M$. This proves the Theorem.

We make the following remark. Let $M$ be homeomorphic to one of the following spaces:$ L(3, 1) \# L(3, 1)$, $L(3, 1) \# \mathbf{RP}^{3}$, $\mathbf{RP}^{3} \# \mathbf{RP}^{3}$, $L(3, 1)$, $\mathbf{RP}^{3}$, or $S^{3}$. By applying Constructions 4.6 and 4.7, we can construct 4-tetrahedron triangulations of the first three spaces. The last three spaces have known minimal triangulation with at most two tetrahedra. See ~\cite{JR:gnus} for more details on the minimal triangulations of $L(3, 1)$ and $\mathbf{RP}^3$. This shows that, in a minimal triangulation of any of these spaces, there is a non-trivial 2-sphere $S$ with $<S> \leq 4$. Therefore, if $M$ is a closed orientable 3-manifold with a minimal triangulation containing a non-trivial normal 2-sphere, then there exists a non-trivial normal 2-sphere with at most 4 tetraquads.

\bigskip

This result can actually be strengthened. Indeed, we can show that, under the same hypothesis, {\bf there exists a non-trivial normal 2-sphere with at most 2 tetraquads}. It uses Jaco and Tollefson's theorem (see ~\cite{JT:gnus}) stated earlier. If F is a vertex 2-sphere and $A$ is a surgery annulus, the theorem says that not only $A$ is inessential, but also that the disks $D_1$ and $D_2$ are isotopic. Precisely, it says that, if $A$ is a surgery annuli, then the ball $B$ described in Lemma~\ref{L3} does not contain a vertex of the triangulation. Let's take a closer look at the different cases in the proof which require adding four extra tetrahedra.

In case 2 of 1-, $M_{1}$ must have at least 2 vertices since the two disks on $S$ with boundary $\partial A$ are normal isotopic. The same reasoning can be used in case 3 of 1-, case 2 of 2-, and case 1 and 2 of 5-. Hence, all the construction requiring 4 extra tetrahedra can be replaced by constructions with 2 extra tetrahedra. 

\bigskip
\begin{figure}[h] ~\label{fig4.23}
\centering{\psfig{figure=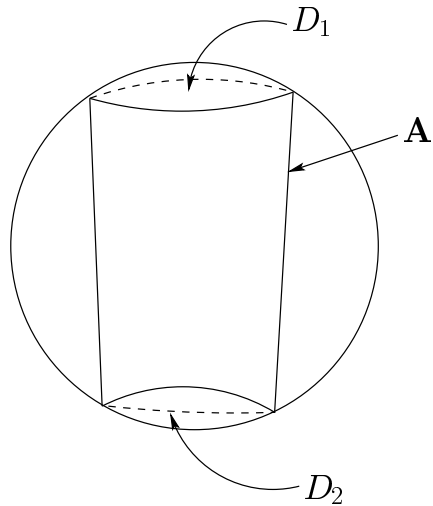, width= 2.2in} \caption{No vertices in the 3-ball bounded by $A \cup D_1 \cup D_2$.}}
\end{figure}

\vspace{.2in}

This shows that, under the same hypothesis, there is a non-trivial normal 2-sphere with at most 2 tetraquads, or it can be shown that $M$ is homeomorphic to one of the 6 spaces described earlier. If $M$ is homeomorphic to either $L(3, 1) \# L(3, 1)$, $L(3, 1) \# \mathbf{RP}^{3}$, or $\mathbf{RP}^{3} \# \mathbf{RP}^{3}$, we can construct a 4-tetrahedron triangulation of these spaces, using the earlier Constructions, with the existence of a normal 2-sphere with 2 tetraquads. Unfortunately, we do not know if every 4-tetrahedron triangulation of these spaces contains such a 2-sphere, so we will assume in the following corollary that $M$ has at least 5 tetrahedra in its triangulation.

\bigskip

\begin{corollary} ~\label{C3}
Let M be a closed orientable 3-manifold equipped with a minimal triangulation with t tetrahedra and $t \geq 5$. Then there is an algorithm to check if M is reducible or not, and this algorithm runs in polynomial time with respect to t.
\end{corollary}

\underline{{\bf Proof :}}
We will assume that the reader is familiar with the terminology used in section 6.1: Casson's algorithm to decompose M into irreducible pieces. What makes Casson's algorithm run in exponential time is that we have to look through, possibly all, the cones $C_{w,i}(M)$ (there are $3^t$ of them) in order to find a non-trivial normal 2-sphere and we have to look through, possibly all, the cones $A_{w,i,l}$ (there are $4t \cdot 3^{t}$ of them) to look for an almost$^{2}$ normal 2-sphere to check if one of the resulting summands is homeomorphic to a 3-sphere. In our case, not only we do not need to look through all the cones $C_{w,i}(M)$, but we don't even need to run the Thompson-Rubinstein algorithm.

1) If there are no non-trivial embedded normal 2-spheres, then a famous result of Kneser (~\cite{Kn:gnus}) tells us that $M$ does not contain any embedded essential 2-sphere. Hence, $M$ is irreducible.

2) Suppose now that $M$ contains a non-trivial normal 2-sphere $S$. Without loss of generality, we can assume that $S$ is a vertex surface. By the main theorem of ~\cite{JT:gnus}, any surgery disk, annulus, or Mobius band is inessential. Hence, when we follow the proof of Theorem~\ref{MT} to cut along $S$ and collapse it, we end up with at most 2 summands for $M$. None of the summands can be homeomorphic to a 3-sphere. Indeed, the existence of a 3-sphere would contradict the minimality of $M$. This tells us that $S$ is essential. Hence, $M$ is reducible. What is important to notice here is that, if $S$ does not have any essential surgery annuli, disks, nor Mobius band, then we must have $<S> \leq 2$ by the earlier constructions.

 Therefore, if $M$ contains a non-trivial normal 2-sphere, then we know that there exists a non-trivial normal 2-sphere $S$ with 1 or 2 tetraquads. For the reader not familiar with the notation in Casson's Algorithm, we advise the reader to first read Chapter 6 and then to come back to this section. Consider the space of normal surfaces represented by a family of type $w$, where $w$ assigns the value 0 for the quadrilateral types in every tetrahedra except one. We noted earlier that surface addition may increase the number of tetraquads. Indeed, if $F_{1}$ has exactly one quadrilateral type in the $j^{th}$ tetrahedra $\Delta_{j}$, and $F_{2}$ has exactly one quadrilateral type in the $i^{th}$ tetrahedra, then $F_{1} + F_{2}$ has a quadrilateral type in $\Delta_{j}$ and $\Delta_{i}$. Hence this space does not represent a cone. Consider now the space of normal surfaces represented by a family, $w_{i}$, of type $w$, where $w$ assigns the value 0 for the quadrilateral types in every tetrahedra except in the $i^{th}$ one. It is now easy to see that this space is a cone in $\mathbf{R}^{7t}$. Every normal surface in this cone has one tetraquad, and conversely, every normal surface with one tetraquad belong to such a cone for some $w_{j}$. We denote it by $C_{w_{j},i}(M)$, where $w_{j}$ represents the type described above with nonzero image in $\Delta_{j}$, and $i$ represents the entry $t_{i}$ being zero. How many such cones are there? For each tetrahedron, there are three different quadrilateral types and so there are three possible $w_{j}$'s. For each type $w_{j}$, there are $4t$ choices for the triangle entries $t_{i}$. Hence, there are $3 \cdot t \cdot (4t) = 12t^{2}$ different cones. 

Similarly, we define the space of normal surfaces having two quadrilateral types in exactly two distinct tetrahedra. If the two tetrahedra are fixed, this space represents a cone. In fact the union of these cones, over all possible pairs of tetrahedra, represents the space of normal surfaces with two tetraquads. We call each cone $C_{w_{j, k},i}(M)$. We count how many of these cones there are. For each pair of tetrahedra, $\{ \Delta_{j}, \Delta_{k} \}$, there are $3 \cdot 3$ possible types $w_{j, k}$. There are ${t \choose 2}$ possible  pairs of distinct tetrahedra. Hence, there are $3 \cdot 3 \cdot {t \choose 2} \cdot (4t) = 36t{t \choose 2}$ different cones $C_{w_{j, k},i}(M)$.

\medskip

We now describe the algorithm to check if a minimal closed orientable 3-manifold $M$ is reducible or not. Fix a cone $C_{w_{j, k},i}(M)$ and look at the convex polyhedron $A = C_{w_{j, k},i}(M) \bigcap$ $\{\sum^{7t}_{i=1} t_{i} = 1 \} $. We maximize $\chi$ on $A$ to obtain a vertex solution S. If $\chi (S) > 0$, then there exists a 2-sphere in $C_{w_{j, k},i}(M)$ and the procedure stops here: $M$ is reducible.
If $\chi (S) \leq 0$, there are no 2-spheres in the cone $C_{w_{j, k},i}(M)$.
Repeat this step with a new cone $C_{w_{r, s},i}(M)$.
If no 2-spheres have been found in any of the cones $C_{w_{j, k},i}(M)$, then $M$ is irreducible and the procedure stops here. Note, we do not need to look through the cones $C_{w_{j},i}(M)$. Indeed, each cone $C_{w_{j},i}(M)$ lies in a cone $C_{w_{j, k},i}(M)$.

What is the complexity of this algorithm? We will see in Casson's algorithm that it takes polynomial time to look for a 2-sphere in each convex polyhedra $A$. Since there are only $36t{t \choose 2}$ cones of the form $C_{w_{j, k},i}(M)$, it will also take polynomial time to go through all such cones to look for a 2-sphere.

There is, though, a problem in decomposing a minimal 3-manifold in irreducible pieces in polynomial time: after cutting $M$ along the non-trivial normal 2-sphere we found, we end up with 2 pieces. Unfortunately, the resulting triangulation of one of these two summands may not be minimal. We hope to find a solution to this problem. For instance, we hope to show that there exists a non-trivial normal 2-sphere with at most two tetraquads in any triangulation of a reducible closed orientable 3-manifold.

\medskip

\newpage
\pagestyle{myheadings}
\markright{  \rm \normalsize CHAPTER 5. \hspace{0.5cm}
 Minimal Triangulations}
\chapter{Normal Disks in Orientable 3-Manifolds with Nonempty Boundary}
\thispagestyle{myheadings}

\section{Collapsing Normal Disks}

\hspace{1cm} This Chapter is essentially a generalization of Chapters 3 and 4 to 3-manifolds with non-empty boundary. The idea is to cut along essential normal disks in a compact orientable 3-manifold with nonempty boundary, to collapse them to points, and to obtain an induced triangulation on the resulting pieces. All the theorems, lemmas, and corollaries are closely related to the ones in Chapters 3 and 4. We would like to point out that Theorem~\ref{MT2} shows that there are no obstructions to collapsing normal disks. Another version of this theorem has been proved by Jaco and Rubinstein ~\cite{JR:gnus} but for irreducible $\partial$-irreducible 3-manifolds only. Indeed, they proved that any triangulation of a compact orientable irreducible $\partial$-irreducible 3-manifold with nonempty boundary can be made into a triangulation where the vertices lie on $\partial M$ and there is exactly one vertex on each connected boundary component. The theorem below generalizes their results and the proof is based on the existence of $\partial$-surgery surfaces. Finally, we would like to say that, even though Proposition~\ref{Pro4} has been proved in ~\cite{JR:gnus}, the proof given here has been obtained independently.

\begin{theorem} \label{MT2}
Let M be a triangulated compact orientable irreducible 3-manifold such that each connected component of its boundary is $not$ homeomorphic to a 2-sphere. Suppose there exists a non-trivial normal disk D properly embedded in M. Then $M \cong M_{1} \#_{\partial}$ ...$\#_{\partial} M_{k} \#_{\partial} r_{1}(S^{1} \times D^{2})$, where $\sum |M_{i}| < |M|$, each $S^{1} \times D^{2}$ factor represents a 1-handle, and each $M_{i}$ is a compact orientable irreducible 3-manifold with nonempty boundary. 
\end{theorem}

\hspace{1cm} Let us define some terms that will be used in the proof of the theorem.

\begin{definition}

$\mathbf{\#_{\partial}}$ is call the boundary connected sum and it is defined as follow: let $M_1$ and $M_2$ be two triangulated compact orientable irreducible 3-manifold with non-empty boundary. Let $D_1$ (resp. $D_2$) be an embedded disk in a connected component of $\partial M_1$ (resp. $\partial M_2$). Consider an orientation reversing homeomorphism $\phi : D_1 \rightarrow D_2$. The 3-manifold $M_1 \#_{\partial} M_2$ is defined to be the quotient space $(M_1 \cup M_2)/(x \sim \phi(x))$.

A surface F with boundary is called \textbf{\textit{properly embedded}} in a 3-manifold M with nonempty boundary if F is embedded and if $F \cap \partial M = \partial F$.

A disk D, embedded in M, is called \textbf{\textit{non-trivial}} normal if it is properly embedded and if D intersects at least one tetrahedron of the triangulation of M in a quadrilateral.

A 3-manifold is called \textbf{\textit{irreducible}} if every embedded 2-sphere bounds a 3-ball.

An irreducible 3-manifold is called $\partial$-\textbf{\textit{irreducible}} if for every properly embedded disk D, $\partial D$ bounds a disk on $\partial M$. Here is an alternate definition: M is $\partial$-$irreducible$ if for every properly embedded disk D, D cuts off a 3-ball from M.

\textbf{\textit{int(M)}} will denote the interior of M. It is an open 3-manifold such that its closure is homeomorphic to M.

{\bf $M \backslash D$} will denote the resulting 3-manifold after cutting M along a normal disk D. In fact, $M \backslash D$ is homeomorphic to the 3-manifold obtained from M after removing an open regular neighborhood of D. The boundary of $M \backslash D$ consists of the boundary of M after doing a compression along D.

The triangulation of a compact orientable 3-manifold with nonempty boundary is defined in the same way as for closed 3-manifolds with the additional property that some of the tetrahedra may have vertices, edges, and/or faces belonging to $\partial M$.
\end{definition}

\medskip

\hspace{1cm} The proof of this theorem is closely related to the one of Theorem~\ref{MT} with some minor changes. We cut $M$ along $D$, and we obtain a 3-manifold $M \backslash D$ with nonempty boundary. Our goal is to collapse each copy of $D$ to a point. Note, this is topologically equivalent to cutting along a disk and leaving the new boundary as it is. We could simply cut along disks and not do any collapsing. What is interesting about the collapsing process, besides the fact that it doesn't change the topology of the manifold, is that it reduces the number of tetrahedra in the triangulation. After cutting M along $D$, we obtain a cell decomposition for $M \backslash D$ composed of the same type of polyhedra as $M \backslash S$ in Theorem~\ref{MT}. The difference is that some of the polyhedra may have vertices, edges, and/or faces belonging to $\partial M$. 

\bigskip
\begin{figure}[h]
\hspace{1.6in} \psfig{figure=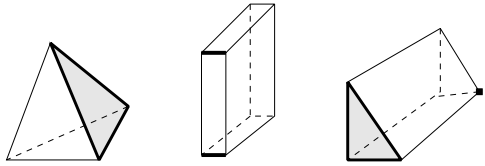, width= 2.5in} \caption{Some polyhedra in $M \backslash D$. The thickened edges and vertices, and the shaded triangles belong to $\partial M$.}
\end{figure}
\bigskip

\hspace{1cm} An irreducible 3-manifold with a 2-sphere boundary component is homeomorphic to a 3-ball. This is the reason why we assume throughout this work that $M$ has connected boundary components which are not homeomorphic to 2-spheres. Moreover, suppose that $M$ is reducible and that a connected component of $\partial M$ is homeomorphic to a 2-sphere. Suppose that $\partial D$ is the curve shown in the figure below. If we cut along this curve and collapse each copy of the curve to a point, we might not obtain a well-defined triangulation for the new boundary components of $M \backslash D$. So, a 2-sphere as boundary may constitute an obstruction to the collapsing process. This is another reason why we are not considering compact orientable 3-manifolds with boundary homeomorphic to a 2-sphere.

\bigskip
\begin{figure}[h]
\hspace{2in} \psfig{figure=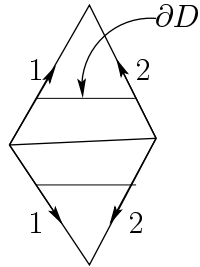, width= 1in} \caption{A triangulation a connected component of $\partial M$ homeomorphic to a 2-sphere.}
\end{figure}
\bigskip

\begin{lemma} \label{L5}
Let $M$ be a compact orientable 3-manifold with a non-empty boundary which not homeomorphic to $S^{2}$. Let D be a properly embedded disk. Let $D'$ be an embedded disk with $\partial D' = \alpha \cup \beta$, such that $D' \cap D = \alpha$ and $D' \cap \partial M = \beta$. We foliate $D'$ as in the figure below, where all the leaves are homeomorphic to arcs except one which is homeomorphic to a point. Let $D = D_{1} \cup D_{2}$, where $D_{1} \cap D_{2} = \alpha$. Consider now the branched surface $D \cup D'$. We cut along it, and we collapse each copy of D to a point and each copy of $D'$ to an edge, where collapsing means taking the quotient space formed by identifying each leaf of the foliation to a point. This is topologically equivalent to cutting along three properly embedded disjoint disks, one parallel to D, one parallel to $D' \cup D_{2}$, and one parallel to $D' \cup D_{1}$, and collapsing each copy to a point. 
\end{lemma}

\bigskip
\begin{figure}[h]
\hspace{1.5in} \psfig{figure=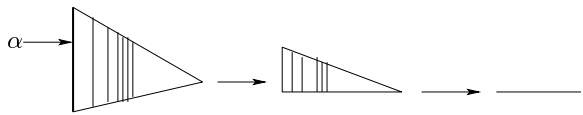, width= 3in} \caption{The collapsing of D'.}
\end{figure}
\bigskip

\hspace{1cm} \underline{\bf{Proof:}} We would first like to mention the following fact: let $M$ be as above, and let $D$ be a properly embedded disk. The 3-manifold obtained by cutting along $D$ is homeomorphic to the 3-manifold obtained by cutting along $D$ and collapsing each copy of it to a simply connected graph. The proof is essentially the same as in Lemma~\ref{L1}, and so we omit it. To simplify the notation, we assume that both $D$ and $D' \cup D_2$ are separating and that $M$ has exactly one connected boundary component.

\bigskip
\begin{figure}[h]
\hspace{1.7in} \psfig{figure=fig5.41.eps, width= 2.6in} \caption{A regular neighborhood of $D \cup D'$.}
\end{figure}
\bigskip

\hspace{1cm} Consider a regular neighborhood of $D \cup D'$. It is homeomorphic to a ball. We remove this ball from $M$. The boundary of the resulting manifold contains one parallel copy of $D$, one parallel copy of $D' \cup D_{1}$, and one parallel copy of $D' \cup D_{2}$. We now collapse $D$ to a point, $D' \cup D_{1}$ to an edge, and $D' \cup D_{2}$ to an edge, which is topologically equivalent to not collapsing anything. This proves the claim. $\Box$

\bigskip

\begin{lemma} \label{L6}
Let M be a compact orientable 3-manifold with nonempty boundary. Let D be a properly embedded disk. Let Q be an embedded quadrilateral with $\partial Q = \alpha \cup \beta \cup \alpha' \cup \beta'$, such that $Q \cap D = \beta \cup \beta'$ and $Q \cap \partial M = \alpha \cup \alpha'$. We foliate Q with arcs as in the figure below. Let $D = D_{1} \cup D'_1 \cup D_{3}$ where $D_{3} \cap D_{1} = \beta$, $D_{3} \cap D'_1 = \beta'$, and $D_{1} \cap D'_{1} = \emptyset$. Consider the branched surface $D \cup Q$. We cut along it, and we collapse each parallel copy of D to a point and each parallel copy of Q to an edge. Here, collapsing means taking the quotient space formed by identifying each leaf of the foliation to a point. The resulting 3-manifold is homeomorphic to the one obtained from M after cutting along two properly embedded disjoint disks, one parallel to D and one parallel to $Q \cup D_1 \cup D'_1$. 
\end{lemma}

\bigskip
\begin{figure}[h]
\hspace{1.5in} \psfig{figure=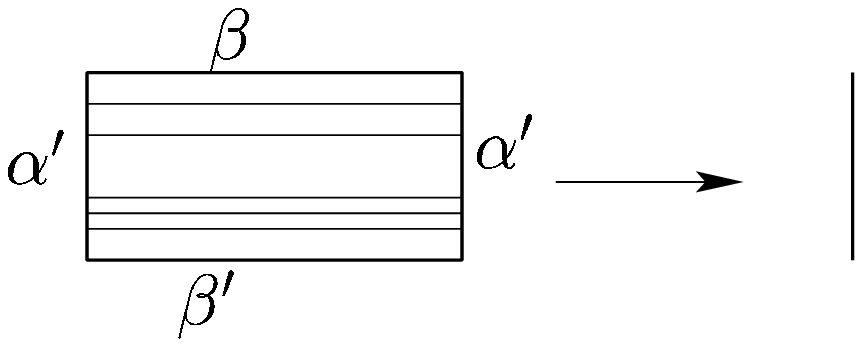, width= 1.6in} \caption{The foliation of $Q$ with leaves parallel to $\beta$ and $\beta'$.}
\end{figure}
\bigskip

\hspace{1cm} \underline{\bf{Proof:}} We need to consider two different embeddings of $Q$. Let $Nbhd(\beta)$ and $Nbhd(\beta')$ be regular neighborhoods, on $Q$, of the respective arcs $\beta$ and $\beta'$. Because $D$ is two-sided, it makes sense to talk about the two sides of $Nbhd(D)$, even if $D$ happens to be non-separating. Now, we see that $Nbhd(\partial \beta)$ and $Nbhd(\partial \beta')$ lie either on the same side of $Nbhd(D)$ or on the opposite sides of $Nbhd(D)$. 

\hspace{1cm} Case 1: $Nbhd(\partial \beta)$ and $Nbhd(\partial \beta')$ lie on the same side of $Nbhd(D)$. To simplify the notation, we suppose that $D$ and $Q \cup D_1 \cup D'_1$ are both separating. The difference with Lemma~\ref{L5} is that the neighborhood of  $D \cup Q$ is not homeomorphic to a ball but to a thickened cylinder. This is where the collapsing of $Q$ to an edge is crucial. Consider parallel copies of $\beta$ and $\beta'$ on $Q$, say $\beta_1$ and $\beta'_1$. Consider parallel copies of $D_1$ and $D'_1$, say $D_2$ and $D'_2$ respectively, such that $D_2 \cap Q = \beta_1$ and $D'_2 \cap Q = \beta'_1$. Consider the subset $Q'$ of $Q$ described in the figure below.

\bigskip
\begin{figure}[h]
\hspace{1.5in} \psfig{figure=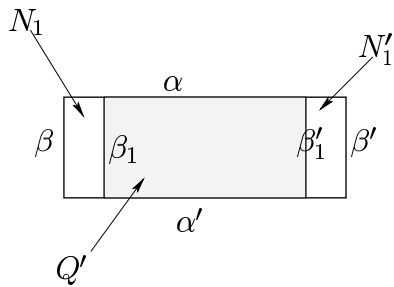, width= 2in} \caption{Partition of $Q$ into 3 subsets.}
\end{figure}
\bigskip

\hspace{1cm} Note, $Q' \cup D_2 \cup D'_2$ is parallel to $Q \cup D_1 \cup D'_1$. We cut along the two properly embedded disjoint disks $D$ and $Q' \cup D_2 \cup D'_2$. From our earlier assumption, we obtain three summands $M_1$, $M_2$, and $M_3$. We are interested in the summand $M_3$ which contains both a copy of $D$ and a copy of $Q' \cup D_2 \cup D'_2$. We collapse $D$ to a point and the other disk to an edge. Indeed, $D_2$ and $D'_2$ are collapsed to points and $Q'$ is collapsed to an edge via the foliation of $Q$. Now, we want to show that the manifold we just obtained is homeomorphic to the one obtained from $M$ after collapsing $D \cup Q$. As we are collapsing $D$ and $Q' \cup D_2 \cup D'_2$, we only need to keep track of $N_1$ and $N_2$. These two disks are collapsed to properly embedded disks as in the figure below. Note, $N_1$ and $N'_1$ now cut off 3-balls from $M_3$, so we remove those 3-balls and we collapse $N_1$ and $N'_1$ via the original foliation of $Q$ without changing the topology of the manifold. The result follows. 

\bigskip
\begin{figure}[h]
\hspace{1.5in} \psfig{figure=fig5.8.eps, width= 2in} \caption{A surgery $\partial$-annulus.}
\end{figure}
\bigskip

\bigskip
\begin{figure}[h]
\hspace{1.5in} \psfig{figure=fig5.9.eps, width= 2.5in} \caption{$M_3$ after collapsing $D$ to a point and $Q' \cup D_2 \cup D'_2$ to an edge.}
\end{figure}
\bigskip

\hspace{1cm} Case 2: $Nbhd(\partial \beta)$ and $Nbhd(\partial \beta')$ lie on opposite sides of $Nbhd(D)$. The argument is essentially the same as in Case 1 and so we will not repeat it. The only difference is that, one of the three resulting pieces is homeomorphic to a solid torus. For more details, see Claim~\ref{C5}. $\Box$

\bigskip

\section{Proof of Theorem~\ref{MT2}}

\hspace{1cm} {\bf Step 1:} We first collapse the (truncated) prisms one at a time. By the same argument as in step 1 of Theorem~\ref{MT}, if the (truncated) prism $P$ is embedded in $int(M)$, we can collapse it without changing the topology of $M \backslash D$. If it is not embedded but it is still a subset of $int(M)$, there are several cases to consider:

\hspace{1cm} {\bf{Case 1:}} None of the leaves of one side of $P$ are embedded, and all the other leaves are. This implies the existence of a surgery annulus or disk. Note, the union of the two sides may not be embedded and in this case the union of the two sides represents the surgery (annulus) disk. We cut along the surgery (annulus) disk and then collapse $P$. Using Lemmas~\ref{L1} and ~\ref{L2}, this is topologically equivalent to cutting and collapsing the disjoint union of an embedded disk $D$ and a 2-sphere S. Because $M$ is irreducible, S bounds a ball, and so we discard the trivial summand. 

\medskip

\hspace{1cm} {\bf Case 2:} Exactly one leaf from one of the sides of $P$ is not embedded. As in case 2 of Theorem~\ref{MT}, there is a surgery Mobius band $B$ . Note, the union of one leaf from one side with another leaf from the other side may not be embedded. In this case the union of the two sides represents a surgery Mobius band. Here, $\partial B \subset D$ and hence $\partial B$ bounds a disk $D'$ in $D$. $B \cup D'$ is homeomorphic to a projective plane, and because $M$ is orientable, a regular neighborhood of $B \cup D'$ is homeomorphic to a punctured $\mathbf{RP}^{3}$. Since $M$ is irreducible, this implies that $M$ is homeomorphic to $\mathbf{RP}^{3}$ which is a contradiction since $\partial M \not= \emptyset$. Therefore, case 2 cannot happen in a compact orientable irreducible 3-manifold with nonempty boundary.

\medskip

\hspace{1cm} {\bf Case 3:} None of the leaves of $P$ are embedded. This implies that the top of $P$ is identified to its bottom without any twist, and hence, $P$ is homeomorphic to a ball. Its sides represent surgery (annuli) disks so we simply cut along them. As in case 3 of Theorem~\ref{MT}, $P$ contributes to an $S^{3}$ summand which will not appear in the decomposition of $M$.

\medskip

\hspace{1cm} {\bf Case 4:} Exactly one leaf in $P$ is not embedded and it does not belong to a side. This means that the top is identified to the bottom with a 1/3 or 2/3 twist. As in case 4 of Theorem~\ref{MT}, $P$ represents a solid torus $T$ with $\partial T = A \cup A'$, where $A$ is an embedded annulus in $D$ and the boundary of $A$ represents two parallel (3,1) curves on $\partial T$. From the Jordan curve theorem, one boundary component of $A$ must bound a disk $D'$  on $D$ whose interior is disjoint from $A$. Hence, $T \cup D'$ represents a punctured $L(3,1)$ which implies that $M$ is homeomorphic to $L(3,1)$. This is a contradiction since $\partial M$ is nonempty. Hence, case 4 cannot happen in an orientable irreducible 3-manifold with nonempty boundary.

\medskip

\hspace{1cm} \underline{Remark:} If the interior of any of the sides of $P$ is embedded in $\partial M$, then case 1 may occur only if $P$ is regular. Indeed, we see that $\partial D$ bounds a disk $D'$ in $\partial M$, and because $M$ is irreducible $D \cup D'$ bounds a ball. Case 2 would imply that $\partial M$ is non-orientable which is impossible since $M$ is orientable. Also, the existence of a surgery annulus would imply that $\partial D$ is disconnected which is not the case. Now, it may happen that both endpoints of the leaves of the sides of $P$ are in $\partial M$. Hence, there are three more cases to consider.

\medskip

\hspace{1cm} {\bf Case 5:} The endpoints only, of all the leaves of one or both sides of a regular prism $P$ belong to $\partial M$:

\bigskip
\begin{figure}[h]
\hspace{1.5in} \psfig{figure=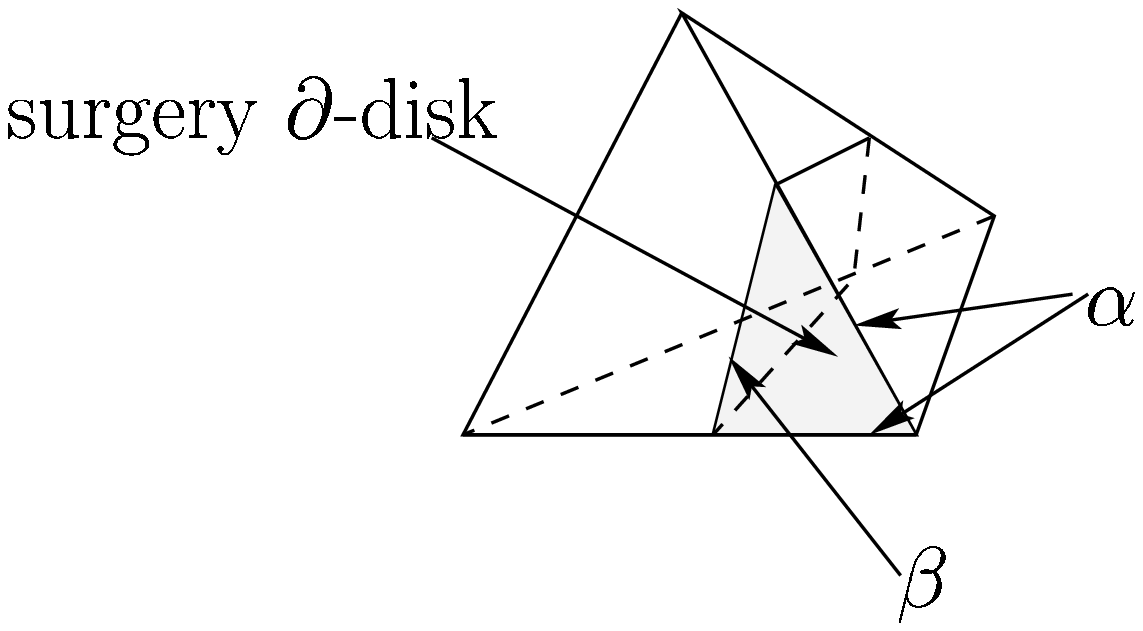, width= 2.5in} \caption{The endpoints of the leaves of one side of $P$ are in $\partial M$.}
\end{figure}
\bigskip

\hspace{1cm} The shaded side of $P$ represents a disk $D'$ with $\partial D' = \beta \cup \alpha$, where $D \cap D' = \beta$ and $\alpha$ is an arc in $\partial M$. By Lemma~\ref{L5}, cutting along $D \cup D'$ and collapsing each copy of $D$ to a point and each copy of $D'$ to an edge is topologically equivalent to cutting along two properly disjoint embedded disks (one parallel to $D$ and and one parallel to $D_{1} \cup D'$) and collapsing each copy to a point. We will call the disk $D'$ a $\partial$-\textbf{\textit{surgery disk}}, and we will say that it is \textbf{\textit{inessential}} if $D' \cup D_{1}$ or $D' \cup D_2$ (see Lemma ~\ref{L5} for notation) cuts off a 3-ball not containing $D$.

\hspace{1cm} {\bf Case 6:} The endpoints only of all the leaves of one or both sides of a truncated prism $P$ belong to $\partial M$

\bigskip
\begin{figure}[h]
\hspace{1.5in} \psfig{figure=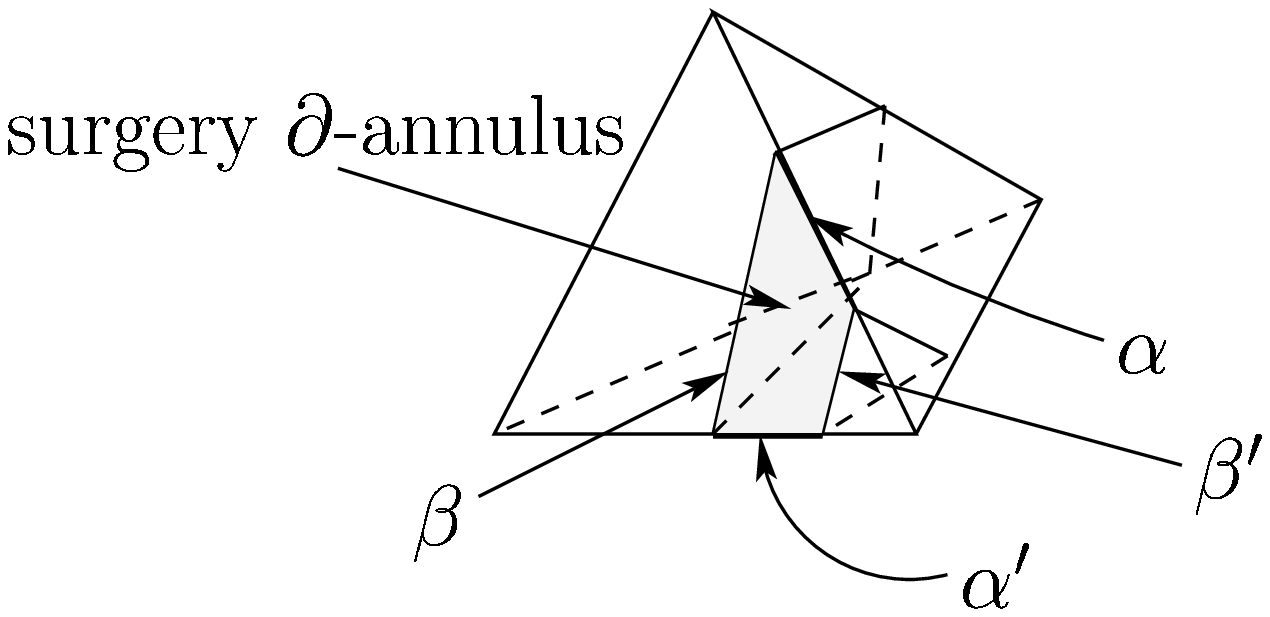, width= 2.7in} \caption{The endpoints of the leaves of one side of $P$ are in $\partial M$.}
\end{figure}
\bigskip

\hspace{1cm} The shaded side of $P$ is homeomorphic to a disk $Q$. Let $\partial Q = \beta \cup \alpha \cup \beta' \cup \alpha'$, where $\alpha$ and $\alpha'$ belong to the edges of the shaded side of $P$ belonging to $\partial M$. Because $\beta$ and $\beta'$ are disjoint embedded arcs on $D$, the Jordan curve theorem tells us that they both bound disjoint disks $D_{1}$ and $D_{2}$ in $D$.

\hspace{1cm} We cut along $Q$, and we collapse $P$. Note that the collapsing of $P$ implies the collapsing of $Q$ to an edge. By Lemma~\ref{L6}, this is topologically equivalent to cutting along two properly embedded disjoint disks, one parallel to $D$ and one parallel to $Q \cup D_{1} \cup D'_1$, and collapsing each copy of a point. We will call the disk $Q$ a $\partial$-\textbf{\textit{surgery annulus}}, and we will say that $Q$ is \textbf{\textit{inessential}} if $Q \cup D_{1} \cup D'_1$ or $D_1 \cup D_3 \cup Q \cup D_1$ (see Lemma~\ref{L6} for notation) cuts off a 3-ball whose interior does not contain $D$.

\hspace{1cm} {\bf Case 7:} There is a collection of (truncated) prisms glue along their tops and bottoms such that the bottom of the first prism and top of the last one are in $\partial M$. There are two cases to consider:

\hspace{1cm} Subcase 1: The collection of (truncated) prisms is homeomorphic to a ball $B$.

\bigskip
\begin{figure}[h]
\hspace{1.5in} \psfig{figure=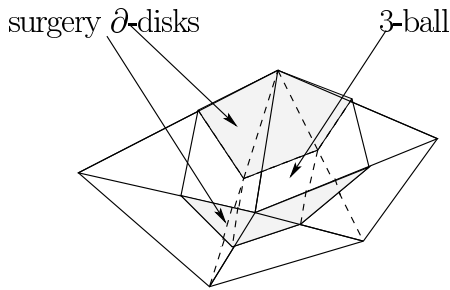, width= 2.3in} \caption{A collection of prisms homeomorphic to a 3-ball.}
\end{figure}
\bigskip

\hspace{1cm} Suppose first that the collection is made of regular prisms. Then there is a ball $B$ such that $\partial B \cap int(M) = D_{1} \cup ...\cup D_{n} \cup D'_{1} \cup ... \cup D'_{n} \cup Q_{1} \cup ... \cup Q_{n}$ where $D_{i}$ and $D'_{i}$ are sides of regular prisms, and $Q_{i}$ are quadrilaterals. Note, $D_{1} \cup ...\cup D_{n}$ and $D'_{1} \cup ... \cup D'_{n}$ both represent surgery $\partial$-disks. As in case 5, we simply collapse each of them to edges, and we discard the ball $B$.

\medskip

\hspace{1cm} Suppose now that the collection is made of truncated prisms. Then $\partial B = D_{1} \cup ...\cup D_{n} \cup D'_{1} \cup ... \cup D'_{n} \cup Q_{1} \cup ... \cup Q_{n} \cup T_{1} \cup ... \cup T_{2n}$ where the $D_{i}$'s and $D'_{i}$'s are sides of prisms, the $Q_{i}$'s are quadrilaterals, and the $T_{i}$'s are triangles. As in case 6, both $ D_{1} \cup ...\cup D_{n}$ and $D'_{1} \cup ... \cup D'_{n}$ represents surgery $\partial$-annuli and so we cut along them, we collapse each copy an edge, and we discard $B$.

\medskip

\hspace{1cm} Subcase 2: The collection of truncated prisms is homeomorphic to a solid cylinder $C$. In this case, $\partial C = A_{1} \cup B_{1} \cup A_{2} \cup B_{2} \cup A_{3} \cup B_{3}$ where, the $A_{i}$'s are $\partial$-surgery annuli, and the $B_{i}$'s are annuli consisting of quadrilaterals and triangles from $D$.

\bigskip
\begin{figure}[h]
\hspace{1.5in} \psfig{figure=fig5.18.eps, width= 2.6in} \caption{A collection of truncated prisms homeomorphic to a solid cylinder.}
\end{figure}
\vspace{.2in}

\hspace{1cm} Collapsing this collection of truncated prisms results in removing $C$ from $M$ and collapsing it to a disk whose boundary consists of three edges from the three $\partial$-surgery annuli. It seems that we created a 1-handle. We are going to see that it is not the case. Consider the intersection of three $\partial$-surgery annuli with $D$. It represents six arcs, $\alpha_{1}$, ..., $\alpha_{6}$, all embedded in $D$. Let $\alpha_{1}$ be an innermost arc in $D$, i.e. the arc which cuts off a disk $D'$ from $D$, whose interior is disjoint from $\partial C$. Now, consider the union of $D'$ with $C$. Topologically, $D' \cup C$ represents a ball, and so cutting along $D'$ and removing $C$ is topologically equivalent to cutting along $D'$ only.

\bigskip
\begin{figure}[h]
\hspace{1.5in} \psfig{figure=fig5.14.eps, width= 2.9in} \caption{Removing the solid cylinder $C$ does not create a 1-handle.}
\end{figure}
\bigskip

\hspace{1cm} {\bf Step 2:} After collapsing all the (truncated) prisms, we look at the tips. This is essentially the same argument as in step 2 of Theorem~\ref{MT}. We look at a maximal collection of tips glued together along their sides. As we have seen, a maximal collection of tips is homeomorphic to a ball, and the union of the sides of the tips of this maximal collection represents surgery disks or $\partial$-surgery disks. Hence, after collapsing the ($\partial$-) surgery surface to an edge, we either obtain an $S^{3}$ summand or a $B^{3}$ summand.

\bigskip
\begin{figure}[h]
\hspace{1.5in} \psfig{figure=fig5.20.eps, width= 2.9in} \caption{A maximal collection of tips.}
\end{figure}
\vspace{.2in}

\hspace{1cm} {\bf Step 3:} We now look at the $I$-bundles. Again, the argument is essentially the same as in Theorem~\ref{MT}. We take a maximal collection of $I$-bundles glued together along their sides. The boundary of this $I$-bundle, $\partial I \times S$, is a subsurface of the disk $D$. We have shown in Chapter 2 the possible different $I$-bundles up to homeomorphism of the total space: either the subsurface is connected, and in this case the $I$-bundle represents an punctured $\mathbf{RP}^{3}$, or it is disconnected, and in this case the $I$-bundle is trivial. Note, the former case cannot occur since $M$ is irreducible and has nonempty boundary. If none of the sides ($I \times \partial S$) of the $I$-bundle are in $\partial M$, then we have an $S^{3}$ summand. If one (or more) of the sides in the maximal collection is in $\partial M$, then we obtain a $B^{3}$ summand.

\vspace{.2in}

\hspace{1cm} {\bf Step 4:} We collapse each truncated tetrahedra to a regular tetrahedra.

\hspace{1cm} {\bf Step 5:} Note, because $M$ is irreducible, there is no need to check the existence of $S^{1} \times S^{2}$ summands. Instead, we need to count the number of $(S^1 \times D^2)$ summands which arise in the collapsing process. We first would like to prove the following claim:

\begin{claim} ~\label{C5}
Let M be as in Theorem~\ref{MT2}. The 3-manifold obtained by cutting along a non-separating disk D and collapsing each copy of D to a point is homeomorphic to a compact orientable irreducible 3-manifold $M'$ with nonempty boundary, such that $M \cong M' \#_{\partial} (S^{1} \times D^{2})$.
\end{claim}

\bigskip
\begin{figure}[h]
\hspace{1.5in} \psfig{figure=fig5.15.eps, width= 2.9in} \caption{Cutting along $D$ and $D'$ results in homeomorphic manifolds.}
\end{figure}
\bigskip

\hspace{1cm} \underline{\bf{Proof \space :}} Let $D$ be the non-separating embedded disk. Because $M$ is irreducible, $\partial D$ is non-separating in $\partial M$. Let $l$ be a non-separating closed curve in $\partial M$ which intersects $\partial D$ exactly once. Consider a regular neighborhood of $D \cup l$. It is homeomorphic to a torus $T$ whose boundary consists of a punctured torus $T'$ in $\partial M$ and a disk $D'$ in $int(M)$.

\bigskip
\begin{figure}[h]
\hspace{1.5in} \psfig{figure=fig5.17.eps, width= 3.5in} \caption{$D'$ is homotopic to $D \cup l$.}
\end{figure}
\bigskip

\hspace{1cm} So if we cut along $D'$, we obtain a 3-manifold $M'$ such that $M \cong M' \#_{\partial} (S^{1} \times D^{2})$. Now, instead of cutting along $D'$, we simply cut along $D$ and call the resulting manifold $M"$. If we can show that $M'$ and $M"$ are homeomorphic, we will have proved the claim. The closed curve $l$ on $\partial M$ now represents an arc on $\partial M"$. Well, $M'$ is obtained from $M"$ by removing a regular neighborhood of the arc $l$. Because $l$ is an arc, its neighborhood is homeomorphic to a ball. Hence, $M'$ is obtained from $M"$ after removing a ball. This proves the claim.
$\Box$

\bigskip

\hspace{1cm} Let's go back to step 5. Suppose that there are $k$ $ \partial$-surgery surfaces. Suppose also that, after the collapsing, we obtain $n$ connected pieces in the decomposition of $M$. Then there are $(k + 2 - n)$ $(S^1 \times D^2)$ summands. In Theorem~\ref{MT}, we did not need to keep track of the 2-spheres which were collapsed. This comes from the fact that the connected sum of two closed manifold does not depend on where the 3-balls are removed and how the 2-spheres boundaries are glued together. On the other hand, when we take the connected sum of two 3-manifolds $M$ and $N$ with nonempty boundaries, i.e. we glue them along disks which are themselves subsets of the boundaries, we do need to keep track on which connected boundary components are the disks. This is not a problem, though. If $D$ is a non-trivial normal disk, then the boundary of $D$ and all the surgery $\partial$-surfaces lie on the same connected boundary component of $M$. We would like to mention that when we cut along $\partial$-surgery annuli or disks, we increase by 1 or 2 the number of vertices on one of the boundary components of $M$. This comes from the fact that $\partial$-surgery surfaces are collapsed to edges. Hence, even if $|M_{1}| +$ ...$+ |M_{k}| < |M|$, the sum of the number of vertices in the boundary components of the $M_{i}'$s may be bigger than the number of vertices in $\partial M$.

\bigskip

\hspace{1cm} We now calculate the running time of the algorithm to cut $M$ along a non-trivial normal disk. Let $|M| = t$. Because there are no more than $2t$ (truncated) prisms, collapsing all of them will take linear time with respect to $t$. Step 2 and 3 don't have any running time, since they only give rise to trivial summands in the decomposition of $M$. Step 4 has a running time linear in $t$. Indeed, step 4 consists of collapsing the triangles in all the truncated tetrahedra to obtain regular tetrahedra. Since there are at most $4t$ triangles to collapse, the result follows. Step 5 requires counting the number of $\partial$-surgery surfaces and the number of summands in the decomposition of $M$. The number of such surfaces is bounded by twice the number of (truncated) prisms, and so step 5 also runs in linear time with respect to $t$. Therefore, the running time for this collapsing process is linear in the number of tetrahedra. On the other hand, to check which summand is homeomorphic to a 3-ball, one needs to run the Rubinstein-Thompson algorithm. It takes exponential time to do this verification.

\bigskip
\bigskip

\hspace{1cm} The following proposition has been proved by Jaco and Rubinstein in~\cite{JR:gnus}. We give here a proof, independently of theirs, which rely on Theorem~\ref{MT2}.

\begin{proposition} \label{Pro4}
Let M be a triangulated compact orientable irreducible $\partial$-irreducible 3-manifold such that none of its connected boundary components are homeomorphic to 2-spheres. Then we can construct a triangulation for M such that all vertices are in $\partial M$ and each connected component of $\partial M$ contains exactly one vertex.
\end{proposition}

\hspace{1cm} \underline{\bf{Proof :}} Let M be a triangulated 3-manifold as in the above statement. 

\hspace{1cm} {\bf Case 1:} Suppose first that there is exactly one vertex in each connected component of $\partial M$. Suppose also that $M$ has $n$ vertices in its interior. Let $v$ be a vertex in $\partial M$. Because $M$ is connected, we know that $v$ and the $n$ vertices are connected in the 1-skeleton, $T^{(1)}$, of $M$. Consider a maximal tree of $T^{(1)}$ containing $v$ and the $n$ vertices only (no other vertices). The boundary of a neighborhood of that tree is homeomorphic to a disk. By construction, the boundary of that disk bounds a trivial disk $D'$ in $\partial M$.

\bigskip
\begin{figure}[h]
\hspace{1.5in} \psfig{figure=fig5.21.eps, width= 2in} \caption{The boundary of a regular neighborhood of a maximal tree of $T^{(1)}$.}
\end{figure}
\bigskip

\vspace{.2in}

\hspace{1cm} We normalize the disk. In the normalization process, we may have to do some surgeries on the disk and hence, we end up with a disk $D$ and possibly some number of 2-spheres. We discard the 2-spheres which have 0-weight and the ones which are composed of triangles only (they trivially bound 3-balls). We now have $k$ normal 2-spheres $S_{1}$, ..., $S_{k}$ and a normal disk $D$. We cut along all of them. Because $M$ is irreducible, we obtain $(k+2)$ pieces (the disks is clearly separating since $M$ is $\partial$-irreducible). By construction, one piece has boundary $S_{1} \cup$ ...$\cup S_{k} \cup D$ and is homeomorphic to a $(k+1)$-punctured 3-sphere. There are also $k$ other pieces each having an $S_{i}$ as boundary. From the irreducibility of $M$, they all bound 3-balls. Finally, there is one  piece with boundary homeomorphic to $\partial M$. We now collapse the disk $D$ to a point. As in the proof of Theorem~\ref{MT2}, there may be some surgery annuli and disks and some $\partial$-annuli and disks. Because $M$ is irreducible and $\partial$-irreducible, all these surfaces and $\partial$-surfaces are inessential. After the collapsing of $D$, we obtain a manifold homeomorphic to $M$ with no vertices in its interior.

\hspace{1cm} {\bf Case 2:} Suppose now that a connected component of $\partial M$ has $n$ vertices. Let $T$ be a maximal tree of $T^{(1)}$ containing these $n$ vertices only. The boundary of a regular neighborhood of $T$ is homeomorphic to a disk. Again, we normalize the disk and we end up with possibly more than one normal disk, say $n$ of them $D_{1}$, ..., $D_{n}$. We then cut along them. Because $M$ is $\partial$-irreducible, all the disks are separating and we obtain $(n+1)$ pieces. The $\partial$-irreducibility of $M$ implies that all pieces except one are homeomorphic to 3-balls or punctured 3-balls.

\bigskip
\begin{figure}[h]
\hspace{1.5in} \psfig{figure=fig5.22.eps, width= 4in} \caption{The boundary of a regular neighborhood of a maximal tree of $T^{(1)}$, before and after normalization.}
\end{figure}
\bigskip
\vspace{.2in}

\hspace{1cm} Consider the summand not homeomorphic to a 3-ball. Note, it is in fact homeomorphic to $M$. When we collapse it, there may be some surgery annuli and disks and some $\partial$-annuli and disks which are all inessential since $M$ is irreducible and $\partial$-irreducible. After the collapsing, we obtain a manifold homeomorphic to M with possibly more or less vertices on its boundary. The point here is that the new triangulation of $M$ has strictly less tetrahedra. Hence, if there are more vertices than at the beginning we can repeat case 1 or 2 until we have the desired triangulation.
The same argument can also be used on each of the connected boundary components. $\Box$
 
\bigskip

\section{$\partial$-Connected Sums of Triangulated 3-Manifolds}

\begin{construction} \label{cons5.1}
Let $M_{1}$ and $M_{2}$ be two triangulated compact orientable irreducible 3-manifolds with nonempty boundaries and with $k_{1}$ and $k_{2}$ vertices respectively. Without loss of generality, we suppose that they have exactly one connected boundary component each. Suppose also that one of them, say $M_{1}$, has at least 2 vertices on its boundary component. We can construct the $\partial$-connected sum $M$ of $M_{1}$ and $M_{2}$ with $(|M_{1}| + |M_{2}| + 1)$ tetrahedra and $(k_{1} + k_{2} - 2)$ vertices. 

\end{construction}

\hspace{1cm} This construction is closely related to Construction~\ref{cons1}. Let $v_1$ and $v_2$ be vertices in $\partial M_{1}$ and $v'$ the  vertex in $\partial M_{2}$. We remove a normal neighborhood of $v_1$ and $v'$ and we glue the resulting 3-manifolds along $\partial Nbhd(v_1)$ and $\partial Nbhd(v')$. As in Construction~\ref{cons1}, we call $T_{1}$ (resp. $T'$) the triangulation of the disk $\partial Nbhd(v_1)$ (resp. $\partial Nbhd(v')$) and $G_{1}$ (resp. $G'$) its dual. Because $T_{1}$ and $T'$ represent disks, we conclude that $G_{1}$ and $G'$ are planar graphs such that the vertices are either 1, 2, or 3-valent. In fact, the 1-valent (reps. 2-valent) vertices come from the triangles in $T_{1}$ and $T'$ having 2 edges (reps. 1 edge) on $\partial T_{1}$ or $\partial T'$. Because $M_{1}$ has two vertices in its boundary, there exists a truncated tetrahedron with the following properties:

\bigskip
\begin{figure}[h]
\hspace{2in} \psfig{figure=fig5.19.eps, width= 1.8in} \caption{A truncated tetrahedron in $M_1$ minus a ball.}
\end{figure}
\bigskip

\hspace{1cm} Let $e$ be the above thickened edge in $T_{1}$. We insert a (truncated) prism along the shaded face of the truncated tetrahedron. This results in inserting a quadrilateral in $T_{1}$ and one vertex, $r_{1}$, in $G_{1}$. We add a copy of $G'$ along the quadrilateral, on one of the two adjacent sides of $e$. See figure below.

\bigskip
\begin{figure}[h]
\hspace{1in} \psfig{figure=fig5.24.eps, width= 3.5in} \caption{Construction of $G'_1$.}
\end{figure}
\vspace{.2in}

\hspace{1cm} Note that $G'$ represents a disk, and hence, $G'_{1}$ still represents a disk. The same procedure can be done to $G'$ to obtain a new graph $G"$. It is not hard to see that $G'_{1}$ and $G"$ are isomorphic graphs, under a color preserving isomorphism. We now glue $M_{1}$ and $M_{2}$ along $T'_{1}$ and $T"$.

\begin{construction} \label{cons5.2}
Let $M_{1}$ be a triangulated compact orientable irreducible 3-manifold with nonempty boundary and $v_{1}$ vertices. Suppose $M_{1}$ has at least 3 vertices on one of its connected boundary components F. Then we can construct a triangulation of $M \#_{\partial} (S^1 \times D^2)$ with $|M_{1}| + 1$ tetrahedra and $v_{1} - 2$ vertices.
\end{construction}

\begin{construction} \label{cons5.3}
Let $M_1$ and $M_2$ be 1-vertex triangulated compact orientable irreducible 3-manifolds with non-empty boundaries. Then we can construct a 1-vertex triangulation of $M_1 \#_{\partial} M_2$ with $|M_1| + |M_2| + 2$ tetrahedra.
\end{construction}

This construction is almost identicalto the first one. Let $f$ be a face of the triangulation of $M_1$ belonging to $\partial M_1$. We now take a tetrahedron, not belonging to the triangulation of $M_1$, and we identify one of its faces to $f$. This results in a triangulation of $M_1$ with $(|M_1| + 1)$ tetrahedra and 2 vertices. We apply Construction~\ref{cons5.1} and the result follows.

\vspace{.2in}

\hspace{1cm} Using these constructions, we prove the following:

\begin{proposition}
Let $M$ be a triangulated compact orientable irreducible 3-manifold, not necessarily $\partial$-irreducible. Then M can be made into a triangulation with no vertices in its interior and exactly one vertex in each connected boundary component.
\end{proposition}

\hspace{1cm} \underline{\bf{Proof :}} Without loss of generality, suppose that $M$ has exactly one vertex $v$ in its interior and 2 vertices, $v_{1}$ and $v_{2}$, in one of its connected boundary components. Because $M$ is connected, there is a disk $D'$ which is the boundary of a neighborhood of a maximal tree of $T^{(1)}$ containing v, $v_{1}$, and $v_{2}$ only. We normalize $D'$ and obtain a certain number of normal disks and normal 2-spheres: $D_{1}$, ..., $D_{k}$, $S_{1}$, ..., $S_{r}$. As in Proposition~\ref{Pro4}, all the 2-spheres bound 3-balls and hence, we will ignore them. We now would like to cut along the disks and collapse them to points. Following the proof of Theorem~\ref{MT2}, we obtain a decomposition of $M$ with some number of pieces say $n$ of them. Some of them may come from the cutting of disks, and some may come from the cutting of $\partial$-annuli or disks. In any case, we are left with n pieces $M_{1}$, ..., $M_{n}$, which we want to glue back together using Constructions~\ref{cons5.2} and ~\ref{cons5.3}. If some of these pieces have a vertex in their interior or more than one vertex in a connected boundary component, we reiterate this construction. The collapsing of disks strictly reduces the number of tetrahedra and hence, there can only be a finite number of iterations. We end up with, say $m$ pieces, all of which have no vertices in their interior and exactly one vertex for each connected boundary component. Note that, to apply Construction~\ref{cons5.2}, we need three vertices in one of the boundary components of one of the summands. But this is not a problem. Indeed, if $M$ has only one vertex in one of its boundary components, we can add two tetrahedra on the boundary to obtain a boundary component with three vertices. See the figure below.

\bigskip
\begin{figure}[h]
\hspace{1in} \psfig{figure=fig5.23.eps, width= 4.5in} \caption{Adding a vertex on $\partial M$ by gluing a tetrahedron on one of the triangles of $\partial M$.}
\end{figure}
\bigskip
\vspace{.2in}

\hspace{1cm} Using the above constructions, we can now glue back all the $m$ pieces together with all the $S^1 \times D^2$ summands, to obtain a 3-manifold homeomorphic to $M$ without vertices in its interior and with exactly one vertex in each of its connected boundary components. $\Box$

\bigskip
\bigskip

\section{Small Disks in Triangulated Compact Orientable 3-Manifolds with Nonempty Boundary}

\hspace{1cm} We are now ready to prove the existence of ``small" non-trivial normal disks in an irreducible 3-manifold with nonempty boundary equipped with a minimal triangulation. 

\begin{lemma} \label{L7}
Let M be a triangulated compact orientable irreducible 3-manifold with nonempty boundary. Suppose that each connected component of its boundary is \textbf{\textit{not}} homeomorphic to a 2-sphere. Let D be a non-trivial normal disk. Then there exists a non-trivial normal disk whose $\partial$-surgery surfaces are all inessential.
\end{lemma}

\hspace{1cm} Contrary to Lemma~\ref{L4}, we do not have a result about vertex normal disks. Indeed, if $D$ is a vertex normal disk, we do not know if any $\partial$-surgery annulus is inessential or not. So, to show the existence of a non-trivial normal disk $D$ with only inessential $\partial$-surgery surfaces, we are going to show that, given $D$, we can always find another disk with fewer essential surfaces.

\hspace{1cm} \underline{\bf{Proof:}} Without loss of generality, we assume that M does not contain any non-trivial normal 2-sphere. Let $D$ be a non-trivial normal disk which has the smallest number of essential $\partial$-surgery surfaces over all non-trivial normal disks. Suppose $D$ has $k$ essential $\partial$-surgery annuli and $n$ essential $\partial$-surgery disks, where $k$ and $n$ are both nonzero. Let $D'$ be one of the essential disks. $D' \cap D$ separates $D$ into two disks, $D_1$ and $D_2$. Because $D'$ is essential, neither $D' \cup D_1$ nor $D' \cup D_2$ cut off a 3-ball. We construct a disk $D"$ made of a normal parallel copy of $D_1$ and a topologically parallel copy of $D'$. See the figures below. We claim that $D"$ has at least one less $\partial$-surgery surface than $D$. Note first that the disk $D'$ is now inessential with respect to $D"$. Moreover, because $D'$ is made of triangles only, any essential $\partial$-surgery surfaces for $D"$ has its boundary on $D_1$. So any essential $\partial$-surgery surfaces for $D"$ is an essential $\partial$-surgery surface for $D$. This shows that $D"$ has less essential surfaces. This is a contradiction, and hence, $D'$ was inessential after all.

\bigskip
\begin{figure}[h]
\hspace{1.7in} \psfig{figure=fig5.25.eps, width= 2.8in} 

\vspace {1cm}
\hspace{1.6in}
\psfig{figure=fig5.26.eps, width= 1.9in} 
\caption{Construction of the disk $D"$ with no essential $\partial$-surgery surfaces.}
\end{figure}
\bigskip

\hspace{1cm} Let $A$ be an essential $\partial$-annulus for $D$. $A \cap D$ separates $D$ into two disks, $D_1$ and $D_2$, and a quadrilateral $Q$. We construct another disk $D"$ made of a topologically parallel copy of $Q$ and 2 normal parallel copies of $D_1$. See the figures below. We claim that $D"$ has less $\partial$-surgery surfaces than $D$. Note, the annulus $A$ is now inessential with respect to $D"$. Moreover, any essential $\partial$-surgery surface for $D"$ has its boundary on $Q$ or $D_1$. Hence, any such essential surface for $D"$ is also essential for $D$. This is a contradiction and hence, $A$
 was inessential after all.

\bigskip
\begin{figure}[h]
\hspace{1.5in} \psfig{figure=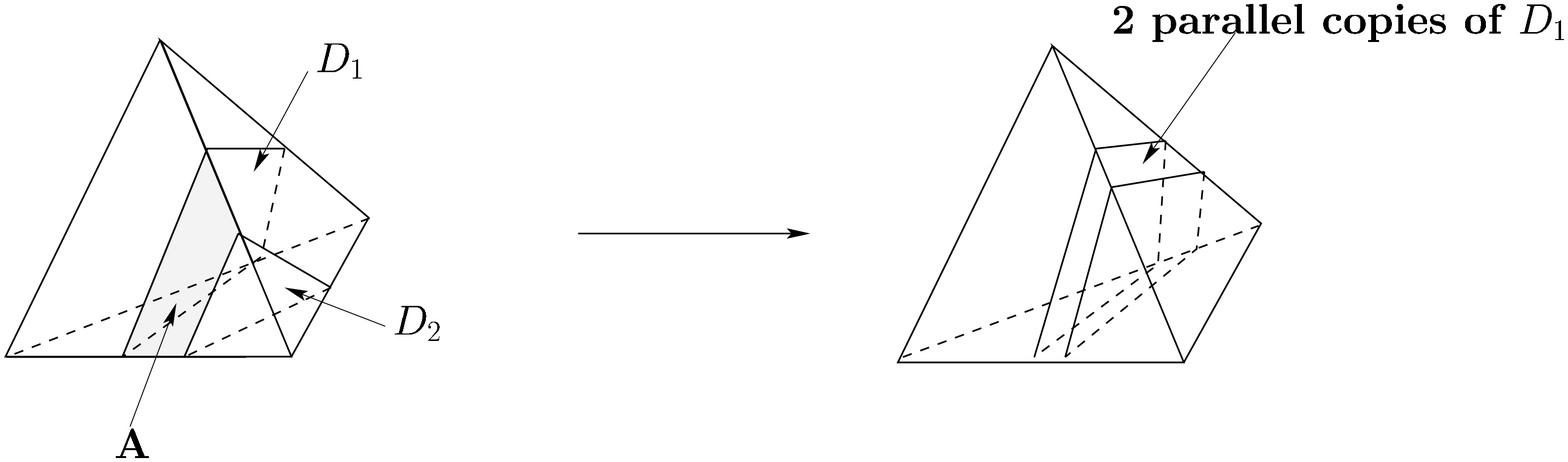, width= 4in} 

\vspace{1cm}
\hspace{1.5in}
\psfig{figure=fig5.28.eps, width= 2.5in}
\caption{Construction of the disk $D"$ with no essential $\partial$-surgery surfaces.}
\end{figure}
\bigskip
\vspace{.2in}

\medskip

\begin{theorem}
Let M be a triangulated compact orientable irreducible 3-manifold such that each connected component of its boundary is $not$ homeomorphic to a 2-sphere. Suppose that M is equipped with a minimal triangulation. Suppose M contains an essential normal disk D. Then there exists a non-trivial normal disk D' with $<$D'$> \leq 2$.
\end{theorem}

\hspace{1cm} \underline{\bf{Proof :}} Let $D$ be given and we assume it has the least number of essential $\partial$-surgery surfaces over all non-trivial normal disks. Suppose $<D> \geq 3$. If $D$ does not have any surgery annuli, surgery disks, $\partial$-disks, or $\partial$-annuli, then we cut along $D$, collapse each copy to a point, and get a decomposition of $M$ with at most two non-trivial summands (one of the summands may be homeomorphic to a 3-ball). In any case, we use the previous constructions to get a 3-manifold homeomorphic to $M$ with less tetrahedra in its triangulation. This contradicts the minimality of $M$, and so we need to assume the existence of ($\partial$) surgery surfaces.

\hspace{1cm} Because $M$ is irreducible, we disregard the existence of surgery annuli and disks. Suppose $D$ has a surgery $\partial$-annulus or disk $T$. If $T$ is essential, then by Lemma~\ref{L7} there is a non-trivial normal disk with smaller number of essential $\partial$-surgery surfaces. This contradicts our assumption. Hence, any $\partial$-annulus or $\partial$-disk must be inessential. Therefore, when we cut along $D$ and collapse each copy of it to a point, we get at most two non-trivial summands in the decomposition of $M$. Moreover, if $D$ is non-separating, then the collapsing process followed by Construction 5.2 leads us to a contradiction.

\hspace{1cm} {\bf Case 1:} Suppose $T$ is a surgery $\partial$-disk. After the collapsing, we end up with two summands. One of the two summands must have at least two vertices in one of its boundary components. This comes from the fact the collapsing a surgery $\partial$-disk creates two extra vertices. We then apply Construction 5.1 and we obtain a 3-manifold $M'$, homeomorphic to $M$, with fewer tetrahedra. Contradiction.

\hspace{1cm} {\bf Case 2:} Suppose $T$ is a surgery $\partial$-annulus. After the collapsing, we end up with two summands. It may happen that the two summands have only one vertex in each of their boundary components. In any case, we apply Construction 5.1 and we obtain a new triangulation of $M$ with less tetrahedra. Contradiction. This proves the theorem.

\bigskip

\hspace{1cm} We now apply this theorem to Casson's algorithm (see Chapter 6 for details) to check if a minimal compact orientable irreducible 3-manifold $M$ with nonempty boundary is $\partial$-irreducible or not.

\begin{corollary}
Let M be a compact orientable irreducible 3-manifold equipped with a minimal triangulation such that all the vertices lie in the boundary of M and such that there is exactly one vertex for each connected boundary component. Then it takes polynomial time (with respect to $|M|$) to check if M is $\partial$-irreducible or not.
\end{corollary}

\hspace{1cm} \underline{\bf{Proof :}} The reason this algorithm runs relatively slowly is that it has to check through all cones $D_{w, i_{1}, ..., i_{n}}(M)$, and there are $4t \cdot 3^{t}$ of them. When the triangulation is minimal though, one only needs to check through a very particular subset of this set of cones.

\hspace{1cm} First of all, if there are no non-trivial normal disks in $M$ then $M$ is $\partial$-irreducible. This is a result of Kneser (see ~\cite{Kn:gnus}) which states that if no non-trivial disks exist, then no embedded essential disks (disks which does not cut off a 3-ball) exist. Secondly, if $M$ contains a non-trivial normal disk, then $M$ contains such a disk $D$ of smallest weight. By Lemma~\ref{L7}, any $\partial$-annuli or disk is inessential. We then cut along $D$ and collapse each copy to a point, and we end up with one or two non-trivial summands. This implies that $D$ is essential, and hence $M$ is $\partial$-irreducible. Now, from Theorem~\ref{MT2}, we know there exist an essential disk with $<D> \leq 2$. Therefore, if $M$ contains an essential disk, then it must contain a non-trivial normal disk with 1 or 2 tetraquads.

\hspace{1cm} \underline {\bf{Algorithm :}}
Let $M$ be given with a triangulation with $t$ tetrahedra. 

\hspace{1cm} All the cones $D_{w, i}(M)$ are determined by the triangulation. We define $D_{w_{j}, i}(M)$  to be the set of normal surfaces such that they have at most one quadrilateral type in the $j^{th}$ tetrahedra and, to each vertex $v_{s}$, the triangle entry $i_{s}$ is zero. Consider the set $C_{w_{j}, i}(M)$ which represents all real positive solutions satisfying the normal surface equations, the quadrilateral property, such that the $i^{th}$ triangle type is zero and all entries corresponding to quadrilaterals different from $w_{j}$ are zero. This set represents a cone in $\mathbf{R}^{7t}$. Similarly, we define $D_{w_{j}, w_{k}, i}(M)$ to be the set of normal surfaces which have at most 2 tetraquads respectively. Consider the cone $C_{w_{j}, w_{k}, i}(M)$. Note, since any cone of the form $C_{w_{j}, i}(M)$ lies in a cone of the form $C_{w_{j}, w_{k}, i}(M)$, it suffices to look at the latter ones. 

\hspace{1cm} We maximize the Euler characteristic function on the convex polyhedra $A = C_{w_{j}, w_{k}, i}(M) \cap \{\sum_{i=1}^{7t} t_{i} = 1 \}$. If $\chi(S) > 0$, then there exists a disk in $C_{w_{j}, w_{k}, i}(M)$ and the procedure stops here: M is $\partial$-reducible. If $\chi(S) \leq 0$, there are no non-trivial normal disks in this cone. We repeat this step with a new cone $C_{w_{j}, w_{k}, i}(M)$.

\hspace{1cm} We will show in Chapter 6 that maximizing the Euler characteristic function in a cone takes polynomial time with respect to the number of tetrahedra. We count the number of cones of the form $C_{w_{j}, w_{k}, i}(M)$. There are ${t \choose 2}$ possible ways of choosing 2 distinct tetrahedra. Having fixed 2 tetrahedra, there are $3^{2}$ ways of choosing a quadrilateral type in each tetrahedron. For each choice of quadrilateral types, there are $4t$ ways to choose a triangle type. Hence, there are $3^{2} \cdot 4t \cdot {t \choose 2}$ different cones to consider. This shows that it takes polynomial time to look for a non-trivial disk in the cones $C_{w_{j}, w_{k}, i}(M)$.

\newpage
\pagestyle{myheadings}
\markright{  \rm \normalsize CHAPTER 6. \hspace{0.5cm}
 Minimal Triangulations}
\chapter{Decomposition of 3-Manifolds into Irreducible Pieces}
\thispagestyle{myheadings}
\section{Casson's Algorithm to Decompose a Closed Orientable 3-Manifold into Irreducible Pieces}

\hspace{1cm} In 1952, Edwin E. Moise (~\cite{Mo:gnus}) proved that any 3-manifold can be triangulated and that, given any two triangulations $K_{1}$ and $K_{2}$ of the same 3-manifold, $K_{1}$ and $K_{2}$ are equivalent, i.e. they have isomorphic subdivisions. As of today, it has actually been proven that we can obtain $K_{2}$ from $K_{1}$ through a series of Pachner moves (~\cite{Pa:gnus}). The paper by Moise was a major breakthrough in low-dimensional topology (~\cite{Mo:gnus}), and it generalized some of Kneser's results (~\cite{Kn:gnus}). Less than 10 years later, Wolfgang Haken continued to develop the theory of normal surfaces which relies entirely on the triangulability of 3-manifolds (~\cite{Hak:gnus}). 

\hspace{1cm} Let $M$ be a closed (compact without boundary) orientable 3-manifold triangulated with $t$ tetrahedra ($|M|=t$). 
We would like to mention that, even though the following lemmas were proved by the author, the actual results are due to Andrew Casson. 

\medskip

\begin{definition}
A normal surface is called  \textbf{\textit{trivial}} if it intersects the tetrahedra in triangles only. Note, a trivial surface must be a union of spheres, which must be themselves boundaries of regular neighborhoods of vertices.
\end{definition}

\bigskip
Throughout the algorithm, we will make use of Theorem~\ref{MT} and Proposition~\ref{Pro1}.


\medskip


\medskip

\begin{definition}

Let \emph{$\overline{N}$}(M) be the set of surfaces in M which intersect each tetrahedron in triangles or quadrilaterals. Precisely, \emph{$\overline{N}$}(M) is the set of surfaces which satisfy the normal surface equations but may not satisfy the quadrilateral property (i.e. they may not be embedded).

A type w is a function which assigns, to each tetrahedron, one of the 3 possible types of quadrilateral.

Let N(M) be the set of normal surfaces.

Let $N_{w}$(M) be the set of normal surfaces of type w.
\end{definition}

It is clear that $N(M)=  \bigcup_{w}  N_{w}(M)$.

\hspace{1cm} The reason we look at $\bigcup_{w} N_{w}(M)$ instead of $N(M)$ is because the solution space of $N(M)$ does \underline{not} form a cone in $\mathbf{R}^{7t}$ whereas each $N_{w}(M)$ does form a cone. For example, if $S_{1}$, $S_{2} \in N(M)$ are not of the same type, their normal sum may not represent an embedded surface, and hence may not even lie in $N(M)$. One could note that $S_{1} + S_{2}$ lies in $\overline{N}(M)$, but the problem is that we don't know how to decompose this space into cones.

\bigskip
\begin{claim} ~\label{C6.3}
$\exists  S^{2} \in N(M) \Leftrightarrow \exists S \in N_{w}(M)$  for some $w$, with $\chi(S)>0$.
\end{claim}

\underline{\bf{Proof:}} $\Rightarrow$ Trivial since $S^{2}$ is embedded.

\hspace{.4in} \space $\Leftarrow$ \space If \space $\exists S \in N_{w}(M)$ for some $w$, then this surface must be embedded. Since $\chi(S)>0$, one of the connected components of $S$ must be homeomorphic to either a $\mathbf{RP}^{2}$ or a 2-sphere. If it is a 2-sphere, then  we are done. If it is a $\mathbf{RP}^{2}$, then since $M$ is orientable a regular neighborhood of this projective plane is homeomorphic to a twisted $I$-bundle over it. This twisted $I$-bundle is homeomorphic to a punctured $\mathbf{RP}^{3}$. Hence, $\partial N(RP^{2})\cong S^{2}$ which has the desired property.

\medskip

Suppose M has a unique vertex $v$. Let $t_{1}, t_{2}, ..., t_{4t}$ be the triangle coordinates corresponding to $v$.
Let $N_{w, i}(M)=\{ S \in N_{w}(M) | t_{i}=0 \}$ 

\bigskip

\begin{claim} ~\label{C6.4}
$\exists$ a non-trivial $S^{2}$ $\Leftrightarrow$ $\exists S \in N_{w, i}(M)$ for some w, $i$, and with $\chi (S)>0$.
\end{claim}

\underline{\bf{Proof:}} $\Rightarrow$ If $S^{2}$ is non-trivial, then it cannot be the boundary of a neighborhood of a vertex and hence, we must have $t_{i}=0$ for some $i$.

\hspace{.4in} \space $\Leftarrow$ if $S \in N_{w, i}(M)$ then $S$ cannot be the boundary of a vertex neighborhood and hence, cannot be trivial. As in claim~\ref{C6.3}, $\chi (S)>0$ implies the existence of a non-trivial $S^{2}$.

\bigskip
\hspace{1cm} We define $\mathbf{C_{w, i}(M)}$ to be the solution space to the surface equations with $t_{i}=0$, and the rest of the variables being \emph{real} non-negative.

\bigskip
\begin{claim} ~\label{C6.5}
$\exists S \in N_{w, i}(M)$ for some $w, i$ with $\chi (S)>0$ \space $\Leftrightarrow$ \space $\exists V \in C_{w, i}(M)$ for some $w , i $, with $\chi (V)>0$, $V$ is on an external ray of the cone $C_{w,i}(M)$, and $V \in \mathbf{Z}^{7t}$.
\end{claim}

\underline{\bf{Proof :}} $\Leftarrow$ Suppose $\exists V \in C_{w, i}(M)$, $\chi (V)>0$, and $V \in \mathbf{Z}^{7t}$. Then , trivially, $V \in N_{w, i}(M)$.

\hspace{.4in} \space $\Rightarrow$ If $S \in N_{w, i}(M)$, then $S \in C_{w, i}(M)$. Suppose now that $S$ is not on some external ray, i.e. suppose that $S= V_{1} + V_{2} + ... + V_{n}$ where $V_{j} \in C_{w, j}(M)$, and the $V_{j}$'s lie on an external ray. Since $S$ has integer entries, we can assume that the $V_{j}$'s have rational entries. Then $\chi (S) = \chi (V_{1} + ... + V_{n})= \chi (V_{1}) +... + \chi (V_{n}) >0$. This implies that $\chi (V_{j}) >0$ for some $j$. We multiply each entry of $V_{j}$ be the least common multiple of the denominator of the entries. Since $\chi (V_{j}) >0$, we conclude that either $V_{j}$ is homeomorphic to a 2-sphere, a projective plane, or it is disconnected. The disconnectedness of the surface would contradict the fact that it lies on an external ray. If it is a projective plane, then we double each entry in the vector representation of the surface to obtain a 2-sphere on the same external ray. Also, $V_{j}$ cannot represent a trivial 2-sphere because $i = 0$ for some $i$. $V_{j}$ or its double gives us the desired surface. $\Box$

\hspace{1cm} Remark: let us define the Euler characteristic function on a vector with rational entries, satisfying the surface equations and the quadrilateral property. Let $S$ be an integer solution of the surface equations representing an embedded surface ( $S=$($a_{1}, a_{2}$, ..., $a_{7t}$) ). Let $a_{1}$, ..., $a_{4t}$ and $a_{4t+1}$, ..., $a_{7t}$ denote the coefficient of the triangle and quadrilateral entries respectively. We then define $\chi$ in the following way: $\chi (S) = \lbrack \sum_{i=1}^{4t}a_{i} \cdot (\frac{\rm 1}{\rm d_{i_{1}}} + \frac{\rm 1}{\rm d_{i_{2}}} + \frac{\rm 1}{\rm d_{i_{3}}}) - \sum_{i=1}^{4t}a_{i} \cdot (\frac{\rm 3}{\rm 2}) + \sum_{i=1}^{4t}a_{i} \rbrack + \lbrack \sum_{i=4t+1}^{7t}a_{i} \cdot (\frac{\rm 1}{\rm e_{i_{1}}} + \frac{\rm 1}{\rm e_{i_{2}}} + \frac{\rm 1}{\rm e_{i_{3}}} + \frac{\rm 1}{\rm e_{i_{4}}}) - \sum_{i=4t+1}^{7t}a_{i} \cdot (\frac{\rm 4}{\rm 2}) + \sum_{i=4t+1}^{7t}a_{i} \rbrack = \lbrack \sum_{i=1}^{4t}a_{i} \cdot (\frac{\rm 1}{\rm d_{i_{1}}} + \frac{\rm 1}{\rm d_{i_{2}}} + \frac{\rm 1}{\rm d_{i_{3}}} - \frac{\rm 1}{\rm 2}) \rbrack + \lbrack \sum_{i=4t+1}^{7t}a_{i} \cdot (\frac{\rm 1}{\rm e_{i_{1}}} + \frac{\rm 1}{\rm e_{i_{2}}} + \frac{\rm 1}{\rm e_{i_{3}}} + \frac{\rm 1}{\rm e_{i_{4}}} - 1) \rbrack$, where the $d_{i}$'s and the $e_{i}$'s represent the degrees of the edges touching the triangle $i$ and the quadrilateral $i$ respectively.

\hspace{1cm} Now let $V$ be a rational solution of the surface equations, satisfying also the quadrilateral property ( $V$= ($\frac{\rm a_{1}}{\rm b_{1}}, \frac{\rm a_{2}}{\rm b_{2}}$, ..., $\frac{\rm a_{7t}}{\rm b_{7t}}) $). Let $N$= $\Pi_{i=1}^{7t}b_{i}$.
We then define $\chi$ as follows: $\chi (V) = \frac{\rm 1}{\rm N}\chi ( N \cdot V)$, where $N \cdot V$ denotes the integer solution $(\frac{\rm a_{1} \cdot N}{\rm b_{1}}, \frac{\rm a_{2} \cdot N}{\rm b_{2}}$, ..., $\frac{\rm a_{7t} \cdot N}{\rm b_{7t}})$. From this definition, it follows that $\chi$ is a linear function. Indeed, let $V_{1}$ and $V_{2}$ be two rational solutions. Then, $\chi (V_{1}) + \chi (V_{2}) = \frac{\rm 1}{\rm MN}  \chi (\frac{\rm 1}{\rm MN}(V_{1} + V_{2})) = \frac{\rm 1}{\rm MN} \chi (MN\cdot V_{1}) + \frac{\rm 1}{\rm MN}  \chi (MN\cdot V_{2}) = \frac{\rm 1}{\rm M}  \chi (M\cdot V_{1}) + \frac{\rm 1}{\rm N} \chi (N\cdot V_{2}) = \chi (V_{1}) + \chi (V_{2})$.

\vspace{.2in}

\hspace{1cm} {\bf Remark 1:} Note that, in the previous two claims, the spaces $N_{w, i}(M)$ are defined for 1-vertex triangulations only. In fact, the procedure to decompose $M$ into irreducible pieces will always be applied to 1-vertex triangulations. We now would like to explain why we are making such an assumption about the triangulation of $M$. Let $M$ be a closed orientable 3-manifold with more than one vertex. We take the boundary $S$ of a regular neighborhood of a maximal tree of the 1-skeleton. We normalize $S$. In the normalization process, we may have to do some compressions on $S$ (on some faces of some tetrahedra) and hence, we may end up with some number of normal 2-spheres (after discarding the 0-weight ones). We also discard the vertex-linking ones since they trivially bound balls. If there are no normal 2-spheres left, then $M$ is homeomorphic to the 3-sphere. So suppose we have, say $k$ non-trivial normal 2-spheres $S_{1}$, ..., $S_{k}$, where $k\geq 1$. We cut along each one of them, and we use Theorem~\ref{MT} to collapse all the $2k$ 2-spheres.

\hspace{1cm} Note, since the boundary components ($\beta_{1}$ and $\beta_{2}$) of a surgery annulus may be embedded in distinct 2-spheres, the use of case 5 of Lemma~\ref{L2} here is essential. After collapsing the $2k$ 2-spheres, we obtain say $n$ ($n\geq 1$) distinct pieces $M_{1}$, ..., $M_{n}$. First of all, if $M_{i}$ is homeomorphic to $\mathbf{RP}^{3}$, $S^{3}$, or $L(3,1)$, then no further decomposition is possible (and necessary). If $M_{i}$ has a 1-vertex triangulation, we stop here. If $M_{i}$ has more than one vertex, we go through the same process again. The process must terminate since, by Theorem~\ref{MT}, it strictly decreases the number of tetrahedra in $M_{i}$. Hence, we end up with a decomposition of $M$, $M \cong M_{1} \#$ ... $\# M_{n}$, where each $M_{i}$ has a 1-vertex triangulation (but $M_i$ does not have to be irreducible of course) or is homeomorphic to $\mathbf{RP}^{3}$, $S^{1} \times S^{2}$, $S^{3}$, or $L(3,1)$. 

\hspace{1cm} Now, the time needed to construct the normal 2-spheres $S_1$, ..., $S_k$ is linear in the number of tetrahedra. Indeed, consider the regular neighborhood $S$ of a maximal tree of $T^{(0)}$. The weight of $S$ is bounded by twice the number of edges. Hence, in the normalizing process of $S$, there could not be more than $6t$ isotopies (3 for each face).

\hspace{1cm} Therefore, the normalization process for such a surface takes linear time in the number of tetrahedra. As we saw at the end of the proof of the Theorem~\ref{MT}, collapsing a normal 2-sphere takes time polynomial in the number of tetrahedra. So if we start with a closed orientable 3-manifold $M$ with more than 1 vertex in its triangulation, we can find a decomposition of $M$ ($M \cong M_{1} \#$ ... $\# M_{n}$ with $n \geq 1$) in polynomial time where each summand has a 1-vertex triangulation or is homeomorphic to $\mathbf{RP}^{3}$, $S^{1} \times S^{2}$, $S^{3}$, or $L(3,1)$.

\vspace{.2in}

\hspace{1cm} \underline{\underline {\bf{The algorithm :}}}

\hspace{1cm} Let $M$ be a closed orientable 3-manifold equipped with a triangulation. By the above remark, we may assume without loss of generality that $M$ is a 1-vertex triangulation with $t$ tetrahedra.

\hspace{1cm} The set of surface equations is completely determined by the triangulation, and hence the cones $C_{w}(M)$ are also determined. There are $3^{t}$ of them, one for each type $w$.

\hspace{1cm} The procedure has to find non-trivial 2-spheres.

\hspace{1cm} We now use claim~\ref{C6.5} to find a non-trivial normal 2-sphere. Given a cone $C_{w, i}(M)$, we look at the set $A= C_{w, i}(M)$ \space  $\bigcap$ \space  $\{ \sum_{i=1}^{7t} t_{i} = 1  \} $. This set is a convex compact polyhedron called the projective solution space with vertices having rational entries (this comes from the fact that the equations involved all have integer coefficients). On this polyhedron, we will maximize the Euler characteristic function defined earlier. 

\hspace{1cm} Linear programming theory tells us that the maximum will be attained at a vertex of this polyhedron and hence, on an external ray of the cone $C_{w, i}(M)$. There are several methods to find such a maximum. See for instance ~\cite{Sch:gnus} Theorem 15.3 page 198. This theorem proves the existence of a method to find the maximum of a linear function on a convex polyhedron with a running time polynomial in the size of a matrix $B$. In our context, $B$ describes the surface equations and the non-negativity of the entries $t_{i}$. The size of $B$ (as defined in Schrijver page 29) happens to be a polynomial with respect to $t$. In fact, $size(B) \leq 182t^{2}$. This upper bound comes from the fact that B is a $(7t) \times (6t + 7t)$ matrix with integer entries smaller than 2. This theorem by Schrijver gives us the desired maximum in polynomial time.

\hspace{1cm} Once again, this vertex will have rational entries, so multiply it by the least common multiple of the denominator of the entries. This new vector (or surface) $S$ with integer entries could have the following Euler characteristic : 

\hspace{1cm} $\chi(S)=1$: $S$ represents a projective plane. We look at the surface $2S$ which must represent a non-trivial 2-sphere. Indeed, if it were trivial, the projective plane $S$ would be the boundary of a normal neighborhood of the vertex of the triangulation which is clearly impossible.

\hspace{1cm} $\chi (S)=2$: $S$ is a non-trivial 2-sphere.

\hspace{1cm} $\chi (S)<0$ : by claim~\ref{C6.4} and claim~\ref{C6.5}, there are no non-trivial $S^{2}$ and the procedure stops here.  

\hspace{1cm} $\chi (S)>2$: this would imply that $S$ is not connected which contradicts the fact that it is on an external ray of $C_{w, i}(M)$.

\hspace{1cm} The procedure cuts along either $S$ (if $\chi (S)=2$) or along $2S$ (if $\chi (S)=1$) and then collapse it using the Theorem~\ref{MT}. 

\medskip
\bigskip

\hspace{1cm} We now need to check the following facts:

\hspace{1cm} \underline{\bf{Fact 1:}} The procedure terminates after cutting along and collapsing finitely many non-trivial 2-spheres:
\smallskip
 
\hspace{1cm} This follows from Theorem~\ref{MT} since cutting along a non-trivial $S^{2}$ strictly reduces the number of tetrahedra in $M$.

\bigskip
\medskip

\hspace{1cm} \underline{\bf{Fact 2:}} If there are no non-trivial $S^{2}$, the procedure terminates:
\smallskip

\hspace{1cm} When maximizing $\chi$ on the polyhedron $A$, if we obtain $\chi (S) \leq 0$ we know, by claim~\ref{C6.5}, that there are no non-trivial $S^{2}$ in $C_{w,i}(M)$.

\bigskip

\bigskip

\hspace{1cm} The procedure needs to perform one more task. After cutting and collapsing along all non-trivial $S^{2}$'s, it needs to check if some of the resulting pieces are homeomorphic to $S^{3}$. To do that, we use the Thompson-Rubinstein theorem which states that one of the resulting pieces is homeomorphic to $S^{3}$ if and only if there exists an almost normal 2-sphere in it. 

\hspace{1cm} Consider a summand $M_{i}$ of $M$. 
First, we check that $H_{1}(M_{i};\mathbf{Z}_{2})=0$. If this is not the case, we know that $M_{i}$ is not homeomorphic to $S^{3}$. This can be done in polynomial time using cell decomposition:

\hspace{1cm} Let $C_{0}(M_i), C_{1}(M_i), C_{2}(M_i)$ be the 0, 1, and 2-chain complexes of $M_i$. Since $M_i$ has a 1-vertex triangulation, $C_{0}(M_i) \cong \mathbf{Z}$, and every edge in the 1-skeleton represents a loop. Hence, $\partial_{1}: C_{1}(M_i) \rightarrow C_{0}(M_i)$ is the zero map, and $Ker \partial_{1} \cong \mathbf{Z^{n}}$ where $n$ is the number of edges in the 1-skeleton. We now orient each face of the 2-skeleton by orienting each edge of the 1-skeleton. After orienting the $2t$ faces, we can find $Im \partial _{2}$. We now have $2t$ linear combinations of the edges of the 1-skeleton and we only need to check if those combinations, viewed as elements of $\mathbf{Z^{n}}$, form a basis for $\mathbf{Z^{n}}$. These $2t$ linear combinations can be viewed as a $(2t) \times n$ matrix $C$. Our problem is now equivalent to finding the rank of $C$. This can be done in polynomial time with respect to $size(C)$ (in our problem, $size(C) \leq 2(2t)(6t)$). See ~\cite{Sch:gnus}, Corollary 3.3a, page 33. If $rank(C) < n$, then $M_i$ is not homeomorphic to $S^{3}$ and the procedure stops. If $rank(C) \geq n$, then $H_1(M_i;\mathbf{Z}_2)=0$ and we need to look for an almost normal 2-sphere.

\medskip
\begin{definition}
$S \subset M$ is called \emph{almost $^{2}$ normal} (resp. $almost$ $normal$) if it is normal and if there exists at most one tetrahedron $\tau$ in which S intersects $\tau$ in triangles and octagons (resp. one octagon) only.

\hspace{1cm} Define $A_{w, i, l}(M)$ to be the set of surfaces of type $w$, with $t_{i}=0$, and one type of octagon being allowed only in the $l^{th}$ tetrahedron.
\end{definition}

\bigskip
\begin{claim} ~\label{C6.7}
$\exists$ an almost normal $S^{2}$ $\Leftrightarrow$ $\exists S \in A_{w, i, l}(M)$ for some $w$, $i, l$ with $\chi (S) - o(S) >0$, where o(S) is the number of octagons in S in the $l^{th}$ tetrahedron.
\end{claim}

\hspace{1cm} \underline{\bf{Proof :}} $\Rightarrow$ If there is an almost normal $S^{2}$, then clearly $S \in A_{w, i, l}(M)$ for some $w, i, l$. Since $S$ is almost normal, then $o(S)=1$. This implies that $\chi (S) - o(S)=1 >0$.

\hspace{.4in} \space $\Leftarrow$ If there is a $S$ in $A_{w, i, l}(M)$ with $\chi (S)- o(S) >0$, then, by definition, $S$ is almost$^{2}$ normal. Since $\chi (S) - o(S) >0$, and o(S)$\geq 1$ (if o(S)=0 then $S$ is normal non-trivial, but we already have cut along all such surfaces), then $\chi (S) > 1$. If $\chi (S) = 2$, then $o(S) = 1$ and we have an almost normal $S^{2}$. If $\chi (S) >2$, then $S$ must have an $S^{2}$ component with $\chi (S^{2}) - o(S^{2}) > 0$. This again implies that $o(S^{2}) = 1$ and we have an almost normal $S^{2}$.

\bigskip

\hspace{1cm} As before, the procedure maximizes the linear function $\chi (\bullet) - o(\bullet)$ (it is linear in the cones $A_{w, i, l}(M)$ for fixed $w$ and $l$ only) over the set $A_{w, i, l}(M) \bigcap$ $\{ \sum_{i=1}^{7t} t_{i} = 1 \}$. This gives us a solution on a rational vertex $W$, which in turn gives us an integer vector $V$. If we obtain $\chi (V) - o(V) \leq 0$, then by claim~\ref{C6.7} there are no almost normal 2-spheres and $M_{i}$ is not homeomorphic to $S^{3}$ and the procedure stops. If $\chi (V) - o(V) > 1$, then $\chi(V) > 2$ which means that $V$ is disconnected which is impossible since $V$ is a vertex. If $\chi (V) - o(V) = 1$, then $M_{i}$ is homeomorphic to $S^{3}$.

\bigskip

\hspace{1cm} \underline{\bf{Fact 3:}} The procedure cannot be looking indefinitely for an almost normal $S^{2}$ if there are none:
\smallskip

\hspace{1cm} As in the previous paragraph, if $\chi (V) - o(V)\leq 0$ we know, by claim~\ref{C6.7}, that there are no almost normal $S^{2}$.

\bigskip
\bigskip
\bigskip

\hspace{1cm} We summarize the steps needed in the algorithm to decompose a closed orientable 3-manifold into irreducible pieces.

\hspace{1cm} Let $M$ be given by a $t$-tetrahedra and $v$-vertex triangulation. 

\hspace{1cm} \underline{\bf{Step 1:}} Construct a (possibly disconnected) normal surface $S$ obtained by normalizing the boundary of a regular neighborhood of a maximal tree of the 1-skeleton. We collapse $S$ using Theorem~\ref{MT}. This may result in a decomposition of $M$: $M$  $\cong  M_{1}$ $\# M_{2}$...$\# M_{k}$ $\# r_{1}(S^{1}\times S^{2})$ $\# r_{2}\mathbf {RP}^{3}$ $\#r_{3}L(3,1)$. For each summand not homeomorphic to $\mathbf{RP}^{3}$, $S^{1} \times S^{2}$, $S^{3}$, or $L(3,1)$, and having more than one vertex in its triangulation repeat step 1. We now have a decomposition of $M$ where each summand is either $\mathbf{RP}^{3}$, $S^{1} \times S^{2}$, $L(3,1)$, or it has a 1-vertex triangulation. Let us call $M$ one of the 1-vertex summands. We will repeat the procedure for each of the other 1-vertex summands.

\hspace{1cm} \underline{\bf{Step 2 :}} Go through each cone $C_{w,i}(M)$ to find a non-trivial 2-sphere. If none are found, go to step 3. If one is found, then cut along it and collapse it using Theorem~\ref{MT}. We may obtain a further decomposition of $M$ and if so, we go through step 1 again for each of the new summands containing more than 1 vertex.

\hspace{1cm} \underline{\bf{Step 3:}} For each of the summands found in $M$ (excluding $\mathbf{RP}^{3}$, $S^{1} \times S^{2}$, and $L(3,1)$) go through the cones $A_{w,i,l}(M)$ to find almost$^2$ normal 2-spheres. If none are found, go to step 4. If one is found, then we know the summand is homeomorphic to $S^{3}$.

\hspace{1cm} \underline{\bf{Step 4:}} We obtain a decomposition of $M$ where each summand is irreducible with a 1-vertex triangulation and not homeomorphic to $S^{3}$, or is homeomorphic to $\mathbf{RP}^{3}$, $S^{1} \times S^{2}$, or $L(3,1)$.

\bigskip
\bigskip
\bigskip

\underline{\underline{\bf{Complexity of the algorithm:}}}
 
\noindent
\begin{enumerate} 
\item[1.] How long does it take for the procedure to find a non-trivial normal $S^{2}$?

\item[2.] How long does it take for the procedure to cut along and collapse a non-trivial normal $S^{2}$?

\item[3.] How many non-trivial normal 2-spheres can there be?

\item[4.] How long does it take for the procedure to look for an almost normal $S^{2}$?
\end{enumerate}

\noindent
\begin{enumerate}
\item[1.] First of all, if $M$ has more than one vertex, we can normalize the boundary of a regular neighborhood of a maximal tree, and collapse it to a point. We have already explained that this normalizing and collapsing process takes linear time in the number of tetrahedra. After this process, we may end up with a decomposition of $M$ into pieces which may not be necessarily irreducible. If some of these pieces have more than one vertex, we iterate again. We end up with a decomposition of $M$ into pieces which are either homeomorphic to $\mathbf{RP}^{3}$ or $L(3, 1)$, or which have 1-vertex in their triangulation. Each iteration is linear in the number $t$ of tetrahedra. Because we collapse a non-trivial normal 2-sphere in each iteration, there could not be more than $t$ of them. We call $M$, one of the summands of $M$ and we continue the algorithm with this piece. The procedure has to check through all cones $C_{w, i}(M)$. There are at most $3^{t} \cdot 4t$ of them.
The procedure then maximizes $\chi$ to find our desired vertex solution $S$. This can be done in O($t^{n}$) time, for some fixed integer $n$, independent of $t$.

\item[2.] We use Theorem~\ref{MT} to cut along $S$ and collapse it. It is shown at the end of the theorem that this takes polynomial time in the number of tetrahedra.

\item[3.] Since cutting along non-trivial $S^{2}$'s strictly reduces the number of tetrahedra, there could not be more than $t$ of them.

\item[4.] The procedure now has to check through each cone $A_{w, i,l}(M)$. There are at most $3^{t}\cdot 4t \cdot t$ of them. Again, we know there exists an algorithm (~\cite{Sch:gnus}) to find our desired vertex surface in $O(t^{n})$ time for some fixed $n$.
\end{enumerate}

\medskip

\indent
\hspace{1cm} This algorithm will decompose a closed orientable 3-manifold into irreducible pieces, none of them being homeomorphic to $S^{3}$. The procedure will be done in ($3^{t} \cdot O(t^{n})$) time, for some $n$ independent of $t$. 

\medskip

\hspace{1cm} Note: in Corollary~\ref{C3}, we see that if $M$ is given by a minimal triangulation, it takes polynomial time to check if it is reducible or not. On the other hand, we do not know how long it takes to decompose a minimal triangulation into irreducible pieces. Indeed, let $M$ be minimal and reducible, and let $S$ be a non-trivial normal 2-sphere $S$ with $<S> \leq 2$. After cutting along $S$ and collapsing it to a point, we have the following decomposition: $M \cong M_{1} \# M_{2}$. The problem which arises is that one of the two resulting summands, say $M_{1}$, may have a 2-vertex triangulation which may not be minimal. If this is the case, Theorem~\ref{MT3} cannot be applied to $M_{1}$.

\bigskip
\bigskip
\bigskip
\bigskip
\bigskip

\section{Algorithm to cut an irreducible compact orientable 3-manifold with nonempty boundary into irreducible $\partial$-irreducible pieces.}

\bigskip
\bigskip
\medskip

\hspace{1cm} This algorithm represents a natural generalization of Casson's algorithm. The idea is to cut an irreducible compact orientable 3-manifold $M$ with non-empty boundary along disks. As in Theorem~\ref{MT2}, we will be assuming that no connected component of the boundary of $M$ is homeomorphic to a 2-sphere (indeed, if this were the case, then $M$ would be homeomorphic to a ball). The goal here is to look for non-trivial normal disks and to cut along them to simplify the manifold. As a generalization of Casson's algorithm, we will look at cones similar to $C_{w,i}(M)$ to find such normal disks.

\medskip

\hspace{1cm} Suppose that $M$ has no vertices in its interior, that it has exactly one vertex on each connected boundary component, and that it does not contain any non-trivial normal 2-spheres. Let $v_{1}$, ..., $v_{k}$ be the vertices of the triangulation and $t$ the number of tetrahedra.

\begin{definition}
A disk D in M is called \textbf{\textit{trivial}}  if $\partial D$ consists of edges such that each of these edges comes from a triangle. D will be called \textbf{\textit{non-trivial}} otherwise. Any normal disk will be assumed to be properly embedded.

Fix a vertex of the triangulation, say $v_1$. Let $D_{w,i}(M)$ be the set of normal surfaces in M of type $w$ with $t_{i}= 0$, where $t_{i}$ corresponds to one of the triangle types around the vertex $v_1$. Moreover, we require the surfaces in $D_{w,i}(M)$ not to have any triangle types around a vertex disjoint from $v_1$.

\end{definition}

\hspace{1cm} The idea here is that we are looking for non-trivial normal disks. Because the boundary of a disk is connected, we only need to have non-zero entries for the triangle and quadrilateral types belonging to the same connected boundary component of $M$.

\medskip
\begin{claim} ~\label{C6.9}
$\exists$ a non-trivial disk D in M \space $\Leftrightarrow$ \space $\exists S \in D_{w, j}(M)$ with $\chi (S) >0$ for some w and j.
\end{claim}

\underline{\bf{Proof :}} $\Rightarrow$ If $D$ is non-trivial, then clearly $t_{i} =0$ for some $i$, where $t_{i}$ corresponds to one of the triangle types around the vertex $v_1$. This implies that $D \in D_{w, j}(M)$ with $\chi (D) >0$, for some $w$ and $j$.

\hspace{.4in} \space \space$\Leftarrow$ If $S \in D_{w,j}(M)$ with $\chi (S) >0$, then one of the connected components of $S$, say $S'$, must be homeomorphic to either $S^{2}$, $\mathbf{RP}^{2}$, or $D^{2}$. The first two cases imply the existence of a non-trivial normal 2-sphere which contradicts our assumption on the triangulation of $M$. Hence, $S' \cong D^{2}$ and it is non-trivial by definition of $D_{w,j}(M)$.

\hspace{1cm} Define $C_{w, i}(M)$ to be the solution space to the surface equations with $t_{i}=0$, and the rest of the variables being \emph{real} non-negative.

\bigskip
\begin{claim} ~\label{C6.10}
$\exists$ a non-trivial disk D in M \space $\Leftrightarrow$ \space $\exists V \in C_{w, j}(M)$, for some $w, j$, with $\chi (V) >0$, V is on an external ray of this cone, and $V \in \mathbf{Z}^{7t}$.
\end{claim}

\underline{\bf{Proof :}} $\Rightarrow$ If $D$ is non-trivial, then $D \in C_{w, j}(M)$ by claim~\ref{C6.9}. Suppose $D$ is not on an external ray, i.e. $D$ = $V_{1} + $ ...$+ V_{p}$ for some $p$, where $V_{i} \in C_{w, i}(M)$ and all the $V_{i}$'s are on an external ray. Because $D$ (viewed as the solution of the normal surface equations representing $D$) has integer entries, all the $V_{i}$'s must have rational entries. Since $D$ is homeomorphic to a disk, we have $\chi(D) = \chi (V_{1} + ...+ V_{p}) = \chi (V_{1}) + ... + \chi (V_{p}) > 0$. This implies that  $\chi (V_{i}) > 0$ for some $i$. Multiply each entry of $V_{i}$ by the least common multiple of the denominators of all its entries. The new surface, $V'_{i}$, may be a sphere, a projective plane, or a disk. The existence of a non-trivial 2-sphere or a  projective plane contradict our assumption on the triangulation of $M$. Also, because $V'_{i}$ is on an external ray it cannot be disconnected, and because $V'_{i} \in C_{w, i}(M)$, it cannot be trivial. Hence, $V'_{i}$ is a normal surface with the desired properties.

\hspace{.4in} \space \space $\Leftarrow$ There is nothing to prove.

\bigskip

\hspace{1cm} \underline{\bf{The algorithm :}} Let $M$ be an irreducible compact orientable 3-manifold with nonempty boundary.

\hspace{1cm} {\bf Step 1:} We first apply Casson's algorithm to look for non-trivial 2-spheres, to cut along them, and to collapse them. Since M is irreducible, there is no need to search for almost normal 2-spheres.

\hspace{1cm} {\bf Step 2:} After collapsing those 2-spheres, $M$ may contain vertices in its interior and it may also contain more than one vertex for each connected boundary component. 

\hspace{1cm} Case 1: Suppose $M$ has a vertex $v$ in its interior. Let $e$ be any edge joining $v$ and $v'$, where $v'$ is any vertex in $\partial M$ (we can find such an edge since $M$ is connected). Because $e$ is embedded, the boundary of its neighborhood is homeomorphic to a disk $D$. We normalize $D$ and obtain a certain number of normal disks and 2-spheres $D_{1}$, ..., $D_{m}$, $S_{1}$, ..., $S_{l}$. Because $M$ is irreducible, each of the 2-spheres bounds a ball which does not contain $e$. Indeed, if it did, then $M$ would be homeomorphic to $S^{3}$ which is a contradiction since $\partial(M) \neq \emptyset$. Hence, we ignore the 2-spheres (or we can simply assume that $M$ does not contain any by applying the collapsing process of Theorem~\ref{MT} to all non-trivial 2-spheres). We now cut along the disks and use Theorem~\ref{MT2} to collapse each of them to points. Note that when more than one disk is collapsed, there could be a $\partial$-annulus glued to two different disks, but by Lemma~\ref{L5}, this does not affect the collapsing process. So after collapsing all the disks, we obtain a decomposition of $M$ into, say $n$ pieces. Consider one of the summands. If it has a 2-sphere as boundary, it must be homeomorphic to a 3-ball by irreducibility and hence we discard it and we look at another piece in the decomposition of $M$. If a piece still has a vertex in its interior, we reiterate the above procedure (note that eventually no vertices in the interior will be found since collapsing disks strictly reduces the number of tetrahedra). If it has no interior vertices, we then go to case 2.

\hspace{1cm} Case 2: Suppose $M$ has more than one vertex in one of its connected boundary component. Let $e$ be any edge, on the boundary of M, joining two disjoint vertices. As in case 1 we take the boundary of a regular neighborhood of $e$, normalize it, and cut along the resulting non-trivial normal disks and spheres. We collapse them using Theorems~\ref{MT} and ~\ref{MT2}. Consider one summand in the resulting decomposition of $M$. If it has a 2-sphere as one of its connected boundary component, we then discard it and we look at another summand. If it has more than one vertex on one of its connected boundary component or if it contains a vertex in its interior, we then reiterate case 1 or case 2 accordingly.

\hspace{1cm} By Remark 1, we know that normalizing the boundary of a regular neighborhood of an edge takes time linear in the number of tetrahedra. Moreover, by Theorems~\ref{MT} and ~\ref{MT2}, we know that each cut-and-collapse operation strictly reduces the number of tetrahedra and takes time linear in the number of tetrahedra. Hence, we eventually obtain, in linear time, a decomposition of $M$ into irreducible pieces (which may not necessarily be $\partial$-irreducible) with a triangulation as in the conclusion of Proposition~\ref{Pro4}. We can now run step 3 for each of the summands obtained.

\hspace{1cm} {\bf Step 3:} We consider a summand of $M$ which we will still call $M$. For each cone $C_{w, i}(M)$, we maximize $\chi$ on the convex polyhedra $C_{w, i}(M) \bigcap$ $\{ \sum_{i=1}^{7t} t_{i} = 1 \}$, and we obtain a vertex solution. By Claim~\ref{C5} above, if it represents a normal surface $S$ satisfying $\chi(S) > 0$, then we know there exists a non-trivial normal disk. We cut along the disk and collapse it. We call $M$ one of the resulting pieces and go to step 2 unless the triangulation satisfies the conclusion of Proposition~\ref{Pro4}. If the vertex solution found satisfies $\chi(S) \leq 0$, we know, by Claim~\ref{C5}, that there are no non-trivial normal disks in this cone. We repeat this step with a new cone. 

\hspace{1cm} {\bf Step 4:} For each summand found, we repeat steps 2 and 3 until no non-trivial normal disks are found. We obtain a decomposition of $M$, and we use Rubinstein's algorithm to check which summand is homeomorphic to a ball.

\bigskip

\hspace{1cm} \underline{\bf{Complexity:}} Let $M$ be a compact irreducible orientable 3-manifold with nonempty boundary. We first run Casson's algorithm to find and collapse all non-trivial normal 2-spheres. We showed in the previous section that this procedure takes time exponential in the number of tetrahedra ($O(3^t) \cdot t^n$). If $M$ contains vertices in its interior or$\backslash$and more than one vertex on a connected boundary component, we showed in step 2 above that it takes linear time to obtain a decomposition of $M$, where each piece has a triangulation as in the conclusion of Proposition~\ref{Pro4}. 

\hspace{1cm} For each summand of $M$ (there could be at most $t$ of them), there are at most $4t \cdot 3^{t}$ different cones, and it takes polynomial time (see previous section) to maximize $\chi$ on the corresponding convex polyhedra. If a non-trivial disk is found, we collapse it. It was shown in Chapter 5 that it takes linear time ($4t$) to do this collapsing. Some of the resulting pieces of the collapsing may not have triangulations satisfying the conclusion of Proposition~\ref{Pro4}. If this is the case, we run step 2 again which takes. For each summand obtained at the end of the procedure, we run Rubinstein's algorithm to check which pieces are homeomorphic to balls. If $M$ does not contain any non-trivial 2-spheres, we conclude that the running time for the above procedure is $L \cdot O(t^{r}) \cdot 3^{t}$, where $r$ and $L$ are independent of $t$.


\bigskip

\medskip

\hspace{1cm} As a direct application of this algorithm, we take a quick look at the unknotting problem. At the beginning of the century, the mathematician Max Dehn found a relation between the triviality of a knot embedded in $\mathbf{R}^3$ and the fundamental group of its complement. It was only in 1957 that C.D.Papakyriakopoulos proved Dehn's result which states that a knot is trivial if and only if the fundamental group of the knot is isomorphic to $\mathbf{Z}$. With the work of H. Seifert in the 1930's, we know today that a knot is trivial if and only if there exists an embedded non-separating disk in its complement.

\hspace{1cm} Moreover, given a polygonal knot in $\mathbf{R}^{3}$, Hass, Pippenger, and Lagarias (~\cite{HLP:gnus}) gave a canonical construction to triangulate the knot complement and find such a disk. So all we have to do is the following: each time we will find a non-trivial normal disk, we will check if it separates the manifold or not. If all of the non-trivial normal disks found separate the manifold, then we know we have a non-trivial knot. If one of them is non-separating , then we have a representation of the trivial knot.



\newpage
\pagestyle{myheadings} 
\markright{  \rm \normalsize BIBLIOGRAPHY. \hspace{0.5cm}
  Minimal Triangulations}


\begin{thebibliography}{99}
\thispagestyle{myheadings}
\addcontentsline{toc}{chapter}{\bf Bibliography}


\bibitem[Ar]{Ar:gnus} Michael Artin, {\it Algebra}.

\bibitem[Ba]{Ba:gnus} David Bachman, {\it A Note on Kneser-Haken Finiteness}, to appear.

\bibitem[Ca]{Ca:gnus} Andrew Casson, {\it Notes on three-dimensional topology}.

\bibitem[Hak] {Hak:gnus}  Wolfgang Haken, {\it Theorie der Normalflachen}. Acta Math., 1961, 105, 245-375.

\bibitem[Ha] {Ha:gnus} Joel Hass, {\it Algorithms for Recognizing Knots and 3-Manifolds}. Knot theory and its applications. Chaos Solitons Fractals 9 (1998), no. 4-5, 569--581.

\bibitem[HLP] {HLP:gnus}  Joel Hass, Jeffrey Lagarias and Nicholas Pippenger, {\it The computational Complexity of Knot and Link Problems}. J. ACM 46 (1999), no. 2, 185--211.

\bibitem[Hat 1] {Hat 1:gnus} Allen E. Hatcher, {\it Notes on Basic 3-Manifold Topology}.

\bibitem[Hat 2] {Hat 2:gnus} Allen E. Hatcher, {\it Algebraic Topology}.

\bibitem[He] {He:gnus} Geoffrey Hemion, {\it On the classification of homeomorphisms of 2-manifolds and the classification of 3-manifolds}. $ Acta$ $Math. 142 (1979), no. 1-2, 123--155$.

\bibitem[J] {J:gnus} Klaus Johannson, {\it Homotopy equivalences of 3-manifolds with boundaries}, volume 761 of Lecture Notes in Mathematics. Springer, Berlin, 1979.

\bibitem[JLR] {JLR:gnus} William Jaco, David Letscher, J. Hyam Rubinstein, {\it Algorithms for Essential Surfaces in 3-Manifolds}. In preparation.

\bibitem[JO] {JO:gnus} William Jaco and Ulrich Oertel, {\it An algorithm to decide if a 3-manifold is a Haken manifold}. Topology 23 (1984), no. 2, 195--209.

\bibitem[JR] {JR:gnus}  William Jaco and J. Hyam Rubinstein, {\it 0-Efficient Triangulations of 3-Manifolds}. Preprint, 2001.

\bibitem[JS] {JS:gnus} William Jaco and Eric Sedgwick, {\it Decision Problems in the Space of Dehn Filings}. Topology, 2001.

\bibitem[JSh] {JSh:gnus} William Jaco and Peter B. Shalen, {\it Seifert fibered spaces in 3-manifolds}, Mem. Amer. Math. Soc., 1979, 220.

\bibitem[JT] {JT:gnus} William Jaco and Jeffrey L. Tollefson, {\it Algorithms for the Complete Decomposition of a Closed 3-Manifold}.  Illinois J. Math. 39 (1995), no. 3, 358--406.

\bibitem[Kn] {Kn:gnus} Kneser, H., {\it Geschlossene Flachen in dreidimensionalen Mannifgfaltigkeiten}, Jber. Ent. Math., Ver., 1929, 28, 248-260.

\bibitem[La] {La:gnus} Marc Lackenby, {\it Taut ideal triangulations of 3-manifolds}.  $Geom.$ $Topol.$ 4 (2000), 369--395. 

\bibitem[Ma] {Ma:gnus} William S. Massey, {\it A Basic Course in Algebraic Topology}.

\bibitem[MR] {MR:gnus} Saburo Matsumoto and Richard Rannard , {\it The Regular Projective Solution Space of the Figure-Eight Complement}. Experiment. Math. 9 (2000), no. 2, 221--234.

\bibitem[Mat] {Mat:gnus} S. V. Matveev, {\it Classification of sufficiently large three-dimensional manifolds}.  (Russian) Uspekhi Mat. Nauk 52 (1997), no. 5(317), 147--174; translation in Russian Math. Surveys 52 (1997), no. 5, 1029--1055.

\bibitem[Mi] {Mi:gnus} J. Milnor, {\it A unique factorization theorem for 3-manifolds}. Amer. J. Math. 84(1962), 1-7.

\bibitem[Mo] {Mo:gnus} Edwin E. Moise, {\it Affine Structures in 3-Manifolds: V. The triangulation Theorem and Hauptvermutung}. Ann. of Math. (2) 56, (1952). 96--114.

\bibitem[Pa] {Pa:gnus} Udo Pachner, {\it Konstruktionsmethoden und das kombinatorische Homomorphieproblem fur Triangulationen kompakter semilinearer Mannigfaltigkeiten}. (German) [Construction methods and the combinatorial homeomorphism problem for triangulations of compact semilinear manifolds] Abh. Math. Sem. Univ. Hamburg 57 (1987), 69--86.

\bibitem[Ra] {Ra:gnus} Richard Rannard, {\it Computing Immersed Normal Surfaces in the Figure-Eight Knot Complement}. Experiment. Math. 8 (1999), no. 1, 73--84.

\bibitem[Ru] {Ru:gnus} J. H. Rubinstein, {\it An algorithm to recognize the 3-sphere}. 

\bibitem[Sch] {Sch:gnus} Alexander Schrijver, {\it Theory of Linear and Integer Programming}.

\bibitem[Sc] {Sc:gnus} Peter Scott, {\it The Geometries of 3-Manifolds}.  Bull. London Math. Soc. 15 (1983), no. 5, 401--487.

\bibitem[Sh] {Sh:gnus} Peter B. Shalen, {\it Representations of 3-manifold groups}. Handbook of Geometric Topology, to appear.

\bibitem[Sp] {Sp:gnus} Edwin Henry Spanier, {\it Algebraic Topology}.

\bibitem[Tho] {Tho:gnus} Abigail Thompson, {\it Thin Position and the Recognition Problem for $S^{3}$}. Math. Res. Lett., 1(5):613-630, 1994.

\bibitem[Th 1] {Th 1:gnus} William P. Thurston, {\it The geometry and topology of 3-manifolds}. 

\bibitem[Th 2] {Th 2:gnus} William P. Thurston, {\it Three-Dimensional Geometry and Topology}. Edited by Silvio Levy. Princeton, N.J. : Princeton University Press, 1997- v. <1 > : ill. ; 24 cm.







\end{thebibliography}
\end{document}